%% file: slops_nocomment.tex
\documentclass[12pt,reqno]{amsart}

\def\ifundefined#1{\expandafter\ifx\csname#1\endcsname\relax} \ifundefined{preambleloaded}
    \input{preamble_v3}
\else
\fi

\DeclareMathOperator{\xxMA}{MA}
\DeclareRobustCommand*{\MA}[1]{\ifmmode\xxMA#1\else{Monge--Ampère}\fi}
\DeclareMathOperator{\xxDH}{DH}
\renewcommand{\DH}{\xxDH}
\newcommand{\ddat}[2]{\left.\frac{d}{d #1}\right|_{#2}}

\DeclareMathOperator{\Fut}{Fut}
\DeclareMathOperator{\ini}{in}

\DeclareMathOperator{\wt}{wt}
\DeclareMathOperator{\Hilb}{Hilb}
\newcommand{\pr}{\mathrm{pr}}

\newcommand{\PSH}{\mathrm{PSh}}

\newcommand{\amb}{{\mathbb{P}^{N_1} \times \mathbb{C}^{N_2}}}
\newcommand{\ambC}{{\mathbb{P}^{N_1} \times \mathbb{C}^{N_2} \times \mathbb{C}}}
\newcommand{\bp}{\bar{\partial}}

\usepackage[normalem]{ulem}

\renewcommand{\bar}{\widebar}
\renewcommand{\tilde}{\widetilde}
\renewcommand{\hat}{\widehat}
\usepackage{amsmath}
\usepackage{xcolor,etoolbox}
\usepackage{listings}

\usepackage{geometry}
\usepackage{microtype}

\AtBeginDocument{
    \setlength{\abovedisplayshortskip}{-4pt plus 5pt minus 0pt}
}

\numberwithin{equation}{section}

\renewcommand{\demph}[1]{\textcolor{violet!75!blue}{\textbf{#1}}}

\title[K-polystability of AC Kähler--Ricci shrinkers]{K-polystability of Asymptotically Conical Kähler--Ricci Shrinkers}
\date{\today}

\author[C. Cifarelli]{Charles Cifarelli}  
\address{\noindent Charles Cifarelli \newline  \indent Mathematics Department, Stony Brook University, Stony Brook, NY 11794, USA}
\email{charles.cifarelli@stonybrook.edu}

\author[C. Esparza]{Carlos Esparza}  
\address{\noindent Carlos Esparza \newline
\indent Department of Mathematics,
UC Berkeley,
970 Evans Hall,
Berkeley, CA 94720, USA
}
\email{esparzac@berkeley.edu}

\begin{document}

\begin{abstract}
    Recently, Sun--Zhang \cite{SunZhang} have developed an algebraic theory for Kähler--Ricci shrinkers showing that they admit the structure of a polarized Fano fibration $(\pi:X \to Y, \, \xi)$. In particular, they conjecture \cite[Conjecture~6.1]{SunZhang} that existence of a Kähler--Ricci shrinker metric is equivalent to a notion of K-stability. We prove one direction of this conjecture, namely that existence of a Kähler--Ricci shrinker metric $g$ implies K-polystability of $(\pi:X \to Y, \, \xi)$, in the case that the Ricci curvature of $g$ decays at infinity.
    As an application, we give a non-existence result: if $M$ is the blowup of a six-dimensional quadric along a two-dimensional subquadric, then the total space $X$ of the cube root of $K_M$ is a polarized Fano fibration not admitting a Kähler--Ricci shrinker.
\end{abstract}
\maketitle

\section{Introduction}

By definition, a \demph{Kähler--Ricci shrinker} (or complete shrinking gradient Kähler--Ricci soliton) is a triple $(X, \, \omega, \, \xi)$,  %
where $(X, \, \omega)$ is a complete Kähler manifold and $\xi = J\nabla f$ is a hamiltonian real holomorphic Killing vector field, such that
\begin{equation}\label{eq:shrinker}
    \Ric_\omega + i\partial \bp f = \omega
    \p
\end{equation}
In recent work by Sun and Zhang \cite{SunZhang}, it was shown that any Kähler--Ricci shrinker $(X, \, \omega, \, \xi)$ is a quasiprojective algebraic variety, and specifically it naturally admits the structure of a \emph{polarized Fano fibration}; see  \cref{thm:KRStoPFF} below. This means in particular that $X$ is a Fano fibration, so that there is an affine variety $Y$ and a fibration map $\pi:X \to Y$ such that $-K_X$ is $\pi$-ample; see \cref{def:pff} for details. In their paper, Sun and Zhang formulated a Yau--Tian--Donaldson type conjecture for polarized Fano fibrations:

\begin{conj}[{\cite[Conjecture 1.2]{SunZhang}}]\label{conj:pffytd}
    A polarized Fano fibration $(\pi: X \to Y, \, \xi)$ admits a Kähler--Ricci shrinker $(\omega, \, \xi)$, unique up to automorphisms of $X$ preserving $\xi$, if and only if it is K-polystable.
\end{conj}

The main goal of this paper is to prove one direction of the correspondence in the \emph{asymptotically conical} case:

\begin{thm}\label{mtheorem}
    Let $(X, \, \omega, \, \xi)$ be a Kähler--Ricci shrinker whose Ricci curvature tends to zero at infinity. Then the associated polarized Fano fibration $(\pi:X \to Y, \, \xi)$ is K-polystable. 
\end{thm}

The condition that $|{\Ric}|_{\omega} \to 0$ at infinity is actually equivalent to $\omega$ being asymptotically conical (see \cref{def:AC}). Indeed, it was proved in \cite{MWconical} that a Ricci shrinker whose Ricci curvature decays at infinity must have quadratic decay of the full curvature. Then in \cite{CDS} it was shown that a Kähler--Ricci shrinker $(X, \, \omega, \, \xi)$ has quadratic curvature decay if and only if it is asymptotically conical. In fact, it was already shown in that paper that such an $(X, \, \omega, \, \xi)$ admits a polarized Fano fibration structure, by showing that $X$ admits a map $\pi: X \to Y$ to its tangent cone at infinity $(Y, \, \omega_Y)$ which is a resolution of singularities (see \cref{thm:quadraticdecaytoAC} below). Given a fixed complex manifold $X$, uniqueness up to automorphisms of Kähler--Ricci shrinkers $(\omega, \, \xi)$ on $X$ with quadratic curvature decay was proved in \cite{CarlosUniq}.

The proof of \cref{mtheorem} follows a combination of the approaches of Berman \cite{BermanQFano} and Collins--Székelyhidi \cite{ColSze2}. In particular, using ideas of \cite{ColSze2}, we show that algebraic and analytic Futaki invariants of a test configuration $\cX$ are both equal to the slope at infinity of the Ding functional along a geodesic ray, which is sufficient to show that existence implies semistability. For polystability, we give an adaptation of the argument in \cite{BermanQFano}. The primary new difficulty is the existence of geodesic rays in this non-compact setting. To solve this, we first construct a good background metric associated to $\cX$, by showing it can be equivariantly embedded in an ambient space $\amb$, which we think of as a barrier function. Using this and the existence of geodesic segments established in \cite{CarlosUniq}, we are able to obtain a geodesic ray by taking a limit.

It is in fact only this last point where the quadratic curvature decay condition is used, and most of the results in this paper hold for general Kähler--Ricci shrinkers. In the general case, the expectation is that a Kähler--Ricci shrinker $(\omega, \, \xi)$ on a polarized Fano fibration $(\pi:X \to Y, \, \xi)$ is \emph{weakly asymptotically conical}, meaning that it has a unique tangent cone at infinity which is a Kähler cone metric on $Y$ (\cite[Conjecture 6.2]{SunZhang}). In particular, we expect that the same strategy outlined here will work in the general case once a suitable existence theory for geodesic \emph{segments} is established in the weakly asymptotically conical setting.

To date, the vast majority of known noncompact Kähler--Ricci shrinkers are in fact asymptotically conical. These include the well-known examples of Feldman-Ilmanen-Knopf \cite{FIK} and many similar constructions living on the total space of vector bundles over a compact Fano manifold \cite{ChiLi, CharlieExamples, DancerWang, FutakiWang, Futaki_Eguchi, Yang_Examples}, see \cite{ZhangLi} for a nice summary. To the knowledge of the authors, the only known example which is not AC or a product is the one constructed on $\Bl_p(\bbP^1 \times \bbC)$ in \cite{BCCD}.

We use \cref{mtheorem} to give a concrete example of non-existence. Specifically, we show that a specific Fano fibration $\pi: X \to Y$ is K-unstable for any polarization $\xi$, and therefore admits no AC Kähler--Ricci shrinker.

\begin{cor}\label{mcor}
    Let $Q^6$ be a smooth quadric hypersurface in $\bbP^7$, and let $M$ be the blowup of $Q^6$ in a two-dimensional linear subquadric. There is a line bundle $L$ on $M$ with $3L = K_M$, and we let $X$ be the total space of $L$. Then $X$ does not admit an asymptotically conical Kähler--Ricci shrinker metric.
\end{cor}

By recent work \cite{CD_AC, ZhangLi}, any Kähler--Ricci shrinker on $X$ is necessarily AC, thus it follows that $X$ does not admit any Kähler--Ricci shrinker whatsoever. The Fano manifold $M$ is part of a broader class which were shown by Delcroix not to admit (compact) Kähler--Ricci shrinker metrics themselves \cite{Del_Examples}. It is likely that the same can be shown for the total space of any root of the canonical bundle in this class, but we present the specific example above for simplicity. For more details and precise statements, see \cref{sec:examples}. 

The paper is organized as follows. In \cref{sec:background} we recall some basics on polarized Fano fibrations and K-stability, PSh metrics on line bundles, and asymptotically conical Kähler metrics that we will use throughout.
\Cref{sec:testconf} is where we prove that a test configuration for a polarized Fano fibration $(\pi:X \to Y, \, \xi)$ can be equivariantly embedded in $\amb$ and also establish some basic algebraic results about such embedded families.
In \Cref{sub:product_testconf}, we also include a proof in the polarized Fano fibration setting of a well-known result whose statement we could not find in the literature: a test configuration whose central fiber is isomorphic to its generic fiber must be a product.
\Cref{sec:geodesics} is dedicated to the construction of geodesic rays, using a suitable background metric on $\amb$ as a barrier function.
In \cref{sec:energyfunctionals} we show that the slope at infinity along our geodesic ray is given by the analytic Futaki invariant, and then use an idea of \cite{ColSze2} to show that this coincides with the algebraic Futaki invariant by doing a further degeneration to a union of hyperplanes with multiplicity.
Finally, in \cref{sec:maintheorem} we use the method of \cite{BermanQFano} to complete the proof of \cref{mtheorem}.

\subsection{Acknowledgements}
Both authors would like to thank Robert Berman, Thibaut Delcroix, Connor Halleck-Dub\'{e}, Yuji Odaka, Jason Starr, Song Sun, Dror Varolin, and Junsheng Zhang for comments and helpful discussions. This work was partially completed while the first author was in residence at SLMath (NSF Grant
DMS-1928930), and both are thankful for the ideal working environment provided there. We are also thankful for the IASM at Zhejiang University and the IGP at USTC for their hospitality during the Summer of 2025, when this work was also partially completed. The first author is supported by NSF Grant DMS-2506521.

\section{Background}
\label{sec:background}

\subsection{Polarized Fano fibrations and K-stability}

Here we recall the polarized Fano fibration framework of \cite{SunZhang}, which generalizes the notion of a polarized affine cone originally defined in \cite{ColSzeStab, ColSze2}.

\begin{defn}[Polarized affine cone]\label{def:polariedaffinecone}
    A \demph{polarized affine cone} $(Y, \, \xi, \, T)$ is a normal affine variety $Y = \Spec(A)$ with a torus action $T$ fixing a unique point $o \in Y$ and equipped with a Reeb field $\xi \in \mathfrak{t}= {\rm Lie(T)}$. A \demph{Reeb field} $\xi$ is an element of $\mathfrak{t}$ such that under the weight decomposition
    \[ A = \bigoplus_{\alpha \in \mathfrak{t}^*} A_{\alpha},\]
    we have $\langle \alpha, \, \xi \rangle > 0$. %
\end{defn}

\begin{defn}[Polarized Fano fibration]\label{def:pff}
    A \demph{polarized Fano fibration} $(\pi:X \to Y, \, \xi)$ is a fibration $\pi:X \to Y$ with the following properties:
    \begin{enumerate}
        \item $\pi:X \to Y$ is a Fano fibration. That is, $X$ and $Y$ are normal varieties, $X$ is klt, and $-K_X$ is $\pi$-ample and $\bbQ$-Cartier.
        \item $X$ and $Y$ are equipped with a $\pi$-equivariant torus action $T$, and $\xi \in \mathfrak{t} = \Lie(T)$.
        \item $(Y, \, T, \, \xi)$ is a polarized affine cone.
    \end{enumerate}
    Moreover, given a polarized Fano fibration $(\pi:X \to Y, \, \xi)$, any $\eta \in \ft$ is called a \demph{Reeb field} if it is a Reeb field on $Y$, in other words if $(Y, \, T, \, \eta)$ is a polarized affine cone.
\end{defn}

\begin{ex} \label{ex:amb}
    The simplest example of a polarized Fano fibration is given by $\amb$ for $N_1, \, N_2 \geq 0$, with a choice of suitable vector field $\xi$. To describe this, first consider the extreme cases where $N_1 = 0$ or $N_2 = 0$. If $N_2 = 0$, then $\amb = \bbP^{N_1}$ is Fano, and any $\xi \in \mathfrak{pgl}(N_1+1)$ makes $(\pi:\bbP^{N_1} \to \{0\}, \, \xi)$ a polarized Fano fibration. If instead $N_1 = 0$, then $\amb = \bbC^{N_2}$ is affine, and a choice of Reeb field $\xi$ which makes $(\pi: \bbC^{N_2} \to \bbC^{N_2}, \, \xi)$ a polarized Fano fibration is exactly as in \cref{def:polariedaffinecone}.

    In general, we say that a linear torus action $T \hookrightarrow \PGL(N_1 + 1) \times \GL(N_2)$ on $\amb$ is \demph{diagonal} if it lies inside the subgroup $T \subset \bbT^{N_1} \times \bbT^{N_2}$ of diagonal matrices, and \demph{compatible} if the associated $T$-action on $\bbC^{N_2}$ is effective and admits at least one Reeb field. Given a compatible $T$-action, any $\xi \in \ft$ whose image in $\Lie(\bbT^{N_2}) \subset \mathfrak{gl}_{N_2}$ is a Reeb field on $\bbC^{N_2}$ will make $(\pi: \amb \to \bbC^{N_2}, \, \xi)$ into a polarized Fano fibration.
\end{ex}

\begin{rmk} \label{rmk:pffEmbed}
    If $(\pi: X \to Y, \, \xi)$ is a polarized Fano fibration, we can use the fact that $-K_X$ is $\pi$-ample to embed it into a product of projective and affine spaces. More precisely, there exists a linear $T$-action on $\amb$, a positive integer $p$, and embeddings $\iota_X: (X, -p K_X) \inj (\amb, \cO(1))$ and $\iota_Y: Y \hookrightarrow \bbC^{N_2}$
    such that the following diagram commutes:
    \begin{cd*}
        X \arrow[d, "\pi"] \arrow[r, "\iota_X"] & \amb \arrow[d, "\pr_2"] \\
        Y \arrow[r, "\iota_Y"] & \bbC^{N_2}
    \end{cd*}
\end{rmk}

The reason to introduce polarized Fano fibrations is that they are the natural algebraic setting in which to study Kähler--Ricci shrinkers. Indeed, by the first main result in \cite{SunZhang}, the complex manifold underlying any Kähler--Ricci shrinker admits a polarized Fano fibration structure:

\begin{thm}[{\cite[Theorem 1.1]{SunZhang}}]\label{thm:KRStoPFF}
    A Kähler--Ricci shrinker $(X, \, g, \, \xi)$ determines a polarized Fano fibration structure $(\pi: X \to Y, \, \xi)$ for $(X, \, \xi)$.
\end{thm}

This result provides the underpinning for the algebraic side of the study of Kähler--Ricci shrinkers. In particular, it follows that if $\pi: X \to Y$ is a polarized Fano fibration with $Y = \Spec(A)$, then the ring
\begin{equation}\label{eq:PFFRing1}
    R = \bigoplus_{m\geq 0} R_m, \qquad R_m =H^0(X, \, -mK_X),
\end{equation}
is a finitely-generated $A$-algebra, and $X = \Proj_{A}R$. Note that in general $X$ is only $\bbQ$-Cartier, and in this case $R_m \neq 0$ only when $m$ is sufficiently divisible. Since we have a torus action $T$ on $\pi:X \to Y$, we can further decompose $R$ into weight spaces, namely
\begin{equation}\label{eq:PFFRing2}
    R_m = \bigoplus_{\alpha \in \mathfrak{t}^*} R_{m,\alpha},
\end{equation}
where $T$ acts on each $R_{m,\alpha}$ with weight $\alpha$. Given this, we can define the \emph{weighted volume} of a vector field $\xi \in \ft$, as long as it is a Reeb field:

\begin{defn}[{\cite[Lem.~5.7, Prop.~5.9]{SunZhang}}] Given a Reeb field $\xi \in \ft$, the expression
    \begin{equation}
        \label{eq:volume}
        \bbW_X(\xi) =  \lim_{m \to \infty} \frac{1}{m^n}\sum_{\alpha \in \ft^*} e^{- \langle \frac{\alpha}{m}, \, \xi \rangle } \dim R_{m, \alpha}
    \end{equation}
    converges, and is called the \demph{weighted volume} of $\xi$.
\end{defn}

The algebraic weighted volume was defined originally in \cite{SunZhang} in terms of valuations on $X$, and subsequently they showed that the characterization \cref{eq:volume} is equivalent. This expression is more convenient for our purposes, so we adopt it as the definition directly.

The reason for the name weighted volume comes from considering the case where $X$ is smooth. Here one can define the \demph{analytic weighted volume} by
\begin{equation}\label{eq:analyticweightedvolume}
    \bbW^\mathrm{an}_X(\xi) = \frac{1}{(2\pi)^n}\int_X e^{-\langle \xi,  \, \mu \rangle } \frac{\omega^n}{n!},
\end{equation}
where $\omega \in c_1(-K_X)$ is a suitable choice of Kähler metric and $\mu: X \to \ft^*$ is a suitably normalized moment map (see \cref{section:hamiltonain} for details). This was introduced in \cite{CDS} as an extension of the well-known notions of weighted volume in the compact \cite{TianZhu02} and affine \cite{MSYsasaki, ColSzeStab} cases. It was proved in \cite{CDS} that \cref{eq:analyticweightedvolume} is well-defined for an open convex cone $\xi \in \Lambda \subset \ft$, and in \cite[Prop.~5.9]{SunZhang} that for any $\xi \in \Lambda$, we have
\[ \bbW^{\mathrm{an}}_X(\xi) = \bbW_X(\xi).\]
One of the key steps below will be to establish a version of this equality even when $X$ is singular; see \cref{thm:volumeAgrees}.

In order to define K-stability for a polarized Fano fibration, we introduce the notion of a special test configuration following \cite{SunZhang}:

\begin{defn}[{\cite[Definition 5.2]{SunZhang}}]\label{def:testconfiguration}
    Let $(\pi:X \to Y, \, \xi)$ be a polarized Fano fibration. A ($T$-equivariant) \demph{special test configuration} $(\Pi: \cX \to \cY,\, T,\,  \eta)$ is a commutative diagram
    \begin{equation}\label{def:testconfiguration1}
    \begin{tikzcd}[column sep = .3in]
        & \cX \arrow[ld, "\Pi"'] \arrow[dr, "\Pi_{\cX}"] & \\
        \cY \arrow[rr, "\Pi_{\cY}"] & & \bbC
    \end{tikzcd}
    \end{equation}
    where $\Pi_{\cX}:\cX \to \bbC$, $\Pi_{\cY}: \cY \to \bbC$ are surjective flat morphisms, satisfying the following properties:
    \begin{enumerate}
        \item $\cX$ is $\bbQ$-Gorenstein.
        \item There is a $\bbC^*$-action on $\cX$ generated by a holomorphic vector field $\eta$ such that $\Pi_*\eta$ also generates a $\bbC^*$-action on $\cY$. Moreover, away from $0 \in \bbC$, there is a $\bbC^*$-equivariant isomorphism between the diagram \cref{def:testconfiguration1} and
        \begin{equation}\label{def:testconfiguration2}
    \begin{tikzcd}[column sep = .075in]
         & X \times \bbC^* \arrow[ld] \arrow[dr] &  \\
         Y \times \bbC^*  \arrow[rr] & & \bbC^*
        \p
    \end{tikzcd}
    \end{equation}
    \item $T$ is a torus action on $\cX$ and $\cY$, equivariant with respect to $\Pi$, commuting with $\eta$, and inducing the trivial action on $\bbC$ via $\Pi_{\cX}$ (i.e. acting only on the fibers of $\Pi_{\cX}$). Moreover, with respect to the identification \cref{def:testconfiguration2} away from $0 \in \bbC$, we have that $\xi \in \mathfrak{t} = \Lie(T)$. We set $T' = T \times \bbC^*$ to be the action generated by $T$ and $\eta$ on $\cX$, so that $\Pi: \cX \to \cY$ is also $T'$-equivariant.
    \item If we set $X_0 \coloneqq \Pi_{\cX}^{-1}(0),\, \cY_0 \coloneqq \Pi_{\cY}^{-1}(0)$, then $(\Pi: X_0 \to \cY_0, \, \xi)$ is a polarized Fano fibration. Note that $\xi \in \mathfrak{t}$ gives rise to a vector field on $X_0$ through the $T$-action. In particular this means that the (scheme-theoretic) central fiber is reduced.
    \end{enumerate}
\end{defn}

\begin{defn}[{\cite[Definition 5.4]{SunZhang}}]\label{def:futaki}
    Given any polarized Fano fibration $(\pi:X \to Y, \, \xi)$ and any $\eta \in \ft$, we can define the \demph{Futaki invariant}
    \begin{equation}\label{eq:futaki1}
        \Fut_\xi(X, \, \eta) = - \left. \frac{d}{dt}\right|_{t=0} \bbW_X\left( \xi + t \eta\right).
    \end{equation}
    Moreover, given a special test configuration $(\Pi:\cX \to \cY, \, T, \, \eta)$, we define the Futaki invariant to be $\Fut_\xi(X_0, \, \eta)$, where $(\Pi_0: X_0 \to Y_0, \, \xi)$ is the central fiber. When the context is understood, we will often simply write $\Fut_{\xi}(\eta)$.
\end{defn}
\begin{lem}[{\cite[Lemma 5.7]{SunZhang}}]
    Given any polarized Fano fibration $(\pi:X \to Y, \, \xi)$ and $\eta \in \ft$, we have that the right-hand side in \cref{eq:futaki1} always exists, and moreover
    \begin{equation}\label{eq:futaki2}
        \Fut_\xi(X, \, \eta) = \lim_{m \to \infty} \frac{1}{m^n}\sum_{\alpha \in \ft^*} \left\langle \frac{\alpha}{m}, \, \eta \right\rangle \, e^{- \langle \frac{\alpha}{m}, \, \xi \rangle } \dim R_{m, \alpha}.
    \end{equation}
\end{lem}
\begin{defn}[{\cite[Definition 5.5]{SunZhang}}]\label{def:Kstability}
    A polarized Fano fibration $(\pi:X \to Y, \, \xi)$ is said to be \demph{K-semistable} if, for any special test configuration $(\Pi:\cX \to \cY, \, T, \, \eta)$, we have that $\Fut_\xi(\eta) \geq 0$. It is said to be \demph{K-polystable} if it is K-semistable, and moreover any special test configuration $(\Pi:\cX \to \cY, \, T, \, \eta)$ such that $\Fut_\xi(\eta) = 0$ is $T$-equivariantly isomorphic to $(X \times \bbC, \, T, \, \eta)$, where $T$ acts trivially on $\bbC$ and $\eta = \eta_0 + \tau\partial_\tau$ for a holomorphic vector field $\eta_0$ on $X$ commuting with $\ft$ and $\tau\partial_\tau$ the standard Euler field on $\bbC$. Such a test configuration is called a \demph{product test configuration}.
\end{defn}

\begin{rmk}
    Given a polarized Fano fibration $(\pi:X \to Y, \, \xi)$ and a function $v: \ft^* \to \bbR_+$, one can formally define 
    \[ \bbW_v^\mathrm{an}(\xi) = \frac{1}{(2\pi)^n}\int_X v(\mu)\, \frac{\omega^n}{n!}, \qquad \bbW_{v}(\xi) = \lim_{m \to \infty} \frac{1}{m^n}\sum_{\alpha \in \ft^*} v(\textstyle\frac{\alpha}{m}) \dim R_{m, \alpha} \p\]
    In the compact case where $Y$ is a point, this gives rise to the theory of \emph{weighted solitons} (or $v$-solitons), which are solutions to the equation 
    \[ \Ric_\omega - \omega = i\partial \bp \log(v(\mu)), \]
    see for example \cite{HanLi, LahdiliWeighted}. The Kähler--Ricci shrinker equation is recovered by taking $v(\mu) = e^{-\langle \xi, \, \mu \rangle}$. In general, if $v$ decays sufficiently fast at infinity, then the expressions above can be made to converge, which was the perspective taken in \cite{CWeighted} in the toric case. It is reasonable to expect that a similar theory for $v$-solitons can be established in this way for general polarized Fano fibrations.
\end{rmk}

\subsection{PSh functions and hermitian metrics}

Let $X$ be an irreducible complex analytic space, and $L \to X$ be a holomorphic line bundle. Throughout the paper we will use local additive notation for hermitian metrics on $L$:

\begin{defn}\label{notation:bad}
     By definition, a \demph{nonnegative singular hermitian metric} on $L$ consists of a collection $\{(U_\alpha, \, \phi_\alpha)\}$ such that the $U_\alpha$ are coordinate trivializations of $L$ covering $X$, and the $\phi_\alpha \in \PSH(U_\alpha)$ are consistent on the overlaps. This data determines a hermitian metric on $L$, denoted by $e^{-\phi}$, by setting
     \[ |v|^2_{e^{-\phi}} \underset{\mathrm{loc}}{=} e^{-\phi_\alpha} |v|^2\]
     for any $v \in L|_{U_\alpha} \cong U_\alpha \times \bbC$. We let $\PSH(X; L)$ be the space of singular hermitian metrics on $L$, and as a shorthand we will often write $\phi \in \PSH(X; \, L)$.
\end{defn}

Any $\phi \in \PSH(X, L)$ has a well-defined curvature $\omega = i\partial \bp \phi$, which is a nonnegative $(1, 1)$-current on $X$.
If we pick a smooth positively curved hermitian metric $\psi \in \PSH(X; L)$ with curvature $\omega$ we can identify $\PSH(X, \omega)$ with $\PSH(X; L)$ via $u \mapsto \psi + u$.

The following lemma will be useful later on when applied to sequences of singular metrics. Since the result is local, and any $\phi \in \PSH(X;L)$ is locally represented by a plurisubharmonic function, for simplicity we state it only in this case.

\begin{lem} \label{thm:notHarnack}
    Let $X$ be an irreducible analytic space and $Z \subsetneq X$ a proper analytic subset.
    Around every point $p \in Z$, there exists
    a neighborhood $U \subset X$ of $p$ and a compact set $K \subseteq X \setminus Z$ such that
    \[
        \phi|_{U} \leq \sup_{K} \phi
    \]
    for every $\phi \in \PSH(X)$.
\end{lem}
\begin{proof}
    Passing to a resolution of singularities of $X$ we can assume that $X$ is smooth. Then since the claim is local we can assume that $X$ is an open subset of $\bbC^n$, and $p = 0$. There exists a $v \in \bbC^n$ such that the circle $U(1) v$ in the complex plane $\bbC v$ avoids $Z$. This can be seen for example by using the parametric transversality lemma \cite[Thm.~6.35]{Lee} together with the fact that $Z$ is a union of submanifolds of real dimension $\leq 2n - 2$. Then since $U(1) v$ is compact and $Z$ is closed, there is a positive distance $\epsilon$ between the two sets; thus for any $q \in \bar{B}_{\epsilon/2}(0)$, the circle $q + U(1) v$ avoids $Z$. Set $K_p = \bar{B}_{\epsilon/2}(0) + U(1) v$. Then by the submean-inequality we have for any $z \in U_p \coloneq B_{\epsilon/2}(0)$
    \[
        \phi(z) \leq \frac{1}{2\pi} \int_0^{2\pi} \phi(z + e^{i\theta} v) \d \theta \leq \sup_{z + U(1) v} \phi \leq \sup_{K_p} \phi
        \c
    \]
    as claimed.
\end{proof}

\begin{cor}\label{thm:corNotHarnack}
    Let $X$ be an irreducible analytic space and $Z \subsetneq X$ a proper analytic subset. If $\phi_k \in \PSH(X) \cap L^{\infty}_\mathrm{loc}(X)$ is a sequence which is uniformly locally bounded on $X \setminus Z$, then $\phi_k$ is uniformly locally bounded on all of $X$.
\end{cor}

\subsection{Kähler cones and AC metrics}

By definition, a Riemannian cone $(Y^\circ, \, g_Y)$ is a smooth manifold $Y^\circ \cong \bbR_+ \times L$ for a compact manifold $L$ with a Riemannian metric $g_Y$ given by
\[ g_Y = dr^2 + r^2 g_L,\]
for a metric $g_L$ on $L$. We say that $g_Y$ is a \demph{Kähler cone} metric if $Y^\circ$ has an integrable complex structure $J$ making $g_Y$ Kähler. In this case the vector field $r \partial_r$ is real-holomorphic, and we define the \demph{Reeb field} $\xi = J\left( r \partial_r \right)$. Moreover, in this setting we can always obtain the Kähler form $\omega_Y$ as the curvature of the hermitian metric $e^{-\frac{r^2}{2}}$ on the trivial bundle:
\[\omega_Y = \frac{i}{2}\partial \bp r^2.\]

Now by a result of van Coevering \cite[Theorem 3.1]{vanCoevering:examples}, any Kähler cone is biholomorphic to the regular part of a normal affine variety $Y$ with a unique isolated singularity $Y = Y^\circ \cup \{o\}$. In particular, any Kähler cone $Y$ together with its Reeb field $\xi$ is a polarized affine cone. It was shown by Collins--Székelyhidi \cite{ColSzeStab, ColSze2} (see also \cite{LWX:metrictangent}) that a polarized affine cone admits a Ricci-flat Kähler cone metric if and only if it is K-polystable.

The main setting of this paper will be that of asymptotically conical Kähler metrics. In general, an asymptotically conical Kähler metric is a Kähler manifold $(X, \, J, \, g)$ such that there exists a Kähler cone $(Y^{\circ}, \, J_Y, \, g_Y)$ with radial function $r$, compact subsets $B_X \subset X$ and $B_Y \subset Y$, and a diffeomorphism $\Psi: X\setminus B_X \to Y\setminus B_Y$ making both $\Psi_*g - g_Y$ and $\Psi_*J - J_Y$ small with derivatives as $r \to \infty$. Notice that any such $g$ inherits from $g_Y$ the property that the norm of the curvature decays quadratically at infinity. If $g$ is a Kähler--Ricci shrinker, however, we have the following key result of Conlon--Deruelle--Sun \cite{CDS}, which says that quadratic curvature decay is actually sufficient for $g$ to be asymptotically conical in an even stronger sense:

\begin{thm}[{\cite[Theorem A]{CDS}}]\label{thm:quadraticdecaytoAC}
    Suppose that $(X, \, g, \, \xi)$ is a Kähler--Ricci shrinker with quadratic curvature decay. That is,
    \begin{equation}\label{eq:quadraticcurvaturedecay}
       \sup_{x \in X} \left| |\Rm_g|_g(x) \, d_g(p, \, x)^2\right| < \infty.
    \end{equation}
    Then there exists a Kähler cone $(Y, \, g_Y)$ and a holomorphic resolution of singularities $\pi: X \to Y$ such that $-K_X$ is $\pi$-ample, $\pi_*\xi = r\partial_r$, and
    \begin{equation}\label{eq:AC}
        \left| \left(\nabla^Y\right)^k  \, (\pi_*g - g_Y)\right|_{g_Y} \mkern-4mu \leq C r^{-2-k}.
    \end{equation}
    In particular, $(\pi:X \to Y, \, \xi)$ is a polarized Fano fibration.
\end{thm}

Given this, we will only work with metrics which are asymptotically conical in this strong sense.

\begin{defn}\label{def:AC}
    Let $(X, \, \xi)$ be a complex manifold with a fixed real holomorphic vector field $\xi$. For the purposes of this paper, we say that a Kähler metric $g$ on $X$ is \demph{asymptotically conical} (AC) if there is a resolution $\pi:X \to Y$ with $\pi_*\xi = r\partial_r$ as in \cref{thm:quadraticdecaytoAC} and a Kähler cone metric $g_Y$ on $Y$ satisfying \cref{eq:AC}. Similarly, we will say that a hermitian metric $\phi \in \PSH(X;L)$ on a line bundle $L \to X$ is AC if its curvature $\omega_\phi = i\partial \bp \phi$ is an AC Kähler metric in this sense.
\end{defn}

\section{Test configurations}
\label{sec:testconf}

The goal of this section is to show that every $T$-equivariant special test configuration $(\Pi: \cX \to \cY,\, T,\,  \eta)$ for a polarized Fano fibration $(\pi:X \to Y, \, \xi)$ can be equivariantly embedded in $\ambC$, and then use the ambient structure to study $\cX$. In \cref{sec:geodesics}, we will use the ambient $\amb$ to construct an associated smooth asymptotically conical subgeodesic ray in the case when $\pi:X \to Y$ is a resolution, which will be a crucial starting point for the constructions in the rest of the paper. Moreover, we will discuss the behavior of some more general families $\cX$ embedded in $\ambC$ which we will need to consider in future sections.

\subsection{Ambient spaces}\label{ss:ambient}

Recall from \cref{ex:amb} that given a compatible $T$-action on $\amb$ and a choice of Reeb field $\xi \in \ft$, we can view $(\pi:\amb \to \bbC^{N_2}, \, \xi)$ as a polarized Fano fibration. Let $(\Pi:\cX \to \cY, \, T, \, \eta)$ be a special test configuration for a polarized Fano fibration $(\pi: X \to Y, \, \xi)$. The goal of this section is to construct a $T$-action on $\amb$ and a $T' = T \times \bbC^*$-equivariant embedding $\cX \hookrightarrow \ambC$ such that the fibers over each $\tau \in \bbC$ give rise to an embedding $(\pi:X_\tau \to Y_\tau, \, \xi) \hookrightarrow (\pi:\amb \times\{\tau\} \to \bbC^{N_2} \times \{\tau\}, \, \xi)$ of polarized Fano fibrations. First, we have a technical lemma which says that, given a special test configuration $(\Pi:\cX \to \cY, \, T, \, \xi)$, the anticanonical bundle $-K_{\cX}$ is isomorphic to the relative anticanonical bundle $-K_{\cX/\bbC}$.

\begin{lem}\label{thm:adjunction}
    Let $(\pi:X \to Y, \, \xi)$ be a polarized Fano fibration and $(\Pi: \cX \to \cY, \, T, \, \eta)$ be a special test configuration. Then $K_{\cX}|_{X_\tau} = K_{X_\tau}$ for all $\tau \in \bbC^*$. Moreover the central fiber $X_0$ is automatically $\bbQ$-Gorenstein, and we have $K_{\cX}|_{X_0} = K_{X_0}$ as well.
\end{lem}
\begin{proof}
    For $\tau \neq 0$ this is clear since $\cX^* \cong X \times \bbC^*$.

    \begin{claim*}
        $\cX_\mathrm{sing} \cap X_0 \subseteq (X_0)_\mathrm{sing}$.
    \end{claim*}
    \begin{pf}
        Let $p \in (X_0)_\mathrm{reg}$ and embed the germ $(\cX, p)$ into some $(\bbC^N, 0)$, so that $(X_0, p) = \cV(J)$ and $(\cX, p) = \cV(I)$. By definition we have $J = I + (\Pi_\cX)$. By the assumption of regularity, the maximal ideal of $\cO_{X_0, p}$ is generated by $n = \dim X_0$ elements, so
        \[
            \fm = (f_1, \dots, f_n) + J = (f_1, \dots, f_n, \Pi_\cX) + I
        \]
        where $\fm$ is the maximal ideal of $\cO_{\bbC^N, 0}$. Thus the maximal ideal of $\cO_{\cX, p}$ is generated by $n+1 = \dim \cX$ elements, showing that $p \in \cX_\mathrm{reg}$.
    \end{pf}
    Since the central fiber $X_0$ is normal by assumption, we know that $(X_0)_\mathrm{sing}$ and therefore also $\cX_\mathrm{sing} \cap X_0$ have codimension 2 in $X_0$. This essentially allows us to apply the adjunction formula to $X_0 \subseteq \cX$:
    By the claim, the manifold $M \coloneq \cX_\mathrm{reg} \setminus (X_0)_\mathrm{sing}$ contains $D \coloneq (X_0)_\mathrm{reg}$ as a divisor so the adjunction formula for $M$ implies
    \[
        K_D = (K_M + D)|_D \cong K_M|_D,
    \]
    since the divisor class of $D$ is trivial. Since $\cX$ is $\bbQ$-Gorenstein, there exists an integer $p \in \bbZ_{\geq1}$ such that $p K_{\cX}$ is a line bundle. Hence we can view the above equality as
    \[ \left.\left(pK_{\cX} \right)\right|_{D} = \left.\left(pK_{X_0} \right)\right|_{D}\]
    Since $X_0$ is normal, this equality extends over the codimension 2 set $X_0 \setminus D$. Thus we see that $X_0$ is $\bbQ$-Gorenstein with $K_{X_0} = K_{\cX}|_{X_0}$.
\end{proof}

With this in place, we can prove the main result of this section:

\begin{prop}\label{thm:testconfigurationscanbeembedded}
    Let $(\pi:X \to Y, \, \xi)$ be a polarized Fano fibration and $(\Pi: \cX \to \cY, \, T ,  \, \eta)$ be a special test configuration. Then there exists a $T'$-action on $\left(\ambC, \, \cO_{\amb}(1)\right)$ and
    $T'$-equivariant embeddings $\iota_{\cX}: \cX \to \ambC$, $\iota_{\cY}: \cY \to \bbC^{N_2}\times \bbC$
    such that the diagram
        \begin{equation}
        \begin{tikzcd}%
         \cX \arrow[dd, "\Pi"'] \arrow[rr, "\iota_{\cX}"] \arrow[dr, "\Pi_{\cX}"] &[+20pt]  &[-5pt]  \ambC \arrow[dd, "\pr_2"] \arrow[dl] \\
         & \bbC & \\
         \cY   \arrow[rr, "\iota_{\cY}"] \arrow[ur, "\Pi_{\cY}"] &  &  \bbC^{N_2} \times \bbC \arrow[ul] \\   
    \end{tikzcd}
    \end{equation}
    commutes, and $\iota_{\cX}^*\mathcal{O}_{\amb}(1) = -pK_{\cX/\bbC} \cong -pK_{\cX}$. Here
    $K_{\cX/\bbC} \coloneqq K_{\cX} - \Pi^*K_{\bbC}$, and we choose the lift of the $T'$-action to $\mathcal{O}_{\amb}(1)$ which restricts to the canonical lift to $-pK_{\cX/\bbC}$ on $\cX$.
\end{prop}

\begin{proof}
    Note that $\cY$ is normal, as $\Pi_{\cY}:\cY \to \bbC$ is flat and the fiber over each closed point $Y_{\tau}\coloneqq \Pi_{\cY}^{-1}(\tau)$ is normal \cite[Prop.~6.8.3]{EGA4-2}. We first claim that $\cY$ is affine. By Sumihiro's theorem \cite{Sumihiro}, we can cover $\cY$ by $T'$-invariant affine open subsets. Choose some such invariant affine open set $\cU \subset\cY$ containing $o \in \cY$, the unique fixed point of the $T'$-action. Notice that $o$ is also the unique fixed point of the $T$-action on $Y_0 \subset \cY$. We claim that in fact $\cU = \cY$. %
    For each $t \in \bbC$, let $o_{\tau} \in Y_{\tau}$ be the unique fixed point for the $T$-action on $Y_{\tau}$. Since $(Y_{\tau}, \, \xi), \, (Y_0, \, \xi)$ are polarized affine cones by assumption, $o_{\tau}$ and $o$ can be identified as the zero sets of $\xi$ restricted to $Y_{\tau}, \, Y_0$. Thus we see that the union $\cup_{\tau \in \bbC} \{o_{\tau}\}$ coincides set-theoretically with the zero set $Z_\xi \subset \cY$ of $\xi$ on $\cY$. By the $\bbC^*$-equivariance of $\Pi_{\cY}: \cY \to \bbC$ and the fact that the $T'$-action must preserve $Z_\xi$, it follows that, for any given $t$, $Z_\xi \setminus\{o\}$ coincides with the orbit $T'\cdot o_{\tau}$, and therefore $Z_\xi = \overline{T'\cdot o_{\tau}}$. From this we see that $Z_\xi \subset \cU$. Finally, fix any point $p \in \cY$, which  lies in some fiber $Y_{\tau}$. Again using that $(Y_{\tau}, \, \xi)$ is a polarized cone we know that $o_{\tau} \in \overline{T \cdot p}$ so that in fact $p \in \cU$.

    Using that $\cY = \cU$ is affine, we now construct $\iota_{\cY}:\cY \to \bbC^{N_2}\times\bbC$.  %
    Write $\cY = \Spec(\cA)$, where we can decompose the ring $\cA$ into weight spaces:
    \[ \cA = \bigoplus_{\beta \in \Lie(T')^* }\cA_\beta. \]
    Of course, under our identification $T' = T \times \bbC^*$, any $\beta \in \frak{t}'^* \coloneqq \Lie(T')$ can be decomposed as  $\beta = (\alpha, \, b) $ for $\alpha \in \frak{t}^*$ and $b \in \bbZ$. In particular, the function $\Pi_{\cY}: \cY \to \bbC$ lies by assumption in the weight space associated to $(0, \,1)$. Now we can choose  $f_1, \dots, f_N \in \cA$ such that each $f_i$ lies in a weight space $\cA_{\beta_i}$, and such that $f_1, \dots, f_N, \, \Pi_{\cY}$ generate $\cA$.
    Therefore $f_1, \dots, f_N, \, \Pi_{\cY}$ induce an embedding $\iota_{\cY}:\cY \to \bbC^N \times \bbC$, equivariant with respect to the $T'$-action on $\cY$ and the diagonal action on $\bbC^N \times \bbC$ with weights $\beta_1, \dots, \beta_N, \, (0, 1)$. By construction $\Pi_{\cY}$ is precisely the composition of $\iota_{\cY}$ with the projection $\bbC^N \times \bbC \to \bbC$, and hence $\iota_{\cY}$ has the desired properties.

    Let us denote $\cX^* \coloneqq \cX \setminus X_0 \cong X \times \bbC^*$. By \cref{def:testconfiguration2}, we clearly have that $\left.-K_{\cX}\right|_{\cX^*}$ is $\Pi$-ample. By \cite[Theorem 1.7.8]{Lazarsfeld1}, we can conclude that $-K_{\cX}$ is globally $\Pi$-ample if we can show that $\left.-K_{\cX}\right|_{X_0}$ is $\Pi$-ample. Since $X_0 = \Pi_{\cX}^{-1}(0)$ is a principal divisor in $\cX$, we get from adjunction (see \cref{thm:adjunction}) that $-K_{X_0} = \left.\left(-K_{\cX} - \cO_{\cX}(X_0)\right)\right|_{X_0}  =\left.-K_{\cX}\right|_{X_0}$. Since $-K_{X_0}$ is $\Pi:X_0 \to \cY_0$ ample by assumption, it follows that $-K_{\cX}$ is $\pi:\cX \to \cY$ ample. By \cite[Remark 1.7.4]{Lazarsfeld1} (see also \cite[p. 120]{Hart}), therefore, we obtain an embedding $\iota_1:\cX \to \bbP^{N_1} \times \cY$ such that $-pK_{\cX} = (\iota_1\circ\pi_1)^*\cO_{\bbP^{N_1}}(1)$, for some $p \geq 1$.

    If we let $\tau$ be the unique $\bbC^*$-invariant holomorphic coordinate on $\bbC$, then $K_{\bbC}$ admits a canonical trivialization given by $d\tau$. This gives rise to an isomorphism $pK_{\cX} \cong pK_{\cX/\bbC}$, and also a lift of the $T'$-action on $\cX$ to $pK_{\cX/\bbC}$. Note that this lift differs from the canonical action on $pK_{\cX}$. Indeed, if $\Omega$ is a local section of $pK_{\cX/\bbC}$ which is a weight vector, then the corresponding section of $K_{\cX}$ can be identified with $(d\tau)^p \wedge \Omega$, whose weight for the $\bbC^*$-action clearly differs by a factor of $p$.

    Let $\cA = \bbC[x_1, \dots, x_{N_2}]/I_{\cY}$ be the affine coordinate ring of $\cY$. Set
    \[\cR \coloneqq \bigoplus_{m\geq0} \cR_m, \qquad \cR_m \coloneqq H^0(\cX, \, -mpK_{\cX/\bbC}),\]
    and note that $\cR_0 = \cA$ (c.f. \cite[Section 3]{SunZhang}). Therefore
    \[ \bigg(\Pi:\cX \to \cY\bigg)  = \bigg(\pr:\Proj \cR \to \Spec \cA \bigg).\]
    Note that since $\Pi:\cX \to \cY$ is projective, we have that $\cR$ is finitely generated as an $\cA$-algebra by elements of $\cR_1$ \cite[II, Cor.~5.16, (b)]{Hart}. In particular, $\cR_1$ is a finitely-generated $\cA$-module. Hence, if we take a minimal generating set $s_0, \dots, s_{N_1} \in \cR_1$ of $\cR_1$ as an $\cA$-module, then $s_0, \dots, s_{N_1}$ will generate $\cR$ as an $\cA$-algebra. If we let
    \[\bbP^{N_1}\times \cY = \Proj \, \left( \cA \otimes \bbC[x_0, \dots, x_{N_1}]\right),\]
    then the embedding $\iota_1:\cX \to \bbP^{N_1}\times \cY$ is given simply by the map
    \[ \iota_1^*: \cA \otimes \bbC[x_0, \dots, x_{N_1}] \to \cR, \qquad \iota_1^*(x_i) = s_i. \]
   
    We claim that $s_0, \dots, s_{N_1}$ can be taken to be weight vectors for the $T'$-action on $\cX$ and its lift to $-pK_{\cX/\bbC}$. Indeed, we can write
    \[\cR_1 = \bigoplus_{\alpha \in \frak{t}^*} \cR_{1,\alpha},\]
    so that
    \[s_i = \sum_{j=1}^{k_i} \sigma^i_{j,\alpha_j^i}, \qquad v\cdot\sigma^i_{j,\alpha_j} = \langle v, \, \alpha_j^i\rangle \sigma^i_{j, \alpha_j^i} \quad \textnormal{for all } v \in \frak{t}.\]
    Then clearly $\sigma^1_{1, \alpha_1}, \dots, \sigma^{N_1}_{N_1, \alpha_{N_1}}$ generate $\cR$, and so the map  $\iota_1^*: \cA \otimes \bbC[x_0, \dots, x_{\tilde{N}_1}] \to \cR $ sending $x_0, \dots, x_{\tilde{N}_1}$ to $\sigma^1_{1, \alpha_1^1}, \dots, \sigma^{N_1}_{N_1, \alpha^{N_1}_{N_1}}$ is surjective. Since $\cA \otimes \bbC[x_0, \dots, x_{\tilde{N}_1}]$ is generated in degree 1, it follows that the induced map $\tilde{\iota}:\cX \to \bbP^{\tilde{N}_1} \times \cY$ is a closed embedding such that $\tilde{\iota}^*\cO_{\bbP^{\tilde{N}_1}\times \cY}(1) = -pK_{\cX/\bbC}$ (see  \cite[{\href{https://stacks.math.columbia.edu/tag/01MX}{Tag 01MX}}]{stacks}). %
    Moreover this map is $T$-equivariant with respect to the diagonal $T$-action on $\bbP^{\tilde{N}_{1}}_{\cY}$, lifting the $T$-action on $\cY$, whose weight on each homogeneous coordinate $x_k$ is given by the corresponding $\alpha_j^i$.
    By relabeling, we assume henceforth that $\iota_1: \cX \to \bbP^{N_1}\times \cY$ is $T$-equivariant, where $T$ acts on $\bbP^{N_1}\times \cY$ by
    \[ \tau\cdot\bigg( [x_0:\dots:x_{N_1}],\,  y \bigg) = \bigg( [\tau^{\alpha_0}x_0:\dots:\tau^{\alpha_{N_1}}x_{N_1}], \, \tau\cdot y \bigg).\]
    The maps $\iota_{\cY}$ and $\iota_{\cX} = \iota_{\cY} \circ \iota_1$ then satisfy the required properties.
\end{proof}

\begin{cor}\label{thm:testconfigurationsarepolarizedfanofibrations}
    Let $(\Pi: \cX \to \cY, \, T, \, \eta)$ be a special test configuration for a polarized Fano fibration $(\pi: X \to Y, \, \xi)$. For any $\epsilon > 0$ sufficiently small, $(\Pi: \cX \to \cY, \, \xi+\epsilon \eta)$ is itself a polarized Fano fibration.
\end{cor}
\begin{proof}
    From the proof of \cref{thm:testconfigurationscanbeembedded}, we can easily see that for any $\epsilon > 0$ sufficiently small, $(\cY, \, \xi + \epsilon \eta)$ is a polarized affine cone. The key point is that $\xi + \epsilon \eta$ is indeed a Reeb field (c.f. \cite[Lem.~7.3]{ColSze2}). Using the notation above, we can see that the induced $T'$ action on $\bbC^{N_2} \times \bbC$ splits as a $T'$ action on $\bbC^{N_2}$ with weights $\beta_1, \dots, \beta_{N_2}$ and the standard $\bbC^*$ action on $\bbC$. Suppose that each $\beta_j$ is given by $\beta_j = (\alpha_j, \, b_j)$ with respect to the decomposition above. Then the $T'$-action on $\bbC^{N_2}$ further splits as a $T$-action with weights $\alpha_1, \dots, \alpha_{N_2}$ and a $\bbC^*$-action with weights $(b_1, \dots, b_{N_2})$. In addition, we must have that $\langle \xi, \, \alpha_j \rangle > 0$. This follows since, by construction, each $f_j$ restricts to any given fiber $Y_{\tau}$ to a weight vector for the $T$-action on $Y_{\tau}$ with weight $\alpha_j$, together with the fact that $\xi$ is a Reeb field on $Y_{\tau}$.  We read off immediately that
    \[\langle \xi + \epsilon \eta  , \, \beta_j \rangle
 = \langle \xi  , \, \alpha_j \rangle + \epsilon b_j,\]
 which is positive as long as $\epsilon$ is sufficiently small. Since we clearly have that $\langle \xi + \epsilon \eta \, , \, (0, \,1) \rangle = \epsilon > 0$, the conclusion follows.

The only missing point is to see that $\cX$ is klt. This however follows by the argument of \cite[Lem.~2.2]{BermanQFano}, since $X_0$ is reduced and klt.
\end{proof}

\subsection{Filtrations and degenerations}%
\label{sub:test_configurations_and_filtrations}

Let $S  = \bbC[x_0, \dots, x_{N_1}] \otimes \bbC[y_1, \dots, y_{N_2}]$ be the ring of homogeneous functions of the ambient space $\amb$, so that
\[
    \amb \cong  \Proj (S) \p
\]
An important point, and possible source of confusion, in this discussion will be that $S$ carries multiple (compatible) gradings. Unless specified otherwise, the default grading of $S$ will be the grading by "$x$-degree" for us.
We will be interested in the situation where we have a $T$-equivariantly embedded variety $X \hookrightarrow \amb$, and a $\bbC^*$-action on $\amb$ commuting with $T$ which gives rise to a degeneration of $X$ in $\amb$.

\begin{convention}
    Given a $\bbC^*$-action $\rho: \bbC^* \times X \to X$ on a variety $X$, we say that a function $f$ has weight $k$ if $f(\rho(\tau)^{-1} x) = \tau^{-k} f$ for every $x \in X, \, \tau \in \bbC^*$.
\end{convention}
Note that this differs from another typical convention, namely declaring that $f$ has weight $k$ if $f(\rho(\tau)^{-1} x) = \tau^k f$. Our definition however is consistent with saying that a function $f$ satisfying $\cL_\eta f = k f$ has weight $k$ for the vector field $\eta$.

The goal of this section is to state the well-known characterization of such $\bbC^*$-degenerations of $X \subset \amb$ in terms of the algebra of $S$, in our current setting. Suppose then that we have such an $X$ and we equip $\amb$ with a $\bbC^*$-action $\rho$. Since any subscheme of $\amb$ is defined by an ideal, there is a corresponding $I \ideal S$. Moreover, the ring $S$ has a grading $S = \bigoplus_w S_w$ by the weights of $\rho$, and this descends to a filtration of $I$: %
    \begin{align}
        F_i S &= \bigoplus_{w \leq i} S_w , \\
        F_i I &= I \cap F_i S,
    \end{align}
    so that $F_i I$ contains those $f \in I$ whose weight components all have weight $w \leq i$. If we degenerate $X$ by $\rho$, we obtain a family $X_{\tau} \subset \amb$, where $X_{\tau} \cong X = X_1$. The functions on $X$ are related to those on $X_{\tau}$ via the map $\theta: S \to S[\tau, \tau^{-1}]$ given by
    \[
        \theta(f)(\tau, p) \coloneq f(\rho(\tau)^{-1} p)
        \p
    \]

    For any $f \in S$ we define $\ini(f)$ to be the term of leading order in $\tau$ in $\theta(f)$, so that $\theta(f) = \tau^{-k} \ini(f) + O(\tau^{-k+1})$. Equivalently, $\ini(f)$ is the component of $f$ with highest $\rho$-weight. The initial term $\ini(f)$ can be viewed as those terms in $f$ contributing to the most negative power in the Laurent expansion with respect to $\tau$ of $\theta(f)$. %
    We are interested in the limit $X_0\subset\amb$ of these varieties $X_\tau$ as $\tau \to 0$. To make sense of this, we define $\cX^* \subset \ambC$ to be the union of all $(X_{\tau}, \, \tau)$, and define
    \begin{equation}\label{eq:degenerationtotalspacedefinition}
        \cX = \overline{\cX^*} \subset \ambC.
    \end{equation}
    Since the projection $\ambC \to \amb$ induces a map $\Pi_{\cX}:\cX^* \to \bbC^*$, it's clear that the closure admits a map $\Pi_{\cX}: \cX \to \bbC$. We define $X_0$ to be the scheme-theoretic fiber in $\cX$ over $0 \in \bbC$, which we view as embedded $X_0 \hookrightarrow \amb$. Since the embedding $X \hookrightarrow \amb$ is $T$-equivariant, there is a natural $T' = T \times \bbC^*$-action on $\ambC$ leaving $\cX$ invariant.

    \begin{defn}\label{thm:flatlimit}
        In this situation, we say that $X_0$ is the \demph{flat limit} of $X$ with respect to $\rho$.
    \end{defn}

\begin{lem} \label{thm:degenerationRings}
   The ideal of the flat limit $X_0 \subset \amb$ is given by the initial ideal $I_0$ associated to $\rho$ and $I$, namely
    \begin{equation} \label{eq:I0ini1}
        I_0 = \left(\ini f \mid f \in I\right) \ideal S.
    \end{equation}
    Moreover, $I_0$ can be identified with
    \begin{equation}\label{eq:I0ini2}
        I_0 \cong \gr I = \bigoplus_i F_i I / F_{i-1} I \ideal \bigoplus_w S_w.
    \end{equation}
\end{lem}

\begin{proof}
    The variety $\cX^* \subseteq \amb \times \bbC^*$ is isomorphic to $X \times \bbC^*$ by
    \[
        R(p, \tau) = (\rho(\tau)p, \tau)
        \c
    \]
    Since $R$ is invertible it defines a pushforward map on functions $R_*: S[\tau, \tau^{-1}] \to S[\tau, \tau^{-1}]$, given by $R_*(g) = g \circ R^{-1}$. Note that for $f \in S$ we have $R_*(f \tau^0) = \theta(f)$.
    $\cX^*$ is cut out by the ideal
    \[
        J^* = R_* I[\tau, \tau^{-1}] \ideal S[\tau, \tau^{-1}]
    \]
    The ideal of $\cX$ is $J \coloneq J^* \cap S[t] \ideal S[t]$ since a (homogeneous) function $f \in S[t]$ vanishes on $\cX^*$ iff it vanishes on $\cX$.

    By construction, functions in the image of $\theta$ have weight 0 for the $T'$-action on $\bbP^{N_1} \times \bbC^{N_2} \times \bbC^*$ which
    preserves $\cX$. Since $I[t, t^{-1}]$ is generated as an abelian group %
    by $\{f \tau^k \mid f \in I, k \in \bbZ \}$, we can conclude that $J^*$ is generated as an abelian group by $R_*(f t^k) = \theta(f) t^k$, which have weight $k$ for the induced action on $\cX$.
    Hence the weight $k$ components of $J^*$ are
    \begin{equation}
        J^*_k = \left\{\tau^k \theta(f) \mid f \in I\right\}
        \p
    \end{equation}
    Any $f \in F_i I$ can be weight-decomposed in $S$, as $f = \sum_{w \leq i} g_w$, where $g_w \in S_w$ and $g_i \neq 0$. Then
    \[
        \theta(f) \tau^k = \sum_{w \leq i} \tau^{k-w} g_w,
    \]
    and we see that $\tau^k \theta(f) \in S[\tau]$ if and only if $k \geq i$.
    This allows us to identify the weight $w$ component of $J = J^* \cap S[\tau]$ as
    \begin{equation}
        J_k = J^*_k \cap S[\tau] = \left\{\tau^k \theta(f) \mid f \in F_k I\right\} = \left\{R_*(\tau^k f) \mid f \in F_k I\right\}
        \p
    \end{equation}
    Algebraically, $I_0$ is obtained from $J$ by taking the image under the quotient map $p: S[\tau] \to S, \,  \tau \mapsto 0$.
    To see \cref{eq:I0ini1} we note that since $\tau \in S[\tau, \tau^{-1}]$ the functions $\tau^k \theta(f)$ for $f \in F_k I \setminus F_{k-1} I$ already form a generating set of $J$. But for such $f$, we precisely have $\tau^k \theta(f) = \ini f$.
    Thus
    \[
        I_0 = \left( p(\tau^k \theta(f)) \mid f \in F_k I \setminus F_{k-1} I \right)_S = (\ini f \mid f \in F_k I \setminus F_{k-1} I)_S.
    \]
    Furthermore, $I_0 = p(J) = J / \tau J$
    We have $R^*(J) = \bigoplus_i \tau^i F_i I$. Therefore $I_0 \cong R^*(J / \tau J) = R^*(J) / \tau R^*(J) = \gr I$.
\end{proof}

\subsection{Product test configurations}%
\label{sub:product_testconf}

The goal of this subsection is to prove
\begin{prop} \label{thm:isoTrivialProduct}
    Let $(\Pi: \cX \to \cY, \, T ,\,  \eta)$ be a test configuration for $(\pi: X \to Y, \, \xi)$ whose central fiber $X_0$ is $T$-equivariantly isomorphic to $X$. Then $\cX$ is
    equivariantly isomorphic to a product test configuration.
\end{prop}

\begin{notation}
    For convenience we will write
    \[
        \bbC[\underline{x}] \coloneq \bbC[x_0, \dots, x_{N_1}] \qtext{and}
        \bbC[\underline{y}] \coloneq \bbC[y_1, \dots, y_{N_2}]
    \]
    and $\bbC[\underline{x}]_m$ for the degree $m$ component of the polynomial ring. We also write $\bbC[X]$ for the ring of regular functions on a variety $X$ and $\bbC[X]^+ \coloneq \bbC[X] / \bbC$ for regular functions modulo constants.
\end{notation}
Recall from \cref{rmk:pffEmbed} that we have a $T$-equivariant embedding $\iota_X: \left(X, -p K_X\right) \to \left(\bbP^{N_1} \times \bbC^{N_2}, \cO(1)\right)$, where the ambient space has homogeneous coordinate ring $S = \bigoplus_{m\geq 0} S_m$, where
\[
    S_m \coloneq \bbC[x_0, \dots, x_{N_1}]_m \otimes \bbC[y_1, \dots, y_{N_2}],
\]
Also recall the definition
\[
    \cR(X) \coloneq \bigoplus_{m \geq 0} R_m(X) \coloneq \bigoplus_{m \geq 0} H^0(X, -m K_X)
    \p
\]
and that the embedding $\iota_X$ corresponds to a $T$-equivariant graded map $S \to \cR$ which is surjective onto the subring of $\cR(X)$ consisting of homogeneous functions.

\begin{defn}\label{def:linearlynormal}
    We say that the embedding is \demph{linearly normal in low weight} if the components of the map $S \to \cR$
    \[
        p: \bbC[\underline{y}]_1 \to \bbC[X]^+_{\wt_\xi \leq w_2}
    \]
    and
    \[
        q: \bbC[\underline{x}]_1 \to R_p(X)_{\wt_\xi \leq w_1}
    \]
    are surjective, where $w_2 = \max_i \wt_\xi(y_i)$ and $w_1 = \max_i \wt_\xi(x_i)$.
    Here $\bbC[X]^+_{\wt_\xi \leq w_2}$ denotes the sum of all $\xi$-weight spaces of $\bbC[X]^+$ of weight $\leq w_2$, and analogous for $R_p(X)_{\wt_\xi \leq w_1}$.

\end{defn}

\begin{ex}
    Consider $\bbP^1$ with the $\bbC^*$-action $\lambda \cdot[u \colonrel v] = [\lambda u \colonrel v]$, and $\bbC$ with its standard $\bbC^*$-action of weight $1$.
    Consider the anticanonical embedding $\iota_{-K}: \bbP^1 \times \bbC \to \bbP^2 \times \bbC$
    \[
        ([u \colonrel v], z) \mapsto ([u^2 \colonrel uv \colonrel v^2], z).
    \]
    Then the $\bbC^*$-action on $\bbP^1 \times \bbC$ is induced by the action on $\bbP^2 \times \bbC$ with weights (2, 1, 0, 1), and we take $\xi$ to be the vector field generated by this action. Then observe that the map
    \begin{equation} \label{eq:justSomeAut}
        ([u \colonrel v], z) \mapsto ([u + vz \colonrel v], z)
    \end{equation}
     cannot be induced by $\PGL(2) \times \GL(1)$.  %

    We see also that the embedding %
    is not linearly normal in low weight, since for example the anticanonical section $v^2z$ of weight 1 cannot be induced by a linear form on $\bbP^2 \times \bbC$. However, we can fix this by appropriately increasing the dimension of our target space.
    Indeed, observe that the modified embedding
    \[
        ([u \colonrel v ], z) \mapsto ([u^2 \colonrel uv \colonrel v^2 \colonrel u v z \colonrel v^2 z \colonrel v^2 z^2], z)
    \]
    obtained by composing $\iota_{-K}$ with the embedding $\bbP^2 \times \bbC \hookrightarrow \bbP^5 \times \bbC$ given by
    \[([a_1:a_2:a_3], z) \to ([a_1:a_2:a_3:a_2z:a_3z:a_3z^2], z),\]
    has the property that any anticanonical section on $\bbP^1 \times \bbC$ whose weight under the given $\bbC^*$-action is no bigger than 2 is induced by a linear form on $\bbP^5 \times \bbC$. Then \cref{eq:justSomeAut} is indeed induced by $\PGL(6)$:
    \[
        \begin{bmatrix}
            1 & 0 & 0 & 2 & 0 & 1 \\
            0 & 1 & 0 & 0 & 1 & 0 \\
            0 & 0 & 1 & 0 & 0 & 0 \\
            0 & 0 & 0 & 1 & 0 & 1 \\
            0 & 0 & 0 & 0 & 1 & 0 \\
            0 & 0 & 0 & 0 & 0 & 1
        \end{bmatrix}
        \begin{bmatrix} u^2 \\ uv \\ v^2 \\ uvz \\ v^2 z \\ v^2 z^2 \end{bmatrix}
            =
        \begin{bmatrix} u^2 + 2 u v z + v^2 z^2 \\ uv + v^2 z \\ v^2 \\ uvz + v^2 z^2 \\ v^2 z \\ v^2 z^2 \end{bmatrix}.
    \]
\end{ex}

In our setting, we can assume without loss of generality that the fibers of a test configuration are embedded in a way that is linearly normal in low weight:
\begin{lem} \label{thm:linearlynormal}
    Let $X, X_0$ be polarized Fano fibrations embedded $T$-equivariantly into $\amb$.
    By performing a $T$-equivariant linear embedding $\bbP^{N_1} \times \bbC^{N_2} \inj \bbP^{N_1'} \times \bbC^{N_2'}$,
    we can arrange for both embeddings $X, X_0 \inj \bbP^{N_1'} \times \bbC^{N_2'}$ to be linearly normal in low weight.
\end{lem}
\begin{proof}
    Since $\xi$ is a Reeb vector, we have that $\bbC[Y]^+_{\wt_\xi \leq w_2}$ is finite-dimensional. In particular, we can extend $p(y_1) , \dots , p(y_{N_2})$ to a basis $v_1, \dots, v_{N_2'}$ of $\bbC[Y]^+_{\wt_\xi \leq w_2}$, where $p$ is the map as in \cref{def:linearlynormal}.
    Then the map
    \begin{equation} \label{eq:ambAffine}
        \bbC[y_1, \dots, y_{N_2'}] \to \bbC[Y]
    \end{equation}
    extending $p$ is defined by $p(y_i) = v_i$. By construction, the linear polynomials in $y_1, \dots, y_{N_2'}$ surject onto $\bbC[Y]^+_{\wt_\xi \leq w_2}$. Similarly, $R_p(X)_{\wt_\xi \leq w_1}$ is finite-dimensional because $R_p(X)$ is a finitely generated $\bbC[Y]$-module. Thus we can add variables $x_{N_1 + 1}, \dots x_{N_1'}$ of $\xi$-weight $\leq w_1$ so that
    \begin{equation} \label{eq:ambProj}
        \bbC[x_1, \dots, x_{N_1'}]_1 \to R_p(X)_{\wt_\xi \leq w_1}
    \end{equation}
    is surjective. Thus the embedding $X \to \bbP^{N_1'} \times \bbC^{N_2'}$ is linearly normal in low weight.

    Repeating the same procedure for $X_0 \hookrightarrow \bbP^{N_1'} \times \bbC^{N_2'}$, we can make the embedding of $X_0$ linearly normal in low weight as well while preserving $w_1, w_2$ and linear normality in low weight of $X$.
\end{proof}

\begin{lem}\label{thm:getamatrix}
    If the embeddings $X, X_0 \inj \amb$ are linearly normal in low weight, then any isomorphism $f: X \to X_0$ can be induced by an element of $G \coloneq \PGL(N_1 + 1) \times \GL(N_2)$.
\end{lem}
\begin{proof}
    First note that the numbers $w_1, w_2$ are a property of the ambient space $\amb$, and thus the same for $X_0$ and $X$.
    The map $f$ induces an isomorphism $Y \to Y_0$, corresponding to a map $\gamma: \bbC[Y_0] \to \bbC[Y]$ of algebras.
    Now since subspace $\bbC[Y]^+_{\wt_\xi \leq w_2}$ is defined independent of the embedding, we have
    \[
        \phi(\bbC[Y_0]^+_{\wt_\xi \leq w_2}) = \bbC[Y]^+_{\wt_\xi \leq w_2}
        \p
    \]
    Thus writing $p_0$ and $q_0$ for the maps from \cref{def:linearlynormal} corresponding to the embedding of $X_0$, the composition
    \[
        \phi \circ p_0: \bbC[\underline{y}]_1 \surj \bbC[Y]^+_{\wt_\xi \leq w_2}
    \]
    is surjective.
    On the other hand, since $p: \bbC[\underline{y}]_1 \surj \bbC[Y]^+_{\wt_\xi \leq w_2}$ is surjective, $\phi \circ p_0$ admits a lift $\Phi: \bbC[\underline{y}]_1 \to \bbC[\underline{y}]_1$ along $p$. Since $\phi \circ p_0$ is surjective we can pick $\Phi$ to be invertible.
    Then we extend $\Phi$ to a map of rings $\Phi: \bbC[\underline{y}] \to \bbC[\underline{y}]$.
    This makes the diagram
    \begin{cd} \label{eq:diagAffine}
        \bbC[\underline{y}] \arrow[d, "p_0", two heads] \arrow[r, "\Phi", dashed]  & \bbC[\underline{y}] \arrow[d, "p", two heads] \\
        \bbC[Y_0] \arrow{r}{\phi}[swap]{\sim}  & \bbC[Y]
    \end{cd}
    commute because two maps out of a polynomial algebra agree if and only if they agree on linear polynomials.

    We now repeat essentially the same argument for $q$ instead of $p$.
    Since $X, X_0 \inj \amb$ are linearly normal in low weights, we have maps
    \begin{cd} \label{eq:diagProj}
        \bbC[\underline{x}]_1 \arrow[d, "q_0", two heads] \arrow[r, "F", dashed]  & \bbC[\underline{x}]_1 \arrow[d, "q", two heads]  \\
        R_1(X_0)_{\leq w_1} \arrow{r}{f^*}[swap]{\sim}  & R_1(X)_{\leq w_1}
    \end{cd}
    where the dashed linear map can be chosen to be an isomorphism, so $F \in \PGL(N_1 + 1)$.

    Finally, we need to verify that $F \times \Phi$ induces $f$, which amounts to checking the commutativity of
    \begin{equation}
    \begin{minipage}{0.4\textwidth}
    \begin{cd*}
        X \arrow[r, "f"] \arrow[d, "\Proj(q)"] & X_0 \arrow[d, "\Proj(q_0)"] \\
        \bbP^{N_1} \arrow[r, "T"]                & \bbP^{N_1}
    \end{cd*}
    \end{minipage}
    \begin{minipage}{0.4\textwidth}
    \begin{cd*}
    X \arrow[r, "f"] \arrow[d, "\Spec(p) \circ \pi"] & X_0 \arrow[d, "\Spec(p_0) \circ \pi_0"] \\
        \bbC^{N_2} \arrow[r, "L"]                & \bbC^{N_2}
    \end{cd*}
    \end{minipage}
    \end{equation}
    This follows immediately from \cref{eq:diagAffine,eq:diagProj}.
\end{proof}

\begin{proof}[Proof of \cref{thm:isoTrivialProduct}]
    By \cref{thm:testconfigurationscanbeembedded,thm:linearlynormal} we can embed $\cX \inj \amb \times \bbC$ such that the embeddings of $X_1$ and $X_0$ are linearly normal in low weights, and we consider the closure $\bar{\cX} \subseteq \bbP^{N_1} \times \bbP^{N_2} \times \bbC$. Since $\cX$ is irreducible, so is $\bar{\cX}$, therefore $\bar{\pi}: \bar{\cX} \to \bbC$ is a flat family. By definition, $\bar{\pi}$ then corresponds to a morphism  $h_{\bar{\pi}}: \bbC \to \Hilb(N_1, N_2, Q)$, to the Hilbert scheme of subschemes of $\bbP^{N_1} \times \bbP^{N_2}$ with the same Hilbert polynomial $Q$ as $\bar{X}_1$.

    By \cref{thm:getamatrix}, $X$ and $X_0$ lie in the same $G$-orbit of $\Hilb(N_1, N_2, Q)_\mathrm{red}$. This $G$-orbit is isomorphic to $G / \Stab_G(X_0) \eqcolon G / H$, so there is a map $h: \bbC \to G / H$. The principal $H$-bundle $P \coloneq (G \to G/H)$ can be pulled back to $\bbC$ and by the Oka principle we know that $h^* P$ is a trivial bundle. Furthermore, the degenerating $\bbC^*$-action $\lambda : \bbC^* \to G$ gives us a left $\bbC^*$-action on $G / H$ and, by construction, $h: \bbC \to G/H$ is $\bbC^*$-equivariant for this action.
    Hence $h^* P \to \bbC$ is a $\bbC^*$-equivariant principal $H$-bundle.
    By the equivariant Oka principle \cite{EquivOka,EquivOkaAlg},  there exists a $\bbC^*$-equivariant lift
    \begin{cd}
         & G \arrow[d] \\
        \bbC \arrow[r, "h"] \arrow[ur, "s", dashed]  & G / H
        \p
    \end{cd}
    By construction, $s$ has the property that $s(\tau) X_0 = X_\tau$. Thus $s$ provides a $\bbC^*$-equivariant isomorphism $X_0 \times \bbC \to \cX$.
\end{proof}

\section{Geodesics}%
\label{sec:geodesics}

    Throughout this section, we denote the unit disk in $\bbC$ by $\bbD$, the punctured disk by $\bbD^*$, and more generally a disk of radius $r$ by $\bbD_r$.
    Given a special test configuration $(\Pi: \cX \to \cY, \, T, \, \eta)$, we write $\cX^*$ for $\Pi^{-1}(\bbC^*) = \cX \setminus X_0$,  and similarly $\cX^*_k = \Pi^{-1}(\bbD \setminus \bbD_{e^{-k}})$ and $\cX_\bbD = \Pi^{-1}(\bbD)$.

\subsection{Function spaces}
To begin, we recall a useful framework from \cite{CarlosUniq}. The idea is to use the fact that the soliton vector field $-J\xi$ is nonvanishing on $X \setminus E \cong Y^\circ$ to cover this set by holomorphic charts in which weight vectors for $\xi$ behave well:
we can cover this region with \demph{equivariant charts}:
\begin{defn}[{\cite[Lem.~2.3]{CarlosUniq}}] \label{def:equivchart}
    There exist finitely many holomorphic charts
    \begin{equation} \label{eq:defEqCh}
        U_\alpha \cong \bbR \times i (-1, 1) \times U'_\alpha, \qquad U'_\alpha \subseteq \bbC^{n-1} \text{ open and precompact}
    \end{equation}
    covering $M \setminus E$, called \demph{equivariant charts}, such that $-J\xi|_{U_\alpha} = \frac{\partial }{\partial x_1^\alpha}$. Here $\bbR \times i(-1, 1)$ is equipped with coordinate $z_1^\alpha = x_1^\alpha + i y_1^\alpha$, and $U_\alpha'$ with $(z_2^\alpha, \dots, z_n^{\alpha})$.
    By slightly shrinking $U'_\alpha$ if necessary, we can arrange for the chart $U_\alpha$ to be contained in a larger chart $V_\alpha$ defined on $V_\alpha \cong \bbR \times i(-2, 2) \times V'_\alpha$, with $V'_\alpha \supseteq \overline{U'_\alpha}$.
\end{defn}

As a consequence, if $f$ is a function on $V_\alpha$ which is a weight vector for $\xi$, then in an equivariant chart $U_\alpha$ we can write
\[ f(z_1, \dots, z_n) = G_\alpha(z^\alpha_2, \dots, z^\alpha_n) e^{\lambda x^\alpha_1}\]
for a smooth bounded function $G_\alpha$ on $U_\alpha'$.

We write
\newcommand{\Hsing}{\bar{\cH}_2}
\begin{equation}
    \Hsing = \left\{ \psi \in \PSH\left(X; -K_X\right) ~\middle|~ \psi = \Theta(r^2) \right\}
\end{equation}
where $\psi = \Theta(r^2)$ means that the trivializations of $\psi$ in equivariant coordinates have a quadratic upper and lower bound.
We will mostly work within the nicer space
\begin{equation}
    \cH^*_2 =
    \left\{ \psi \in \Hsing \cap C^\infty ~\middle|~
            \exists \text{ cone metric } \omega_C = \frac{i}{2} \partial\bar{\partial} \tilde{r}^2 \text{ on $(C, \xi)$ s.t.\ }
            |i \partial\bar{\partial} \psi - \omega_C |_{\omega_C} = O(r^{-2}),
    \right\}
\end{equation}
and moreover we further denote
\begin{equation}
   (\Hsing)^T = \left\{ \psi \in \Hsing ~\middle|~ T-\textnormal{invariant} \right\}, \qquad (\cH^*_2)^T = \left\{ \psi \in \cH^*_2 ~\middle|~ T-\textnormal{invariant} \right\}.
\end{equation}

 The definition of $\cH_2^*$ here is slightly different from that in \cite[Def.~2.5]{CarlosUniq}, but the set $(\cH_2^*)^T$ is the same with both definitions, as can be seen from the proof of \cite[Lem.~2.6]{CarlosUniq}.

 We want to remark here that there is an easy way to obtain the asymptotic cone metric $\omega_C$ associated to any given $\phi \in \cH^*$. Indeed the proof of  \cite[Lem.~2.6]{CarlosUniq} shows:
\begin{lem}\label{thm:metricvsconeC0}
    For any $\phi \in (\cH_2^*)^T$,
     the radial function $r_\phi^2$ of the asymptotic cone metric is the unique function of $-J\xi$-weight one such that in every equivariant chart
    \[
        -\log \frac{e^{-\phi}}{\Omega \wedge \bar{\Omega}} = r_\phi^2/2 + O((\log r)^2)
    \]
\end{lem}

Write $\gamma_t$ for the time $t$ flow of $-J \xi$. Then in the situation of the lemma we see that
\[
    - e^{-2s} \log \frac{\gamma_s^* e^{-\phi}}{\Omega \wedge \bar{\Omega}} = r_\phi^2 + e^{-2s} O((\log r_\mathrm{BG})^2 + 2 s)
    \p
\]
Therefore $r_\phi^2$ is the locally uniform limit
\begin{equation} \label{eq:obtainr2}
    r_\phi^2 = -\lim_{s \to \infty} e^{2s} \log \frac{\gamma_s^* e^{-\phi}}{\Omega \wedge \bar{\Omega}}
    \p
\end{equation}
We need the following consequence of this formula later:
\begin{cor} \label{thm:comparePhiR}
    If $\phi_1 \geq \phi_2$ are metrics in $(\cH_2^*)^T$ then the radial functions of their asymptotic cone metrics satisfy $r_1^2 \geq r_2^2$.
\end{cor}

\subsection{Ambient metrics}

One of the main consequences of \cref{thm:testconfigurationscanbeembedded} for our purposes is that it will allow us to construct a well-behaved smooth subgeodesic ray associated to a given special test configuration $(\Pi: \cX \to \cY, \, T, \, \eta)$.

\begin{defn}\label{def:subgeodesic}
    Let $(\pi:X \to Y, \, \xi)$ be a polarized Fano fibration. Given a special test configuration $(\Pi: \cX \to \cY, \, T, \, \eta)$, a \demph{subgeodesic ray} is defined simply to be a $T' = T \times S^1_{\eta}$-invariant $\Phi \in \PSH(\cX_{\bbD}, \, -K_{\cX/\bbC})$.

 \end{defn}
Given a subgeodesic ray, one can use the isomorphism $\cX_{\bbD}^* \cong X \times \bbD^*$ to produce an $S^1$-invariant family $\phi_\tau \in \PSH(X, -K_X)$ for $|\tau| \in (0,1)$. One can more generally define subgeodesic rays in this fashion which a priori may not close up to the total space of some test configuration, but the more restricted setting is sufficient for our purposes.

Suppose now that we have an embedding $X \hookrightarrow \amb$ as in \cref{rmk:pffEmbed}, so that the vector field $\xi$ is induced from a Reeb field $\xi = \xi_1 + \xi_2 \in \ft_1 \oplus \ft_2$, on $\amb$. Our next goal is to construct a suitable fixed background metric $\omega_\mathrm{amb}$ on $\amb$ which will be useful throughout the rest of the paper.

\begin{lem}\label{thm:conicalsmoothing}
    Let $T_2 \subset \bbT^{N_2}$ be the torus generated by $\xi_2$ on $\bbC^{N_2}$. Then given a Kähler cone metric $\omega_{\xi_2}$ on $\bbC^{N_2}$, smooth away from the origin, with Reeb field $\xi_2$ and radial function $r: \bbC^{N_2} \to \bbR_+$,
    there exists a $T$-invariant %
    smooth function $h: \bbC^{N_2} \to \bbR_{\geq 0}$ such that $h \equiv r^2/2$ on $\{r > 4\}$ and $\omega_\mathrm{AC} = i \partial\bar{\partial} h$ is a Kähler metric on $\bbC^{N_2}$.
\end{lem}
\begin{proof}

    We set
    \[
        h = V(r^2)/2 + \frac1{C} \chi(r) |z|^2
        \p
    \]
    For an appropriate choice of convex increasing function $V$ (see \cite[Lem.~4.3]{VCcrepant}), $\tilde{\omega} \coloneq i \partial\bar{\partial} V(r^2)/2$ is a nonnegative closed $T$-invariant $(1,1)$-form that agrees with $\omega_{\xi_2}$ for $r > 2$.
    Then we pick $\chi$ to be a bump function supported on $[0, 3]$ and equal to $1$ on $[0, 2]$. Finally picking $C$ large we can ensure that
    \[
        \omega_\mathrm{AC} = i \partial\bar{\partial} h = \tilde{\omega} + \frac{1}{C} i \partial\bar{\partial} \left(\chi(r) |z|^2 / 2\right)
    \]
    is positive.
\end{proof}

Throughout the rest of the paper, we fix any such choice $\omega_{\xi_2}$, as the only key property that we need is the choice of Reeb field. Note that it is always possible to choose such a cone metric, see for example \cite{MSY}.

\newcommand{\psiamb}{\psi}
Given this, we can construct an ambient metric that will be suitable for our purposes:
\begin{defn} \label{def:psiamb}
    Given a smooth Kähler cone metric on $(\bbC^{N_2},\, \xi_2)$, we can construct an associated Kähler metric $\omega_\mathrm{amb}$ on $\amb$ by applying the lemma above, and setting
    \begin{equation}\label{eq:ambientmetric}
        \omega_\mathrm{amb} \coloneq \frac1{p} \omega_\mathrm{FS} + \omega_\mathrm{AC}
        \p
    \end{equation}
    Thus $p \omega_\mathrm{amb}$ is the curvature of the $T$-invariant hermitian metric $p \psi$  on $\cO(1)$ defined by
    \[
        e^{-p \psiamb} \coloneq e^{-p h} h_\mathrm{FS}
        \p
    \]
\end{defn}

The metric $\omega_\mathrm{amb}$ on $\amb$ is clearly weakly asymptotically conical, in the sense that it has a unique (Kähler) tangent cone at infinity, namely $(\bbC^{N_2}, \, \omega_{\xi_2})$.
\begin{notation}
    We will write $\lambda_t$ for the action of $\lambda: \mathbb{C}^* \to \bbP^{N_1} \times \bbC^{N_2} \times \bbC$ evaluated at $e^{-t}$, as well as the corresponding actions on $\PGL(N_1 + 1; \bbC)$ and $\GL(N_2; \bbC)$. We denote by $\eta = \frac{d}{dt}\lambda_t|_{S^1}$ the corresponding infinitesimal generator.
\end{notation}

\begin{defn}\label{def:psit}
    Pulling back by $\lambda_t$, the restrictions $e^{-\psiamb}|_{X_t}$ give us a family of hermitian metrics
    \[
        e^{-\psi_t} \coloneq \lambda_t^* e^{\psiamb} = \lambda_t^* \left(e^{-h} h_\mathrm{FS}^{1/p}|_{X_t}\right)
    \]
    on $-K_X \cong \frac{1}{p}\bigl(\cO(1)|_X\bigr)$.
\end{defn}

\newcommand{\proofstep}[1]{\noindent\textit{#1}}

\begin{prop} \label{thm:phitasymp}
    $\psi_t$ is a smooth family of metrics and for all $t \in (0, \infty)$ we have $\psi_t \in (\cH_2^*)^T$. In fact, for any compact interval $J \subseteq (0, \infty)$, all derivatives, in directions of $J$ or $M$, of $\psi_t$ will be uniformly $O(r^2)$ in equivariant charts.
\end{prop}

The proof of \cref{thm:phitasymp} is essentially an application of the following elementary fact.
\begin{lem} \label{thm:weightFrac}
    Let $K \subseteq \bbR^n$ be compact and let $f_i, g_i, i = 1, \dots, s$ be smooth functions on $(0, \infty) \times K$, which are weight vectors for $\frac{\partial}{\partial x^\alpha_1}$, i.e.
    \begin{equation} \label{eq:pdeweight}
        \frac{\partial f_i}{\partial x_1} = w_i f_i \quad\text{and}\quad \frac{\partial g_i}{\partial x_1} = w_i g_i \quad\text{for $i = 1, \dots s$}
    \end{equation}
    for some constants $w_i \in \bbR$. Further assume that $g_i > 0$. Then
    $ F \coloneq \dfrac{\sum_i f_i}{\sum_i g_i}$
    and all of its derivatives are bounded. In particular, all derivatives of
    $\log \sum_i g_i$
    are bounded.
\end{lem}

\begin{proof}
    We have $f_i = e^{\langle w_i, \, x \rangle } \tilde{f}_i$, where $\tilde{f}_i$ is a function only on $K$ and similarly for $g_i$. Thus by compactness of $K$, the functions $f_i/g_i = \tilde{f}_i / \mkern1mu \tilde{g}_i$ are all bounded by some constant $C$. Therefore $\left(\sum_i g_i\right) F = \sum_i f_i \leq C \sum_i g_i$, showing that $F \leq C$.

    Now applying the quotient rule we see that
    \[
        \partial_1 F = \frac{\sum_i \partial_1 f_i}{\sum_i g_i} - \frac{\sum_{ij} f_i \partial_1 g_j}{\sum_{ij} g_i g_j}
    \]
    is a difference of fractions satisfying the hypothesis of the theorem, since $\partial_1 f_i$ and $\partial_1 g_i$ still satisfy \cref{eq:pdeweight}, albeit for different weights. Then we induct on the order of the derivative.
\end{proof}
\begin{proof}[Proof of \cref{thm:phitasymp}]
    To see that $\psi_t$ and its derivatives are uniformly $O(r^2)$, it is sufficient to check that this is true for $\lambda_t r^2$ and for $\rho_t \coloneq \psi_t - \lambda_t^* h$, since the latter agrees with $\psi_t - \lambda_t^* r^2/2$ outside of a compact set.
    On an equivariant chart $U_\alpha \subseteq Y^\circ$, the map $X \hookrightarrow \bbP^{N_1}$ is given by holomorphic sections $f_0, \dots, f_{N_1}$ of $K_{U_\alpha} \cong \cO_{U_\alpha}$, which are weight vectors for $\xi$. Thus
    \[
        \rho_{t, \alpha} = \frac1{p} \log \sum_i |f_i^t|^2,
    \]
    where $f_i^t \coloneq \sum_j (\lambda_t)_{ij} f_j$ and $(\lambda_t)_{ij}$ are the matrix coefficients of the action of $\lambda_t$ on $\bbC^{N_1+1}$. Since $\lambda$ and $\xi$ commute, we have that $f_i^t$ and $|f^t_i|^2$ are also weight vectors for $\xi = \frac{\partial }{\partial x_1^\alpha}$.
    By the definition of equivariant charts (cf.\ \cref{def:equivchart}), the functions $f_i$ extend to $(0, \infty) \times [-1, 1]i \times \overline{U_\alpha'}$,
    thus we can apply \cref{thm:weightFrac} to conclude that all derivatives of $\rho_{t, \alpha}$ are bounded.

    On the other hand, on $U_\alpha$ we have $\lambda_t^* r^2 = e^{2 x_1^\alpha} g_t$ where $g_t$ does not depend on $x_1^\alpha$, i.e.\ has $\xi$-weight zero. Again, $g_t$ extends to a smooth function on $(0, \infty) \times [-1, 1]i \times \overline{U'_\alpha} \times J$ and therefore its derivatives are bounded. Thus all derivatives of $e^{2 x_1^\alpha} g_t$ are uniformly $O(e^{2 x_1^\alpha}) = O(r^2)$.

    We are left with showing that $\psi_t \in \cH^*_2$, which amounts to showing that $|i \partial\bar{\partial} \rho_t|_{i \partial\bar{\partial} r^2 / 2} = O(r^{-2})$. Since $i \partial\bar{\partial} r^2 / 2$ has $\xi$-weight 2 in equivariant charts, this is equivalent to saying that the coefficients of $i \partial\bar{\partial} \rho_t$ are $O(1)$ in equivariant charts, which was already established above.
\end{proof}

\subsection{Normalized hamiltonian potentials}\label{section:hamiltonain}
Let $(\pi:X \to Y, \, \xi)$ be a polarized Fano fibration with $X$ smooth, and suppose $X$ is equipped with a smooth $T$-invariant hermitian metric $e^{-\phi}$ on a line bundle $L \to X$, with positive curvature $\omega = i\partial \bp \phi > 0$. Since $H^1(X) = 0$ \cite{Wylie} (in fact, $X$ is simply connected, see \cite[Cor.~1.3]{Simplyconnected}, \cite[Prop.~3.10]{SunZhang}), there exists a moment map $\mu: X \to \ft^*$, which by definition satisfies
\begin{equation}\label{eq:momentmapdefinition}
    d\langle \mu , \, \eta \rangle = - \eta \iprod \omega ,
\end{equation}
for any $\eta \in \ft$. In particular, $\mu$ determines a choice of Hamiltonian potential $\theta_\eta$ for any $\eta \in \ft$ by
\begin{equation}\label{eq:hamiltonianpotentialmomentmap}
    \theta_\eta = \langle \mu , \, \eta \rangle.
\end{equation}
In general $\mu$ is unique only up to a translation in $\ft^*$, which affects each $\theta_\eta$ by the addition of a constant. If however we are given a specific choice of lift of the $T$-action to the total space of $L \to X$, then we can normalize our choice for $\mu$ by defining
\begin{equation}\label{eq:hamiltonianlift}
    \theta_\eta \coloneqq -\frac12 \frac{\cL_{-J\eta} e^{-\phi}}{e^{-\phi}}.
\end{equation}
It's straightforward to verify that $\theta_\eta$ is indeed an $\omega$-Hamiltonian potential.

\begin{defn}\label{thm:canonicalnormalizationdefinition}
    Let $(\pi:X \to Y, \, \xi)$ be a polarized Fano fibration with $X$ smooth, and $e^{-\phi}$ be a smooth $T$-invariant hermitian metric on $-K_X$ with positive curvature $\omega$. Then we say that a Hamiltonian potential $\theta_\eta$ for $\eta \in \ft$ is \demph{canonically normalized} if it satisfies \cref{eq:hamiltonianlift} for the canonical lift of the $T$-action to $-K_X$. In this case, this is equivalent to choosing $\theta_\eta$ such that
    \begin{equation}
        \Delta_{\omega} \theta_\eta + \theta_\eta + \frac{1}{2} J\eta \cdot\theta_\eta = 0,
    \end{equation}
    see \cite{CDS, SunZhang}.
\end{defn}

If $X$ or $e^{-\phi}$ are not smooth, then defining Hamiltonian potentials requires a bit of care. In this paper however we will only need to consider the case when $e^{-\phi}$ is induced from an embedding of $X \hookrightarrow \amb$, where $\amb$ is equipped with a smooth hermitian metric $e^{-\psi}$ on $\cO_{\amb}(1)$. In this case, we can simply define $\theta_\eta$ as a restriction:

\begin{rmk}\label{thm:canonicallynormalizedfromtheambientspace}
    Let $(\pi:X \to Y, \, \xi)$ be a (potentially singular) polarized Fano fibration with a $T$-equivariant embedding $\iota:X \hookrightarrow \amb$ such that $\iota^*\cO_{\amb}(1) = - pK_X$. In particular, there is a unique lift of the $T$-action to $\cO_{\amb}(1)$ making the induced map $-pK_{X} \to \cO_{\amb}(1)$ equivariant for the canonical lift of $T$ to $-pK_X$. Suppose that $\amb$ is equipped with a smooth hermitian metric $e^{-\psi}$ on $\cO_{\amb}(1)$ with positive curvature $p\omega$, and we set $\omega_X = \omega|_{X}$, which is the curvature of a hermitian metric $\left( e^{-\psi}|_X\right)^{\frac{1}{p}}$ on $-K_X$. Then for any $\eta \in \ft$, we can define an $\omega_X$-Hamiltonian potential $\theta_\eta$ by letting $p \theta_\eta$ be the Hamiltonian potential on $\amb$ defined by \cref{eq:hamiltonianlift} with respect to $e^{-\psi}$, and restricting to $X$. Then on the smooth locus of $X$ one can readily verify that the induced moment map satisfies \cref{eq:momentmapdefinition}. By a slight abuse of notation, we will also call $\theta_\eta$ the canonically normalized Hamiltonian potential in this situation.
\end{rmk}

\begin{lem} \label{thm:extendForms}
    Let $\theta_\xi \in C^\infty(\amb)$ be the $\omega_\mathrm{amb}$-hamiltonian obtained from the action of $\xi$ on $\cO_\amb(1)$ and let $\alpha \in \Omega^p(\amb)$ be any differential form that is a $(-J \xi)$-weight vector. Then the form $e^{-\theta_\xi} \alpha$ extends smoothly to $\bbP^{N_1} \times \bbP^{N_2}$.
    As a consequence, for any $\eta \in \ft$, the forms $e^{-\theta_\xi}\omega_\mathrm{amb}^n,\, \theta_\eta e^{-\theta_\xi}\omega_\mathrm{amb}^n,\, \eta \iprod e^{-\theta_\xi}\omega_\mathrm{amb}^n$ extend smoothly to $\bbP^{N_1} \times \bbP^{N_2}$.
\end{lem}

\begin{proof}
    Recall that $\xi = \xi_1 + \xi_2$ where $-J\xi_2 = \sum 2 a_i \Re \left(z_i \frac{\partial}{\partial z_i}\right)$ for some $a_i > 0$, and $z_i$ are euclidean coordinates on $\bbC^{N_2}$.
    Outside of a compact subset, the Hamiltonian $\theta_{\xi_2}$ is the sum of a constant and a function of $(-J \xi_2)$-weight 2, as can be seen from the equation $\d \theta_{\xi_2} = \pm \xi_2 \iprod \omega_\mathrm{amb} = \pm \xi_2 \iprod \omega_\mathrm{AC}$. Meanwhile $\theta_{\xi_1}$ (and its derivatives) are bounded. Thus for the Reeb vector field $\xi$ we obtain $\theta_\xi \geq \epsilon r^2$ for $r > 1$ and $\epsilon$ sufficiently small.

    Thus it suffices to show that, in the euclidean coordinates on $\bbC^{N_2}$, forms which are weight vectors for $\xi$ will have coefficients that have polynomial growth, and their derivatives also have polynomial growth.
    Indeed, all coordinate derivatives of the components of $e^{-\theta_\xi} \alpha$ will be of the form $e^{-\theta_\xi}$ times a function of polynomial growth. Since the transition functions of the standard charts of $\bbP^{N_2}$ are rational functions, we see that in all standard charts for $\bbP^{N_2}$ we have $(e^{-\theta_\xi} \alpha)_x \to 0$ as $x \to \bbP^{N_2} \setminus \bbC^{N_2}$, with derivatives. Thus extending $e^{-\theta_\xi} \alpha$ by zero on $\bbP^{N_2} \setminus \bbC^{N_2}$ is smooth.

    We have
    \begin{equation} \label{eq:coordVectWeight}
        \left[-J\xi_2, \, \frac{\partial}{\partial z_i} \right] = -a_i \frac{\partial }{\partial z_i}, \qquad \left[-J\xi_2, \, \frac{\partial}{\partial \bar{z}_i} \right] = -a_i \frac{\partial }{\partial \bar{z}_i}
        \p
    \end{equation}
    The coordinate components of $\alpha \in \Omega^p(\amb)$ can be obtained as $(T \iprod \alpha)|_{T \bbP^{N_1}}$ for some multivector $T = \partial_{z_{i_1}} \wedge \dots \wedge \partial_{z_{{i_k}}} \wedge \partial_{\bar{z}_{j_1}} \wedge \dots \wedge \partial_{\bar{z}_{j_{m-k}}}$. Here we think of the components as taking values in $\Omega^{p-m}(\bbP^{N_1})$. By \cref{eq:coordVectWeight}, $T$ is a weight vector for $\xi_2$, thus ${T \iprod \alpha}|_{T \bbP^{N_1}}$ is a weight vector for $\xi$. By \cref{eq:coordVectWeight} again, we see that the coordinate derivatives of this are also weight vectors.
    Finally, if $\beta: \bbC^{N_2} \to \Omega^{p-m}{\bbP^{N_1}}$ has $\xi$-weight $w$ then $r^{-w} \beta$ is $\xi$-invariant and thus $\cL_{\xi_2} r^{-w} \beta = - \cL_{\xi_1} r^{-w} \beta$. Since the latter is bounded on a compact set, we conclude that $\beta$ has polynomial growth.

    Then to see that, say, $e^{-\theta_\xi} \omega_\mathrm{amb}^n$ extends to $\bbP^{N_1} \times \bbP^{N_2}$ we decompose
    \[
        \omega_\mathrm{amb}^n = \sum_{k=0}^n \binom{n}{k} \omega_\mathrm{FS}^k \wedge \omega_\mathrm{AC}^{n-k}
    \]
    and notice that by what we just proved $e^{-\theta_\xi} \omega_\mathrm{AC}^{n-k}$ extends since $\omega_\mathrm{AC}^{n-k}$ has a fixed weight outside a compact set.
\end{proof}

\subsection{Geodesic rays}

\begin{defn}
    Let $(\pi:X \to Y, \, \xi)$ be a polarized Fano fibration with $\phi \in \overline{\cH}_2$, and $(\Pi:\cX \to \cY, \, T, \, \eta)$ be a special test configuration. Then a \demph{geodesic ray} emanating from $\phi$ is a subgeodesic ray $\Phi$ solving the homogeneous Monge-Amp\`ere equation with $\eta$-invariant boundary data:
    \begin{equation}\label{eq:geodesicray}
        \begin{cases}
             (i\partial\bp \Phi)^{n+1} = 0  & \textnormal{on } \cX_{\bbD} \\
              \Phi|_{X_1} = \phi
              \p
        \end{cases}
    \end{equation}
    Similarly, given $\phi_1, \phi_k \in \Hsing$ and an annulus $\Omega_k = \bbD\setminus\bbD_{e^{-k}}$, a \demph{geodesic segment} is an $S^1$-invariant $\Phi \in \PSH(X \times \Omega_k; -K_{X})$ satisfying
    \begin{equation}\label{eq:geodesicsegment}
        \begin{cases}
             (i\partial\bp \Phi)^{n+1} = 0  & \textnormal{on } X \times \Omega_k \\
              \Phi|_{X_1} = \phi_1  &  \\
              \Phi|_{X_{e^{-k}}} = \phi_k
              \p
        \end{cases}
    \end{equation}
\end{defn}

The key point of this section is to use the existence of geodesic segments in the asymptotically conical setting established in \cite[Thm.~2.9]{CarlosUniq} to prove existence for geodesic rays.

\begin{thm}[Existence of geodesic rays with smooth initial data]\label{thm:geodesicrays}
    Let $(\pi: X \to Y, \, \xi)$ be a polarized Fano fibration such that $X$ is smooth, and let $\phi \in (\cH_2^*)^T$. Given any special test configuration $(\Pi: \cX \to \cY, \, T, \, \eta)$, there exists a geodesic ray $\Phi$ emanating from $\phi$.
\end{thm}

In order to prove this, we will crucially make use of an ambient metric $\psiamb$ as constructed in \cref{def:psit}, which depends on a choice of initial cone metric $\omega_{\xi_2}$ on $(\bbC^{N_2}, \, \xi_2)$. For technical reasons, we need to adapt our choice of initial cone to the given $\phi$.

\begin{lem}\label{thm:referenceissmaller}
    Given any initial $\phi \in \mathcal{H}^*_2$, we can choose a reference $\psiamb$ on $\mathcal{O}_{\amb}(1)$ as in \cref{def:psit} with the property that $\phi \geq \psi_0$, where $\psi_0 = \psiamb|_{X_{e^{-0}}}$.
\end{lem}

\begin{proof}
    Let $\omega_\phi = i \partial \bp \phi$ be the AC metric on $X$ associated to $\phi$, and let $r_\phi$ be the radial function for the corresponding asymptotic cone on $Y$. Choose any cone metric $\tilde{\omega}_{\xi} = i \partial \bar{\partial} \tilde{r}^2/2$ on $(\bbC^{N_2}, \, \xi_2)$. Then $\tilde{r}|_Y$ is the radial function of the cone metric $\tilde{\omega}_\xi|_Y$ on $(Y, \, \xi_2)$, and so it is uniformly equivalent to $r_\phi$:
    \[
        \epsilon r_\phi \leq \tilde{r} \leq \epsilon^{-1} r_\phi
        \p
    \]
    Now we pick $r = \epsilon \tilde{r}/\sqrt2$, so that $r_\phi^2/2 \geq r^2$, and perform the construction from \cref{def:psit} to produce a metric $\tilde{\omega}_\mathrm{amb}$ on $\amb$.

    Now by \cref{thm:metricvsconeC0} and construction of $\psi$, we have
    \[
        \phi - \psi_0 = \frac{r_\phi^2}{2}  - \frac{r^2}{2} + O((\log r)^2) \geq \frac{r^2}{2} + O((\log r)^2)
        \p
    \]
    It follows that $\phi \geq \psi_0$ outside of a compact set $K \subset X$. Finally, we use our freedom to modify $\psi_0$ by adding a constant.
    In particular, if $C > 0$ is sufficiently large, then replacing $\psi_0 \to \psi_0 - C$, we can further ensure that
    $\inf_{K } (\phi - \psi_0) \geq 0$, so that $\phi \geq \psi_0$ globally.
\end{proof}

Fixing now such a choice of $\psiamb$, we can trivially extend this to a metric on $\cO_{\amb \times \bbC}(1)$ and therefore as a metric on $-p K_\cX$ by restricting to $\cX$. Let $\psi_t \coloneq (\lambda_t^* \psiamb)|_{X_0}$ be the associated subgeodesic ray as in \cref{def:psit}.

For any $k$ we can invoke \cite[Thm.~2.9]{CarlosUniq} to obtain a geodesic segment $\hat{\Phi}^k$ in $\Hsing$ connecting $\Phi_0 = \phi$ to $\psi_k$.
Using the isomorphism $\cX^* \cong X \times \bbC^*$ we can view this as a metric in $\PSH\left(\cX^*_k; -K_\cX\right)$. Since $\Phi_0 \geq \psiamb$, by the construction of $\hat{\Phi}^k$ in \cite{CarlosUniq} using Perron's method, we know that $\hat{\Phi}^k \geq \psiamb$ on $\cX^*_k$.
Then we define
\begin{equation}\label{eq:phik}
    \Phi^k \equiv \begin{cases}
        \hat{\Phi}^k & \text{on $\cX^*_k$} \\
        \psiamb & \text{on $\cX_\bbD \setminus \cX^*_k = \Pi^{-1}(\bar{B}_{e^{-k}})$}
        \p
    \end{cases}
\end{equation}
\begin{lem}
    The metrics $\Phi^k$ are a nondecreasing sequence in $\PSH(\cX_{\bbD}; -K_{\cX})$.
\end{lem}
\begin{proof}
    Observe that $\Phi_k = \lim_{\epsilon \searrow 0} u_\epsilon $, where
    \[
        u_\epsilon \coloneq \begin{cases}
            \max\{\psiamb + \epsilon,\, \hat{\Phi}^k\} & \text{on $\cX^*_k$} \\
            \psiamb + \epsilon & \text{on $\cX_\bbD \setminus \cX^*_k$}
        \p
        \end{cases}
    \]
    Since $\psiamb$ is continuous, $\hat{\Phi}^k$ is upper semicontinuous and both agree on $\Pi^{-1}(\partial \bbD_{e^{-k}})$, we have $\psiamb + \epsilon > \hat{\Phi}^k$ in a neighborhood of $\Pi^{-1}(\partial \bbD_{e^{-k}})$, showing that the metrics $u_\epsilon$, and therefore their decreasing limit are PSh.

    To see that $\Phi^k$ is nondecreasing in $k$, simply note that $\Phi^k|_{\cX^*_{k+1}}$ has the same boundary values as $\Phi^{k+1}|_{\cX^*_{k+1}}$. Thus $\Phi^k|_{\cX^*_{k+1}}$ is a competitor in the supremum in Perron's method for $\hat{\Phi}^{k+1}$, so $\hat{\Phi}^{k+1} \geq \Phi^k|_{\cX^*_{k+1}}$. This proves that $\Phi^{k+1} \geq \Phi^k$ on $\cX^*_{k+1}$, and both metrics agree on the rest of $\cX_\bbD$.
\end{proof}

\begin{prop} \label{thm:GeodesicRay}
    The increasing limit $\lim_{k \to \infty} \Phi^k$ exists, and its usc regularization defines a PSh metric $\Phi$ on $\cX$ with the property that $(i \partial\bar{\partial} \Phi)^{n+1} = 0$.
\end{prop}
\begin{proof}
    We first produce a locally uniform bound for $\Phi^k$ on $\cX^*$.
    Observe that
    \[
        \dot{\psi}_t = \frac{d}{d t} (\psi_t - \psi_0) = -\frac{d}{d t} \log \frac{e^{-\psi_t}}{e^{-\psi_0}} = e^{\psi_t} \frac{d}{d t} (\lambda_t^* e^{-\psi})|_X
        \p
    \]
    Moreover, recalling that $\eta$ denotes the infinitesimal generator of $\lambda$ then
    \[
        \theta_\eta = -\frac12 e^\psi \cL_{-J\eta} e^{-\psi} = -\frac12 e^\psi \frac{d}{d t}\bigg|_{t=0} \lambda_t^* e^{-\psi}
        \c
    \]
    so we have
    \begin{equation} \label{eq:psitdot}
        \dot{\psi}_t = -2 (\lambda_t^* \theta_\eta)|_X
        \p
    \end{equation}
    Since $\d \theta_\eta = - J \eta \iprod \omega_\mathrm{amb}$, we see that $\d \theta_\eta$ has $\xi$-weight 2, hence $\lvert \theta_\eta \rvert = O(r^2)$.
    Moreover, $\lambda^*_t r \leq e^{C t} r$ for some $C > 0$. Thus for $t \geq 1$ we have
    \begin{equation} \label{eq:psitdiff}
        \left\lvert\psi_t - \psi_0\right\rvert \leq 2 \int_0^t |\lambda_s^* \theta_\eta| \d s \leq C' \int_0^t \left( 1 + e^{2 C s} r^2 \right) \d s \leq C'(t + e^{2 C t} r^2)
    \end{equation}
    for some $C' > 0$.
    Since $\psi_0 - \Phi_0 = O(r^2)$, we observe that
    \begin{equation} \label{eq:PhikGrowth}
        \left\lvert \hat{\Phi}^k_k - \hat{\Phi}^k_0 \right\rvert = \left\lvert \psi_k - \Phi_0 \right\rvert \leq O(e^{2 C t} r^2)
        \c
    \end{equation}
    uniformly in $k$.

    For any $p \in X$ we consider the (immersed) complex curve given by
    \begin{align}
        c: \bbC &\to X  \times \bbD \\
        \tau &\mapsto (\gamma_{-C \tau}(p), e^\tau)
    \end{align}
    \c
    where $C$ is the same constant as above. By working in an equivariant chart $U$ on $X$, we can view $\hat{\Phi}^k$ as a function on $U \times \Omega_k$. Then $\hat{\Phi}^k \circ c$ is a PSh function on $\bbC$ which is independent of the imaginary direction because $\hat{\Phi}^k$ is invariant in $\xi$ and $\eta$.
    Thus $\Phi \circ c$ is convex on the real line $\bbR \subset \bbC$, so for $t \in \bbR$
    \begin{equation} \label{eq:convex}
        \hat{\Phi}^k_t(\gamma_{-C t}(p)) \leq \frac{t}{k} \hat{\Phi}^k_k(\gamma_{-C k}(p)) + \frac{k - t}{k} \hat{\Phi}^k_0(p)
        \p
    \end{equation}
    By the definition of $\cH_2^*$, we know that $\hat{\Phi}^k_0 = \phi = O(r^2)$.
    Furthermore equation \cref{eq:PhikGrowth} implies that
    \[
        \hat{\Phi}^k_k(\gamma_{-C k}(p)) - \hat{\Phi}^k_0(\gamma_{-C k}(p)) = O\left(e^{2 C k} r^2(\gamma_{-C k}(p))\right) = O(r^2(p))
    \]
    uniformly in $k$ and locally uniformly in $p$.
    Thus we can conclude from \cref{eq:convex} that
    \[
        \hat{\Phi}^k_t(\gamma_{-C t}(p)) = O\left(r^2(p)\right) \c
    \]
    In other words, writing $q = \gamma_{-C t}(p)$, we see that
    \begin{equation} \label{eq:PhiUnifBound}
        \hat{\Phi}^k_t(q) \leq O(r^2(\gamma_{C t}(p))) = O(e^{2 C t} r^2) \p
    \end{equation}
    so we see the sequence $\Phi^k$ on $\cX^*$ is locally uniformly bounded. Then \cref{thm:corNotHarnack} implies that the sequence is locally uniformly bounded on all of $\cX$.
    Since $\Phi^k$ is an increasing sequence, the limit $\Phi$ exists and is PSh by \cite[I.4.2]{agbook}.

    To see that the limit $\Phi$ satisfies the homogeneous \MA/ equation, we first observe that $(i \partial\bar{\partial} \Phi|_{\cX^*})^{n+1} = 0$. This follows from the continuity of the Bedford--Taylor \MA/ operator along increasing sequences \cite[Thm.~3.23]{GZ}, since $(i \partial \bar{\partial} \Phi^k)^{n+1} \equiv 0$ on $\cX^*_k$ by construction. Then since the \MA/-measure $(i \partial\bar{\partial} \Phi)^{n+1}$ is non-pluripolar, this is in fact sufficient to conclude that $(i \partial\bar{\partial} \Phi)^{n+1} = 0$ on all of $\cX$.
\end{proof}
Later on, we will need a quantitative form of the locally uniform bound on $\Phi$ that we obtained in the proof above:
\begin{lem}\label{thm:PhiBound}
    Using the same notation as above, there exists a $\kappa > 0$ uniform in $t = -\log|\tau|$ such that
    \[\Phi - \psi \leq |\tau|^\kappa r^2 + O(\log|\tau|) \p\]
\end{lem}

\begin{proof}
    Taking the limit as $k \to \infty$ of \cref{eq:convex} we obtain
    \begin{equation} \label{eq:convex2}
        \hat{\Phi}_t \leq \lim_{k \to \infty} \frac{t}{k} \gamma_{-C(k-t)}^* \psi_k + \frac{k-t}{k} \gamma_{C t}^*\Phi_0
        \p
    \end{equation}
    Recall that we have $\psi_k = \lambda_k^* \psi_0 = \lambda_k^* \left(\frac1p \phi_\mathrm{FS} + h\right)$. With a computation similar to \cref{eq:psitdiff}, but simpler one can check that
    \begin{equation}\label{eq:FSpullbacklinear}
        |\gamma_s^* \phi_\mathrm{FS} - \phi_\mathrm{FS}| \leq C_1 s
        \qtext{and}
        |\lambda_s^* \phi_\mathrm{FS} - \phi_\mathrm{FS}| \leq C_2 s
        \p
    \end{equation}
    Hence
    \[
        \frac{t}{k} \gamma_{-C(k-t)}^* \lambda_k^*\phi_\mathrm{FS} \leq \frac{t}{k} \phi_\mathrm{FS} + C_1t + C_2 Ct \left|1-\frac{t}{k}\right| \leq \frac{t}{k} \phi_\mathrm{FS} + C_3 t
        \c
    \]
    where we interpret any inequalities of this form as a statement about $\frac{t}{k}\left(\gamma_{-C(k-t)}^* \lambda_k^*\phi_\mathrm{FS} - \phi_\mathrm{FS}\right) \in C^\infty(\bbP^{N_1})$. On the other hand, since $\frac{r^2}{2} - C_2\leq h \leq \frac{r^2}{2} + C_2$, we have for every $t$ (pointwise)
    \[
        \lim_{k \to \infty} \frac{t}{k} h \circ \lambda_k \circ \gamma_{-C(k - t)} = \lim_{k \to \infty} \frac{t}{k} r^2 \circ \lambda_k \circ \gamma_{-C(k - t)} \leq \lim_{k \to \infty} \frac{t}{k} e^{2C k} e^{-2C(k-t)} r^2 = 0
        \p
    \]
    Putting these things together, we obtain from \cref{eq:convex2}
    \begin{align} \label{eq:mainest}
        \hat{\Phi}_t &\leq \lim_{k \to \infty} \left[ \frac{t}{k} \phi_\mathrm{FS} + C_3 t + \frac{k - t}{k} \gamma_{Ct}^* \Phi_0 \right]
            = C_3 t + \gamma_{Ct}^* \Phi_0 \\
            &\leq \gamma_{Ct}^* \psi_0 + O(e^{2 C t} r^2) + C_3 t
    \end{align}
    where we have used that $\Phi_0 - \psi_0 = O(r^2)$.
    Now similar to the computations above, we have
    \[
        \gamma_{-Ct}^* \psi_t = \gamma_{-Ct}^* \lambda_t^* \phi_\mathrm{FS} + h \circ \lambda_t \circ \gamma_{-Ct} = \phi_\mathrm{FS} + O(t) + O(r^2) = \psi_0 + O(t) + O(r^2) \p
    \]
    Pulling back both sides by $\gamma_{Ct}^*$, we obtain $\gamma_{Ct}^*\psi_0 = \psi_t + O(t) + O(e^{2 C t} r^2)$. Thus from \cref{eq:mainest} we can obtain
    \[
        \hat{\Phi}_t \leq \psi_t + O(t) + O\left(e^{4 C t} \lambda_t^* r^2\right) \p
    \]
    Pushing forward by $\lambda_t$ we can rephrase this as
    \begin{equation}
        \Phi \leq \psiamb + O(e^{4 C t}) r^2 + O(t)
    \end{equation}
    on $\cX^*$, from which the conclusion follows.
\end{proof}

\section{The Ding functional}
\label{sec:energyfunctionals}

Throughout this section, we restrict to the setting relevant for asymptotically conical Kähler shrinkers. Complex geometrically, this means that we fix a polarized Fano fibration $(\pi:X \to Y, \, \xi)$ such that $\dim_{\bbC}Y = \dim_{\bbC}X$ and $\pi:X \to Y$ is a resolution.

Fix a reference $T$-invariant $\phi_*\in \PSH(X; -K_X)$ with curvature $\omega_{*} = i \partial \bp \phi_*$.
Given another $T$-invariant smooth hermitian metric $\phi \in \PSH(X; -K_X)$ with curvature $\omega_\phi = i\partial \bp \phi$, we can define the \demph{energy} $\cE_\xi(\phi)$ by
\begin{equation}\label{eq:energysmooth}
    \cE_\xi(\phi) = \cE_{\xi}(\phi, \, \phi_*) = \int_0^1\int_X (\phi - \phi_*) e^{-\theta_{\xi, t}} \, \omega_t^n \wedge dt,
\end{equation}
where $\omega_t \coloneqq t\omega_\phi + (1-t)\omega_*$, and $\theta_{\xi, \, t}$ is the canonically normalized $\omega_t$-hamiltonian for $\xi$. A hermitian metric $\phi \in \PSH(X; -K_X)$ also naturally defines a measure on $X$. Indeed, in local coordinates $U \subset X_\mathrm{reg}$, we can write
\[e^{-\phi} = e^{-\phi_U} \, \Omega_z \otimes \overline{\Omega}_z, \]
where $\phi_U\in \PSH(U)$ and $\Omega_z = dz^1 \wedge \dots \wedge dz^n$ is the standard holomorphic volume form in $U$. We have therefore an associated volume form $dV_\phi$ on $X_\mathrm{reg}$, given locally by:
\begin{equation}\label{eq:dV}
    dV_{\phi} \underset{\mathrm{loc}}{=} e^{-\phi_U}\,  \Omega_z \wedge \overline{\Omega}_z \p
\end{equation}
Using this we can define another functional of anticanonical hermitian metrics:
\begin{equation}\label{eq:Lfunctional}
    \mathcal{L}(\phi) = -\log\left( \int_{X_\mathrm{reg}} \!\!dV_{\phi} \right) \p
\end{equation}
The condition that $X$ is klt is equivalent to the local integrability of $dV_{\phi}$ near the singularities of $X$.
Moreover, if $X$ is smooth it was shown in \cite{CarlosUniq} that $\cE_\xi(\phi)$ and $\mathcal{L}(\phi)$ are well-defined as long as $e^{-\phi}$ is AC (see \cref{def:AC}). In this paper, we will in fact only need to deal with the case where $X$ is smooth. Given these, we define the \demph{Ding functional} on smooth AC metrics by:
\begin{equation}\label{eq:ding}
    \cD_\xi(\phi) = \cD_{\xi}(\phi, \, \phi_*) = -\frac{1}{(2\pi)^n n! \bbW^{\mathrm{an}}(\xi)} \cE_{\xi}(\phi, \, \phi_*) + \cL(\phi).
\end{equation}
It is well-known that Kähler--Ricci shrinkers are formally the critical points of $\cD_\xi$. Indeed, if $X$ is smooth then a Kähler metric $\omega_\mathrm{KRS}$ on $X$ satisfies \eqref{eq:shrinker} with vector field $\xi$ if and only if a suitably normalized $\phi_{\rm KRS} \in \PSH(X; -K_X)$ satisfies 
\begin{equation}\label{eq:shrinker2}
    e^{-f}\omega_\mathrm{KRS}^n = dV_{\phi_{\rm KRS}},
\end{equation}
where $f = \theta_{\xi}$ is the canonically normalized $\omega_\mathrm{KRS}$-Hamiltonian potential for $\xi$. It immediately follows that, formally, $\delta \cD_{\xi}(\phi_{\rm KRS}) = 0$.

We recall some key basic properties of $\cD_{\xi}$ established in \cite{CarlosUniq}:
\begin{lem}[{\cite[Prop.~3.27, Lem.~4.7]{CarlosUniq}}]
    \label{thm:Dingprops}
    The functionals $\cE_\xi, \cL$, and $\cD$ extend in a well-defined way to all of $\Hsing$ only taking finite values. Moreover
    \begin{enumerate}
        \item \label{thm:Dingprops:convex}
            $\cL$ and $\cE_\xi$ are convex along subgeodesics. Hence along a subgeodesic $\Phi$, $\cD_{\xi}$ has a well-defined slope at infinity.
        \item \label{thm:Dingprops:affine}
            $\cE_\xi$ is affine along geodesics, and therefore $\cD_\xi$ is convex along geodesics.
        \item \label{thm:Dingprops:crit}
            If $\Phi$ is a subgeodesic starting from a solution $\phi_\mathrm{KRS}$ to the Kähler--Ricci shrinker equation \eqref{eq:shrinker2} then $\frac{d^+\!\!}{d t} \cD_{\xi}(\Phi_t) \geq 0$.
    \end{enumerate}
\end{lem}

\begin{cor} \label{thm:DingBounded}
    If $X$ admits a Kähler--Ricci shrinker $\phi_\mathrm{KRS}$, then $\cD_\xi$ is bounded below on $\cH^*_2$ by $\cD_\xi(\phi_\mathrm{KRS})$.
\end{cor}
\begin{proof}
    Choose an arbitrary $\phi \in (\cH^*_2)^T$. By \cite[Thm.~2.9]{CarlosUniq} there exists a geodesic segment $(\Phi_t)_{t \in [0, 1]}$ in $\Hsing$ connecting $\phi_\mathrm{KRS}$ to $\phi$. By \cref{thm:Dingprops} $\cD_\xi(\Phi_t)$ will be monotone increasing in $t$.%
\end{proof}

\subsection{The slope at infinity}

Let $(\Pi: \cX \to \cY, \, T, \, \eta)$ be a special test configuration together with an embedding in $\amb$ as in \cref{thm:testconfigurationscanbeembedded}. Let $\psiamb$ be a metric on $\amb$ with corresponding subgeodesic ray $\psi_t$ of metrics on $X$ obtained as in \cref{def:psit}. By \cref{thm:Dingprops} we have well-defined slopes at infinity for $\cE_\xi(\psi_t)$ and $\cL(\psi_t)$:

\begin{prop} \label{thm:Slopes}
Along the subgeodesic ray $\psi_t$, the slopes at infinity of $\cE_{\xi}$ and $\cL$ are given by:
    \begin{align}
        \lim_{t \to \infty} \frac{\d^\pm}{\d t} \cE_\xi(\psi_t, \, \psi_0) &= -2\int_{X_0} \theta_\eta e^{-\theta_\xi} \omega_\mathrm{amb}^n \\
        \lim_{t \to \infty} \frac{\d^\pm}{\d t} \cL(\psi_t) &= 0
    \end{align}
    where $\theta_\xi, \, \theta_\eta \in C^\infty(\amb \times \bbC)$ are the normalized $\omega_\mathrm{amb}$-Hamiltonians for $\xi$ and $\eta$ respectively.
\end{prop}
To prove this, we need the following result, which can be readily deduced from \cite[Thm.~4.7]{FedererVarieties} (see also \cite{KingVarieties}):
\begin{lem} \label{thm:wcont}
    Let $\cX \subseteq \bbP^{N_1} \times \bbC^{N_2} \times \bbC$ be an $(n+1)$-dimensional algebraic variety and let $\alpha \in \Omega^n(\amb \times \bbC)$ be an $n$-form that extends smoothly to $\amb \times \bbC$. Then the function
    \[
        \tau \mapsto \int_{X_\tau} \alpha
    \]
    is continuous in $\tau$.
\end{lem}
Here $X_\tau$ is the scheme-theoretic fiber of the projection map $\cX \to \bbC$, so the integral over $X_\tau$ is computed as an integral over $(X_\tau)_\mathrm{reg}$, accounting for the multiplicity of each irreducible component.
\begin{proof}[Proof of \cref{thm:Slopes}]
    To compute the slope at infinity of $\cE_\xi$, recall from \cref{eq:psitdot} that
    \[ \dot{\psi}_t = -2 (\lambda_t^* \theta_\eta)|_X
        \p\]

    Since $\psi_t \in \cH_2^*$, we can apply \cite[Prop.~3.26]{CarlosUniq} to obtain
    \begin{align}
        \frac{d}{d t} \cE_\xi(\psi_t, \, \psi_0) &= \int_X \dot{\psi}_t e^{-\theta_\xi} \omega_\mathrm{amb}^n = -2\ \int_X \lambda_t^* \left[ \theta_\eta e^{-\theta_\xi} \omega_\mathrm{amb}^n \right] \\
        &= -2 \int_{X_{e^{-t}}} \theta_\eta  e^{-\theta_\xi} \omega_\mathrm{amb}^n = -2\langle [X_{e^{-t}}], \, \theta_\eta  e^{-\theta_\xi} \omega_\mathrm{amb}^n \rangle
        \p
    \end{align}
    Since $\theta_\eta  e^{-\theta_\xi} \omega_\mathrm{amb}^n$ extends smoothly to $\bbP^{N_1} \times \bbP^{N_2}$ by \cref{thm:extendForms}, we can apply \cref{thm:wcont} to take the limit $t \to \infty$ and obtain
    \[\lim_{t \to \infty} \frac{d}{d t} \cE_\xi(\psi_t, \, \psi_0) = \lim_{\tau \to 0}-2\langle [X_{\tau}], \, \theta_\eta  e^{-\theta_\xi} \omega_\mathrm{amb}^n \rangle = -2\langle [X_{0}], \, \theta_\eta  e^{-\theta_\xi} \omega_\mathrm{amb}^n \rangle,\]
    as desired.

    Since $\cL(t) \coloneqq \cL(\psi_t)$ is convex by \cref{thm:Dingprops}, we can compute its slope $\cL^\mathrm{NA}$ at infinity by

         \[  \cL^\mathrm{NA} \coloneqq \inf \left\{ 2\ell \phantom{\big|} \middle|~ \cL(t) - 2 \ell t \leq C \right\} \p\]
    As in \cite[Prop.~3.8]{BermanQFano}, we can identify this with
    \begin{equation}
        \inf \left\{2 \ell ~\middle|~ \int_{\bbD^*} e^{(\ell - 1) \log |\tau|^2 - \cL(-\log |\tau|)} i \d \tau \wedge \d \bar{\tau} < \infty \right\} \p
    \end{equation}

    Indeed if $2\ell > \cL^\mathrm{NA}$, then using convexity again we can see that $\cL(t) -2(\ell-1)t \leq (2 -\varepsilon)t + C$, from which we can see that
    \begin{equation}\label{eq:inf2}
     \int_{\bbD^*} e^{(\ell - 1) \log |\tau|^2 - \cL(-\log |\tau|)} i \d \tau \wedge \d \bar{\tau} \leq C \int_{\bbD^*}|\tau|^{2-\varepsilon}i d\tau \wedge d\bar{\tau} < \infty \p
     \end{equation}
    The same argument shows that the infimum $2\ell = \cL^\mathrm{NA}$ is precisely the value where \cref{eq:inf2} blows up.

    Recall that, by the definition \cref{eq:Lfunctional} of $\cL$, we have for any $t \in \bbR$
    \[
        \cL(\psi_t) = -\log \left(\int_{X}  dV_{\psi_t} \right)
            = - \log \left(\int_{X_{e^{-t}}} \!\! (\lambda_t)_* \d V_{\psi_t} \right)
        \c
    \]
    where $dV_{\psi_t}$ is the associated volume form \cref{eq:dV} on $X = X_\mathrm{reg}$ to $\psi_t \in \PSH(X;-K_X)$.
    Analogous to \cref{eq:dV}, we can associate to $\psi \in \PSH(\cX; - K_{\cX / \bbC})$ an $(n, n)$-form $\d V_\psi$ on $\cX_\mathrm{reg}$.
     Since $\lambda_t$ lifts to $-K_{\cX/\bbC}$ canonically, we have
     \[ dV_{\lambda_t^*\psi} = \lambda_t^*dV_{\psi} \c\]
     for any $t$.
     Hence, for any $t = - \log|\tau|$, we have
    \[
        \cL(\psi_t) = -\log \left(\int_{X_{\tau}} \!\!dV_\psi \right) \p
    \]
     We then observe that
    \begin{align}
        \int_{\bbD^*} e^{(\ell - 1) \log |\tau|^2 -\cL(-\log |\tau|)} i \d \tau \wedge \d \bar{\tau} &= \int_{\bbD^*} |\tau|^{2(\ell -1)} \left(\int_{X_\tau}  dV_{\psi}\right) i \d \tau \wedge \d \bar{\tau} \\
                &= \int_{\cX_\bbD^*} |\tau|^{2(\ell -1)} \,  dV_{\psi} \wedge i \Pi^* (\d \tau \wedge \d \bar{\tau}) \\
                &= \int_{(\cX_\bbD)_\mathrm{reg}} \!\!|\tau|^{2(\ell -1)} \,  dV_{\psi} \wedge i \Pi^* (\d \tau \wedge \d \bar{\tau})
                \label{eq:int_thresh}
        \p
    \end{align}

    Now $\mu \coloneqq dV_{\psi} \wedge i \Pi^* (\d \tau \wedge \d \bar{\tau})$ is just the smooth volume form associated to the image of $e^{-\psi}$ in $\PSH(\cX; -K_{\cX}) \cap C^\infty$ under the isomorphism $-K_{\cX/\bbC} \cong -K_{\cX}$. Then again just as in \cite[Prop.~3.8]{BermanQFano}, we can identify $1 - \frac{1}{2}\cL^\mathrm{NA}$ with the log canonical threshold of $(\cX, \, 0, \, X_0)$ near $X_0$, which vanishes since $X_0$ is klt.

    Indeed,  \cite[Thm.~7.5]{KollarPairs} implies that the pair $(\cX,\, X_0)$ is purely log terminal,
    and hence there exists a log resolution $\sigma: \tilde{\cX} \to \cX$ of $(\cX, \, X_0)$ with snc exceptional divisor $X_0' + \sum_i E_i$ such that
    \[
        \sigma^*(K_\cX + X_0) = K_{\tilde{\cX}} + X_0' - \sum_i a_i E_i \c
    \]
    where $X_0'$ is the proper transform of $X_0$ and $a_i > -1$. We can write $\sigma^* X_0 = X'_0 + \sum_i b_i E_i$ for some $b_i \geq 0$, so
    \begin{equation} \label{eq:discrep}
        \sigma^*(K_\cX + (1 - \ell) X_0) = K_{\tilde{\cX}} + (1- \ell) X_0' - \sum_i (a_i + \ell b_i) E_i
        \p
    \end{equation}
    Since $\tau^{\ell -1}$ has a pole of order $\ell-1$ along $X_0$, it follows that $\sigma^*(|\tau|^{2 (\ell - 1)} \mu)$ is a volume form on $\tilde{\cX}$ with poles along $X_0'$ and $E_i$. Thus the integral \cref{eq:int_thresh} is locally finite if and only if $\ell - 1 > -1$ and $a_i + \ell b_i > -1$ (see \cite[Lem.~3.7, Prop.~3.8]{BermanQFano}, and the references therein). Since $a_i > -1$, both conditions are clearly satisfied for $\ell > 0$, and the first holds if and only if $\ell > 0$. Moreover, local integrability of $|\tau|^{2(\ell -1)}\mu$ on $\cX$ is actually sufficient in this case. To see this, notice that since $(\cX)_\mathrm{sing} = (X_0)_{\rm sing}$ is $T'$-invariant and that $D =  X_0' + \sum_i E_i$ is snc, we can cover $\tilde{\cX}_\bbD = \sigma^{-1}(\cX_\bbD)$ by charts $U_\alpha = (z_1, \dots, z_{N_1 + N_2 +2})$ which are equivariant for the whole $T'$-action, and such that $ D \cap U_\alpha$ is contained in the union of the coordinate hyperplanes. Then a direct calculation shows that
    \[ \int_{U_{\alpha \cap \{ r \geq R\}}} \!\!\!\! \sigma^*\left(|\tau|^{2(\ell -1)} \mu \right) = O\left(R^{k}e^{-R^2} \right) \c \]
    for some $k > 0$, since $\psi = O(r^2)$ in these charts.
\end{proof}

\subsection{The Futaki invariant}

The goal of this section is to prove:
\begin{thm} \label{thm:volumeAgrees}
    Let $(\pi: X \to Y,  \, \xi)$ be a polarized Fano fibration, embedded in $\amb$ as in \cref{thm:testconfigurationscanbeembedded} and let $\omega = i \partial\bar{\partial} \psiamb$ be the restriction of an ambient metric as in \cref{def:psiamb}. Let $\theta_\xi$ be the canonically normalized $\omega$-Hamiltonian potential of $\xi$. Then
    \begin{equation} \label{eq:volumeAgrees}
        \bbW_X(\xi) = \bbW_X^\mathrm{an}(\xi) \coloneqq \frac{1}{(2\pi)^n}\int_X e^{-\theta_\xi} \frac{\omega^n}{n!}
        \p
    \end{equation}
    As a consequence, we have that for all $\eta \in \ft$,
    \begin{equation} \label{eq:futakiAgrees}
        \Fut_\xi(X, \eta) = \frac{1}{(2\pi)^n}\int_X \theta_\eta \, e^{-\theta_\xi} \frac{\omega^n}{n!}.
    \end{equation}
\end{thm}

To prove this, we use an argument originally due to Collins--Székelyhidi \cite{ColSze2} in the affine case. The idea is to show that one can choose a well-behaved degeneration of $X$ to a union $X_0$ of linear subspaces in $\amb$ with multiplicity, and moreover which preserves the weighted volume functional. On $X_0$ one can argue directly to see the equality, which yields the result.

In the course of the proof, it will be useful to deal with weighted volumes of more general $S$-modules:
\begin{defn} \label{def:Wmodule}
    Let $M$ be a finitely generated graded $S$-module equipped with a $T$-action that is compatible with the module structure, so there is a decomposition $M = \bigoplus_{m, \alpha} M_{m, \alpha}$ by degree and weight.
    We define the volume functional of $M$ as
    \begin{equation} \label{eq:volModule}
        \bbW_M(\xi) = \lim_{m \to \infty} \frac{1}{m^n} \sum_{\alpha \in \ft^*} e^{-\langle \frac{\alpha}{m}, \, \xi \rangle} \dim_\bbC M_{m, \alpha},
    \end{equation}
    when the limit exists.
\end{defn}

\begin{notation}
    For a graded module $M$, we write $M^{(p)}$ for the same moduel with a "dilated" grading
    \[
        (M^{(p)})_k = \begin{cases}
            M_{k/p} & p | k \\
            0 & \text{otherwise}
        \end{cases}
    \]
    In the case of $S$-modules, which might be (multi-)graded for several of the different gradings on $S$, $M^{(p)}$ only means we dilate the "main" grading of $S$, i.e.\ the grading by $x$-degree.
\end{notation}

\begin{lem}
    Let $(X, -p K_X) \inj (\amb, \cO(1))$ be a $T$-equivariant embedding of a polarized Fano fibration
    and let $I \ideal S$ be the ideal corresponding to $X$. Then $\bbW_X = \bbW_{(S / I)^{(p)}}$, i.e.
    \begin{equation} \label{eq:volumeI}
        \bbW_X(\xi) = \lim_{m \to \infty} \frac{1}{(p m)^n}\sum_{\alpha \in \ft^*} e^{- \langle \frac{\alpha}{p m}, \, \xi \rangle} \left( \dim S_{m, \alpha} - \dim I_{m, \alpha} \right)
        \p
    \end{equation}
\end{lem}
\begin{proof}
    By definition, $I_m$ is the kernel of
    \[
        S_m = H^0\left(\amb, \cO(m)\right) \to H^0(X, \cO(m)|_X) \cong H^0(X, -p m K_X)
        \p
    \]
    Thus the claim will follow once we know that this map is surjective for $m \gg 0$.
    To that end, consider the long exact sequence in cohomology of
    \[
        0 \to \scI_X \otimes \cO(m) \to \cO(m) \to i_* \cO(m) \to 0
    \]
    where $i: X \to \amb$ is the inclusion. We get
    \[
        H^0\left(\amb, \cO(m)\right) \to H^0(X, \cO(m)|_X) \to H^1(\amb, \scI_X \otimes \cO(m))
    \]
    but the rightmost term vanishes for $m \gg 0$ by Serre vanishing for $\bbP^{N_1}_{\bbC^{N_2}}$ \cite[III.5.3]{Hart}.
\end{proof}

A finitely generated $S$-module $M$ with a $T$-action corresponds to a $T$-equivariant sheaf $\cF$ on $\amb$, or equivalently, a $\bbC^* \times T$-equivariant sheaf $\hat{\cF}$ on $\bbC^{N_1 + 1} \times \bbC^{N_2}$. A useful property of $\bbW_M$ is that it doesn't see changes to $\cF$ on sets of dimension $\leq n-1$, or equivalently, changes to $\hat{\cF}$ on sets of dimension $\leq n$:
\begin{lem} \label{thm:FutakiLowerDim}
    Let $M$ be as in \cref{def:Wmodule}. Then $\bbW_M(\xi) = 0$ for all Reeb vector fields $\xi$ if $\dim \supp M \leq n$.
\end{lem}
\begin{proof}
    Recall that $\supp M \coloneq \{\fp \in \Spec S \mid M_\fp \neq 0\}$.
    For this, we need to work with an additional grading on $S$ coming from the affine structure on $\bbC^{N_2}$. In particular, we can view $S$ as a bi-graded ring by setting
    \[
        S_{m, k} = \bbC[x_0, \dots, x_{N_1}]_m \otimes \bbC[y_1, \dots, y_{N_2}]_k
        \p
    \]
    Consequently we can filter $S$ as $F_k S_m \coloneq \bigoplus_{j \leq k} S_{m, j}$.
    After picking a set of homogeneous generators $s_i \in M_{m_i}, i = 1, \dots, r$ for $M$, this induces a filtration $M_m = \bigcup_k F_k M_m$ on every $M_m$ by
    \[
        F_k M_m = \sum_i F_k S_{m - m_i} \cdot s_i
        \p
    \]
    \begin{claim} \label{thm:BiHilbertPoly}
        There exists an $L > 0$ such that for $m \geq L$ and $k \geq L$ the function $P: \bbN \times \bbN \to \bbN$ given by
        \[
            P(k, m) \coloneq \dim ((\gr M)_{m, k}) = \dim F_k M_m - \dim F_{k-1} M_m
        \]
        is a polynomial in $k$ and $m$ of total degree (in $m$ and $k$ combined) at most $n-2$.
    \end{claim}
    \begin{pf}
        Since $\gr M$ is bigraded and $S$ is generated by elements of bidegrees $(1, 0)$ and $(0, 1)$, we can deduce from \cite[Cor.~5.8.19]{KreuzerRobbiano} that there is a $L \geq 0$ such that $P(m, k) \coloneq \dim_\bbC (\gr M)_{m, k}$ is a polynomial for $m, k \geq L$, which we will denote by $P^*$. We may assume $P^* \neq 0$ since otherwise our claim is trivial.
        On the other hand, $\supp M = \cV(J) \subseteq \Spec S$ where $J = \operatorname{Ann}_S(M) \ideal S$, and $M$ is a finitely generated $S/J$-module. Since $\dim S/J = \dim \supp M \leq n$ the classical result on the degree of the %
        usual Hilbert polynomial of an affine algebra
         \cite[Thm.~4.6.36]{KreuzerRobbiano} shows
         that
        \begin{equation} \label{eq:SP}
            S_P(t) \coloneq \sum_{m+k \leq t} P(m, k)
        \end{equation}
        is a polynomial of degree $\leq n$ for $t \gg 0$.
        We can decompose $P^*$ as $P^* = Q + O(m^{s - 1} + k^{s-1})$ where $s$ is the total degree of $P^*$ and $Q \neq 0$ is the polynomial consisting of the terms in $P^*$ of total degree $s$. %
        Since $P^* \geq 0$ for all $(m, k)$ with $m, k \geq L$, it's easy to see that $Q(x, y) \geq 0$ for all $x, y > 0$.
        Fix $t > L$ and observe that
        \begin{align}
            S_P(2^a t) &\geq \sum_{\substack{m + k \leq 2^a t \\ m, k \geq L}} \big[ Q(m, k) + O(m^{s-1} + k^{s-1}) \big]  \\
                       &= O(2^{(s+1) a}) + \sum_{\substack{m + k \leq 2^a t \\ m, k \geq L}} Q(m, k)
        \end{align}
        We want to show that the sum involving $Q$ is bounded below by $\epsilon \cdot 2^{(s+2)a}$. Since $S_P$ is degree $n$, this will show our claim that $n \geq s+2$.
        As $a \to \infty$ we claim
        \[
            2^{-(s+2) a} \sum_{\substack{m + k \leq 2^a t \\ m, k \geq L}} Q(m, k) = \sum_{\substack{m + k \leq 2^a t \\ m, k \geq L}} Q(2^{-a} m, 2^{-a} k) 2^{-a} 2^{-a}
                \to \iint_{\substack{x + y \leq t \\ x, y > 0}} Q(x, y) \d x \d y .
        \]
        Indeed, substituting $x = 2^{-a} m, \, y = 2^{-a} k$, the last sum can be seen as a Riemann sum approximation over the domain $\{x, y \geq 2^{-a}L, \, x + y \leq t\}$. Since $Q$ is smooth up to the boundary of $\{ x, y \geq 0 \}$, the contribution from $\{ 0 \leq x, y \leq  2^{-a}L, \, x + y < t\}$ vanishes as $a \to \infty$.
        Finally, the integral is positive because $Q \geq 0$ and $Q$ cannot vanish on the whole domain of integration since it is a nonzero polynomial.

    \end{pf}

    \smallskip
    \begin{claim}
        There exists $C$ such that $P(k, m) \leq C m^{n-1}$ for $m \geq 1$ and $k \leq L$.
    \end{claim}
    \begin{pf}
        For any fixed $k$, we have that $\bigoplus_m (\gr M)_{m, k}$ is a finitely generated $\bbC[x_1, \dots, x_{N_1}]$-module. Thus for $m \gg 0$, the function $m \mapsto P(k, m)$ is a polynomial with positive leading term. Consequently for large $t$,
        \[
            S_P(t) = \sum_{m \leq t -k} P(k, m)
        \]
        is a polynomial of degree $1 + \deg P(k, -)$ in $t$.
        Since $S_P(t)$ has degree $\leq n$ in $t$ for large $t$, we conclude that $\deg P(k, -) \leq n-1$ and the claim follows.
    \end{pf}

    Note that if $F_k M_{m, \alpha} \neq \{0\}$ then necessarily $\langle \alpha, \xi \rangle \geq k w_2 + m w_1$ where $w_1$ and $w_2$ are the minimum weights of the action of $\xi$ on $\bbC^{N_1+1}$ and $\bbC^{N_2}$ respectively. Recall that $w_2 > 0$ because $\xi$ is a Reeb vector field, while $w_1$ could be negative.
    Using \cref{thm:BiHilbertPoly} we obtain
    \begin{align}
            m^{-n} \sum_{\alpha \in \ft^*} e^{-\langle \frac{\alpha}{m}, \, \xi \rangle} \dim M_{m, \alpha}
                &= m^{-n} \sum_{\alpha \in \ft^*} e^{-\langle \frac{\alpha}{m},  \, \xi \rangle} \sum_{k \geq 0} (\dim F_k M_{m, \alpha} - \dim F_{k-1} M_{m, \alpha}) \\
                &= m^{-n} \sum_{k \geq 0} e^{\frac{-k w_2 - m w_1}{m}} P(m, k)
            \p
    \end{align}
    Then we use that for large $m$ we have $P(m, k) \lesssim m^{n-1}$ when $k \leq L$ and $P(m, k) \lesssim m^{n-2} + k^{n-2}$ otherwise to continue the inequality as
    \begin{align}
            \cdots &\lesssim e^{-w_1} m^{-1} \left[ \sum_{0 \leq k \leq L} e^{- w_2 k / m} + \sum_{k \geq L} e^{-w_2 k / m} \frac{k^{n-2}}{m^{n-2}} m^{-1} \right] \\
                   &\eqsim m^{-1} \left( L + \int_{L/m}^\infty e^{-w_2 x} x^{n-2} \d x + o(1) \right) = O(m^{-1})
            \c
    \end{align}
    where again we have identified the sum as a Riemann integral as $m \to \infty$ and we crucially use that $w_2 > 0$. Here we use $\eqsim$ and $\lesssim$ to mean equality or inequality up to a constant factor. Taking the limit $m \to \infty$ we are done.
\end{proof}

Next, we show that degenerating $X$ along a good $\bbC^*$-action preserves both the algebraic and analytic weighted volume functionals.
\begin{lem} \label{thm:FurtherDegen}
    Let $\rho: \bbC^* \to \bbT^{N_1 + 1}\times \bbT^{N_2} \subseteq \Aut(\amb, \cO(1))$ be a linear $\bbC^*$-action which commutes with $\xi$.
    Let $X$ be a subscheme of $\amb$ and $X_0$ its flat limit with respect to $\rho$ as in \cref{thm:flatlimit}, with corresponding ideal $I_0 \ideal S$.
    Then the algebraic and analytic volumes of $(X,\, \xi)$ agree with those of the central fiber $(X_0, \, \xi)$:
    \[
        \bbW_X(\xi) = \bbW_{X_0}(\xi)
        \qtext{and}
        \int_X e^{-\theta_\xi} \omega^n = \int_{X_0} e^{-\theta_\xi} \omega^n
        \p
    \]
    where $\bbW_{X_0} := \bbW_{(S / I_0)^{(p)}}$.
\end{lem}

\begin{rmk}
    The character $\bbW_{(S/I_0)^{(p)}}$ is essentially just the volume functional $\bbW_{X_0}$ on $X_0$. This is not literally true in most cases, as the scheme $X_0$ will be neither normal nor reduced, so defining $K_{X_0}$ is not possible. If $X_0$ happens to be normal then it is automatically $\bbQ$-Gorenstein and $\bbW_{X_0} = \bbW_{(S / I_0)^{(p)}}$.
    By abuse of notation, we will keep writing $\bbW_{X'}$ for $\bbW_{(S / I')^{(p)}}$ when $X'$ is the scheme corresponding to the ideal $I' \ideal S$.
\end{rmk}

\begin{proof}

    Since $\rho$ commutes with $\xi$, it follows that $\xi$ remains tangent to each $X_\tau \subset \amb$.  Consider the integral
    \[
        \int_{X_\tau} e^{-\theta_\xi} \omega^n = \int_X \rho(\tau)^*(e^{-\theta_\xi} \omega)
        \p
    \]
    Let $\eta$ be the infinitesimal generator of $\rho$. Writing $\tau = u + i v$ we compute
    \[
        \ddat{u}{u + i 0} \int_X \rho(u)^* (e^{-\theta_\xi} \omega) = \int_{X_u} \cL_\eta(e^{-\theta_\xi} \omega^n) = \pm \frac{1}{n+1} \int_{X_u} \eta \iprod (J\xi \iprod \omega^{n+1}) \pm \int_{X_u} \d (e^{-\theta_\xi} \eta \iprod \omega^n)
    \]
    where we have used $\d (e^{-\theta_\xi} \omega^n) = \pm (J \xi \iprod \omega)\wedge \omega^n = \frac{1}{n+1} J\xi \iprod \omega^{n+1}$. Both integrals on the RHS vanish separately: The first one is proportional to the integral of the restriction $(J\xi \iprod \eta \iprod \omega^{n+1})|_{X_\mathrm{reg}}$. But since $J \xi$ is tangent to $X_{u, \mathrm{reg}}$ that restriction is equal to $J \xi \iprod (\eta \iprod \omega^{n+1})|_{X_{u, \mathrm{reg}}}$, which is zero for degree reasons. Meanwhile $e^{-\theta_\xi} \eta \iprod \omega^n$ extends smoothly to $\bbP^{N_1} \times \bbP^{N_2}$ by \cref{thm:extendForms} and then
    \[
        \int_{X_u} \d( e^{-\theta_\xi} \eta \iprod \omega^n ) = \int_{\overline{X}_u} \d ( e^{-\theta_\xi} \eta \iprod \omega^n ) = 0
    \]
    because integration over $\overline{X}_u$ is a closed current.

    By \cref{thm:wcont} the integrals
    \begin{equation} \label{eq:anFutDeg}
        \int_{X_\tau} e^{-\theta_\xi} \omega^n = \langle [X_\tau], e^{-\theta_\xi} \omega^n \rangle
    \end{equation}
    are continuous in $\tau \in \bbC$. So since we just saw that they are constant for $u \in \bbR_{>0}$ we conclude that they are constant for $u \in \bbR_{\geq 0}$. Thus the values at $\tau = 0$ and $\tau = 1$ agree, which is the claimed equality of the analytic volume functionals.

    Since $\rho$ commutes with $T$, we have a multigrading
    \[
        S = \bigoplus_{m, \alpha, w} S_{m, \alpha, w},
    \]
    where $m$ is the degree of homogeneous functions, $\alpha \in \ft^*$ the weight of the $T$-action and $w \in \bbZ$ the weight of $\rho$.
    The ideal $I$ is homogeneous for the first two gradings, so it inherits these gradings, but since $\rho$ does not preserve $I$, the grading by $\rho$ only descends to a filtration, see \cref{sub:test_configurations_and_filtrations}. This filtration is compatible with the gradings in the sense that every $F_i I$ has a decomposition
    \[
        F_i I = I \cap \bigoplus_{w \leq i} \bigoplus_{m, \alpha} S_{m, \alpha, w}
              = \left( \bigoplus_{m, \alpha} I_{m, \alpha} \right) \cap \left( \bigoplus_{m, \alpha}  \bigoplus_{w \leq i} S_{m, \alpha, w} \right)
              = \bigoplus_{m, \alpha} (I_{m, \alpha} \cap F_i I)
        \p
    \]
    Consequently $I_0 \cong \gr I$ inherits a grading by $m$ and $\alpha$. But $\dim (\gr I)_{m, \alpha} = \dim I_{m, \alpha}$, so the RHS of \cref{eq:volModule} is the same for $M = (S/I)^{(p)}$ and $M = (S / I_0)^{(p)}$.
    Thus $\bbW_X = \bbW_{(S/I)^{(p)}} = \bbW_{(S/I_0)^{(p)}}$.
\end{proof}

\begin{proof}[Proof of \cref{thm:volumeAgrees}]

    \textit{Step 1: Degenerate $(X, \, \xi)$ to a union of hyperplanes.}
    Let $I \ideal S$ be the ideal corresponding to $X$.
    We are looking for a $\bbC^*$-action $\rho: \bbC^* \to \bbT^{N_1 + 1} \times \bbT^{N_2} \subseteq \GL(N_1 + 1; \bbC) \times \GL(N_2; \bbC)$ %
    such that $I_0 \coloneq (\ini_\rho f \mid f \in I)$ is a monomial ideal (see \cref{thm:degenerationRings}).

    Naively, we would want to assign weights $a_0, \dots, a_{N_1 + N_2} \in \bbR$ to the variables $x_0, \dots, x_{N_1}, \, y_1, \dots, y_{N_2}$ which are linearly independent over $\bbQ$. Then different monomials would always have different weights, ensuring that the initial term is always a monomial. %
    Of course, unless $a_i \in \bbZ$ this does not define for us a $\bbC^*$-action, but
    \[
        x_0^{k_0} \cdots x_{N_1}^{k_{N_1}} y_1^{k_{N_1 + 1}} \cdots y_{N_2}^{k_{N_1 + N_2}}
            \mapsto a_0 k_0 + \cdots + a_{N_1 + N_2} k_{N_1 + N_2}
    \]
    does define a monomial ordering $>$ on $S$ (cf.\ \cite[\S 15.2]{Eisenbud}).
    And by \cite[Prop.~15.16]{Eisenbud} any monomial ordering $>$ can be approximated with respect to any ideal $I$ by a grading by integers $\lambda$, in the sense that $\ini_{>}(I) = \ini_\lambda(I)$. Since $\ini_> (I)$ is a monomial ideal by the choice of $a_i$, this is sufficient for our purpose.
    Then by \cref{thm:FurtherDegen}, if we can prove the theorem in the case of $X = \cV(I_0)$ we will have proven it for our original $X$.

    \textit{Step 2: $X$ is supported on a union of hyperplanes.}
    We are reduced to the case where $X$ is a, possibly singular, scheme defined by the monomial ideal $I$. Equivalently, $X$ is supported on a union of $n$-planes in $\amb$.
    Let $X_1 \subseteq X$ be an irreducible component and let $X' \subseteq X$ be the union of the remaining irreducible components of $X$. At the level of ideals, $X_1 = \cV(Q_1), \, X' = \cV(Q_2)$ and $I = Q_1 \cap Q_2$. Then
    \[
        \dim(S_{m, \alpha} / I_{m, \alpha}) = \dim(S_{m, \alpha} / (Q_1)_{m, \alpha}) + \dim(S_{m, \alpha} / (Q_2)_{m, \alpha}) - \dim(S_{m, \alpha} / (Q_1 + Q_2)_{m, \alpha})
        \p
    \]
    Using \cref{thm:FutakiLowerDim} we can see that $\bbW_{S / (Q_1 + Q_2)}(\xi) = 0$, so the third term will not contribute to \cref{eq:volumeI}: Indeed
    \[
        \dim \supp (S / (Q_1 + Q_2)) = \dim C(X_1 \cap X') = 1 + \dim X_1 \cap X' \leq n
        \c
    \]
    where $C(X_1 \cap X') \subseteq \bbC^{N_1 + 1} \times \bbC^{N_2}$ is the cone over $X_1 \cap X'$.
    By induction we thus obtain $\bbW_X = \sum_j \bbW_{X_j}$ where $X_j$ are the irreducible components of $X$.
    Since the current of integration on $X$ is just integration on $X_\mathrm{reg} = \bigcup_j (X_j)_\mathrm{reg}$, accounting for multiplicity, we also have $\bbW^\mathrm{an}_X(\xi) = \sum_j \bbW^\mathrm{an}_{X_j}(\xi)$.

    \textit{Step 3: $X$ is supported on a single hyperplane.}
    Thus we are reduced to the case where the scheme $X$ is supported on a single hyperplane in $\amb$. Equivalently, the cone over $X$ is a $(n+1)$-dimensional hyperplane in $\bbC^{N_1+1} \times \bbC^{N_2} \cong \Spec S$. Let $I$ be the ideal defining $X$, and $\fp = \sqrt{I} \ideal S$ the prime ideal defining the supporting hyperplane.
    The linear retraction $\bbC^{N_1 + 1} \times \bbC^{N_2} \surj C(X)_\mathrm{red}$  corresponds to a map $A \coloneq S / \fp \inj S$, splitting the projection $S \to S / \fp$. Using it we can view $S / I$ as an $A$-module. Localizing at $\fp$, we see that $(S / I)_\fp$ is a $A_\fp = \Frac(A)$-vector space. Let $m_i / u_i$ be a basis of $(S / I)_\fp$ as a $\Frac(A)$-vector space. Then $m_i / 1$ is also a basis and if we decompose $S / I \ni m_i = \sum_j m_{ij}$ as a finite sum of weight vectors, then $\{m_{ij} / 1\}_{ij}$ is a generating set for $(S / I)_\fp$. Thus we can pick a basis of $(S/I)_\fp$ consisting of images of weight vectors $v_1, \dots, v_k \in S / I$ under the localization map.
    This gives us commuting maps
    \begin{cd}
        A^{\oplus k} \arrow[r, hook] \arrow[d, "(v_i)_i"]  & \Frac(A)^{\oplus k} \arrow[d, "\sim" arrowdecor, "(v_i / 1)_i"]  \\
        S/I \arrow[r]  & (S/I)_\fp
    \end{cd}
    Furthermore, all maps are $T$-equivariant,
    where $T$ acts on the $j$-th basis vector in $A^{\oplus k}$ and $\Frac(A)^{\oplus k}$ with weight $\wt(v_j)$.
    The diagram shows that $A^{\oplus k} \to S / I$ has to be injective, giving us a $T$-equivariant short exact sequence
    \[
        0 \to A^{\oplus k} \to S / I \to M \to 0
        \p
    \]
    Since localization is exact and the left map in the sequence is an isomorphism after localizing at $\fp$ by construction, we have that $M_\fp = 0$. Thus $\supp(M) \subsetneq \cV(\fp)$ since it does not contain the generic point of $\cV(\fp)$, and hence $\dim\supp(M) \leq n-1$. Since $\bbW$ is additive under short exact sequences, \cref{thm:FutakiLowerDim} shows that
    \begin{align}
        \bbW_{S/I}(\xi) = \bbW_{A^{\oplus k}}(\xi)
            &= - \lim_{m \to \infty} m^{-n} \sum_{\alpha \in \ft^*} \sum_{j = 1}^k e^{-\langle \frac{\alpha}{m}, \, \xi \rangle } \dim A_{m, \alpha - \wt(m_j)} \\
            &= - \sum_{j = 1}^k \lim_{m \to \infty} m^{-n} \sum_{\alpha \in \ft^*} e^{-\frac{1}{m}\langle \alpha + \wt(v_j), \, \xi \rangle} \dim A_{m, \alpha} \\
            &= - k \lim_{m \to \infty} m^{-n} \sum_{\frac{\alpha}{m} \in \ft^*} e^{-\langle \frac{\alpha}{m}, \, \xi \rangle} \dim A_{m, \alpha}
    \end{align}
    and therefore $\bbW_X = \bbW_{(S / I)^{(p)}} = k\bbW_{X_\mathrm{red}}$.
    Furthermore, a similar argument shows that the leading terms of the Hilbert polynomials of $A^{\oplus k}$ and $S / I$ agree. Thus the multiplicity of the scheme $C(X) = \Spec(S / I)$, and therefore also of $X$, is $k$.
    Thus
    \begin{align}
        \frac1{(2\pi)^n} \int_X e^{-\theta_\xi} \frac{\omega^n}{n!}
            = \frac{k}{(2\pi)^n}\int_{X_\mathrm{red}} e^{- \theta_\xi} \frac{\omega^n}{n!}
            = k \bbW_{X_\mathrm{red}}(\xi)
            = k \bbW_{X}(\xi)
    \end{align}
    The middle equality follows from the same Riemann-Roch argument as in \cite[Prop.~5.9]{SunZhang} applied to $\cO(1)|_{X_{\rm red}}$ with curvature $p \omega = \omega_\mathrm{FS} + p \omega_\mathrm{AC}$, since $X_\mathrm{red} \cong \bbP^{n-\ell} \times \bbC^\ell$ is smooth.
    Explicitly, we have
    \[
        \frac{k}{(2\pi)^n}\int_{X_\mathrm{red}} e^{- \theta_\xi} \frac{\omega^n}{n!}
            = \frac{k}{(2\pi)^n}\int_{X_\mathrm{red}} e^{- \theta^{p \omega}_{\xi / p}} \frac{p^{-n} (p\omega_{\rm amb})^n}{n!}
            = p^{-n} \bbW_{S / \fp}(\xi / p) = \bbW_{(S / \fp)^{(p)}}(\xi) \p
    \]

\end{proof}

\section{Proof of the Main Theorem}
\label{sec:maintheorem}

\subsection{K-Polystability of AC Kähler--Ricci Shrinkers}

Throughout this section, we fix the metric $\psi_0 \in \PSH(X; -K_X)$ as our basepoint for defining the energy $\cE_\xi$ and the Ding functional $\cD_{\xi}$. Thus, for any $\phi$, we will have
\[
    \cD_\xi(\phi) = \cD_\xi(\phi, \, \psi_0) = -\frac{1}{(2\pi)^n n! \bbW^\mathrm{an}(\xi)} \cE_\xi(\phi, \, \psi_0) + \cL(\phi)
    \p
\]

\begin{thm}[Semistability] \label{thm:semistab}
    A smooth polarized Fano fibration $(\pi:X \to Y,\, \xi)$ which admits an AC Kähler--Ricci shrinker is K-semistable.
\end{thm}

\begin{proof}
    Combining \cref{thm:volumeAgrees,thm:Slopes} and the definition \eqref{eq:ding} of $\cD_\xi$, we obtain
    \begin{align}
        \Fut_{\xi}(X_0,\, \eta) &= \frac{1}{(2\pi)^n n!}\int_{X_0} \theta_\eta e^{-\theta_\xi} \omega_\mathrm{amb}^n \\
            &= - \frac{1}{2 (2\pi)^n n!} \lim_{t \to \infty} \frac{d}{d t} \cE_\xi(\psi_t, \psi_0) \\[0.3em]
            &=\frac{\bbW^\mathrm{an}(\xi)}{2} \lim_{t \to \infty} \frac{d}{d t} \cD_\xi(\psi_t)
        \p
    \end{align}
    Since the Ding functional $\cD$ is bounded below by \cref{thm:DingBounded}, its slope at infinity is non-negative.
\end{proof}

The rest of this section is devoted to the proof of
\begin{thm}[Polystability]
    In the setting of \cref{thm:semistab}, if
    \[
        \Fut_{\xi}(X_0, \, \eta) = 0
    \]
    then $(\Pi:\cX \to \cY,\,  T, \, \eta)$ is a product test configuration.
\end{thm}

By \cref{thm:isoTrivialProduct} it will be sufficient to construct a biholomorphism $X \to X_0$. For this we will use the geodesic ray constructed in \cref{thm:GeodesicRay}.

\begin{lem}
    Let $(\Pi:\cX \to \cY, \, T , \, \eta)$ be a special test configuration for $(\pi:X \to Y, \, \xi)$ with vanishing Futaki invariant. If $\Phi_t$ is the geodesic ray constructed in \cref{thm:GeodesicRay} with initial data $\Phi_0 = \phi_\mathrm{KRS}$, then $\cL(\Phi_t)$ is constant in $t$.
\end{lem}
\begin{proof}
    Recall that $\Phi$ is constructed as the increasing limit of $\Phi^k \in \PSH(\cX; -K_{\cX})$. The slope of the convex function $\cL(\Phi^k_t)$ is nondecreasing, and since $\cL(\Phi^k_t)$ agrees with $\cL(\psi_t)$ for $t > k$, its slope is zero at infinity by \cref{thm:Slopes}. Hence $\cL(\Phi^k_t)$ is nonincreasing, so
    \begin{equation} \label{eq:Lineq}
        \cL(\Phi^k_t) \leq \cL(\Phi^k_0) = \cL(\Phi_0)
        \p
    \end{equation}
    Next note that
    \begin{equation} \label{eq:Llimit}
        \cL(\psi_t) \leq \cL(\Phi^k_t) = - \log \int_X e^{-\Phi^k_t} \nearrow - \log \int_X e^{-\Phi_t} = \cL(\Phi_t)
    \end{equation}
    by monotone convergence. Thus taking limits of \cref{eq:Lineq} we see that $\cL(\psi_t) \leq \cL(\Phi_t) \leq \cL(\Phi_0)$, and therefore $\cL(\Phi_t)$ also has zero slope at infinity.

    On the other hand, for $t \in [0, k]$, the family $(\Phi^k_t)_t$ is a geodesic, so $\cD(\Phi^k_t)$ is convex by \cref{thm:Dingprops}(2). Thus we observe that for $t \leq k$
    \begin{equation} \label{eq:cmpDingSlope}
        \cD_\xi(\Phi^k_t) \leq \frac{k-t}{k} \cD_\xi(\Phi^k_0) + \frac{t}{k} \cD_\xi(\Phi^k_k)
            = \frac{k-t}{k} \cD_\xi(\phi_\mathrm{KRS}) + \frac{t}{k} \cD_\xi(\psi_k)
        \p
    \end{equation}
    By \cref{thm:Slopes} and \cref{thm:volumeAgrees}
    \[
        \lim_{k \to \infty} \frac{1}{k} \cD(\psi_k) =  \frac{2}{\bbW^\mathrm{an}_X(\xi)} \Fut_\xi(X_0, \, \eta) = 0
    \]
    Thus by \cref{thm:EcontInc} and the convergence in \cref{eq:Llimit} we have, for every $t \in [0, \infty)$
    \[
        \cD_\xi(\Phi_t) = \lim_{k \to \infty} \cD_\xi(\Phi^k_t) \leq \cD_\xi(\phi_\mathrm{KRS})
        \p
    \]
    But $(\Phi_t)_t$ starts at a shrinker, so by \cref{thm:Dingprops}(3), $\cD_\xi(\Phi_t)$ is increasing, hence the bound above shows it is constant.

    By \cref{thm:Dingprops}(2) we know that $\cE_\xi(\Phi_t)$ is affine. Thus $\cL(\Phi_t) = \cD_\xi(\Phi_t) + \frac{1}{(2\pi)^n n! \bbW^{\mathrm{an}}_X(\xi)}\cE_\xi(\Phi_t)$ is also affine, hence a constant because it has zero slope at infinity.
\end{proof}

Since they primarily rely on the ``capacity convergence theorem'' \cite[Thm.~4.26]{GZ}, the proofs of \cite[Thm.~3.22,  Thm.~3.25]{CarlosUniq} apply with trivial modifications to nondecreasing sequences as well. Thus we have
\begin{prop}
    \label{thm:EcontInc}
    Let $\phi_j, \phi'_j$ be nondecreasing sequences in $\PSH(X; -K_X)$ converging to $\phi_\infty, \phi'_\infty$. Then
    \[
        \cE_\xi(\phi_j, \phi'_j) \to \cE_\xi(\phi_\infty, \phi'_\infty)
    \]
    as $j \to \infty$.
\end{prop}

\begin{lem}
    If $\Phi_t \in \Hsing^T, t \in [0, \infty)$ is a geodesic ray along which $\cL(\Phi_t)$ is constant, then $\Phi_t = (F_t)_* \Phi_0$ for a continuous family of biholomorphisms $F: [0, \infty) \times X \to X$.
\end{lem}
\begin{proof}
    The family $F$ was constructed \cite[\S~4.4,~4.5]{CarlosUniq}, and it was shown that $i \partial\bar{\partial} \Phi_t = (F_t)_* i \partial\bar{\partial} \Phi_0$.
    This means $f_t \coloneq \Phi_t - (F_t)_* \Phi_0$ is pluriharmonic, hence it can be locally (actually even globally since $\pi_1(X) = 0$ \cite{Simplyconnected, SunZhang}) written as $f_t = \Re g_t$ for a holomorphic $g_t: X \to \bbC$. Since $\xi (f_t) = 0$ by the $T$-invariance of $\Phi_t$, we see that $\xi(g_t)$ is an imaginary-valued holomorphic function, hence a constant $i a$, $a \in \bbR$. This allows us to conclude that $-J\xi(f_t) = \Re (-J \xi \iprod \d g_t) = a$ and therefore $f_t(\gamma_{-s}(p)) = f_t(p) - s a$. But since $f_t$ is bounded on $\{r \leq 1\}$ and $r(\gamma_{-s}(p)) = e^{-s} r(p)$, $a$ has to be zero, so $f_t$ is constant in the direction of $-J\xi$. Now let $q$ be the maximum of $f_t$ on the compact set $\{r \leq 1\}$. Then $f(\gamma_{-1}(q)) = f(q)$, so $\gamma_{-1}(q)$ is an interior local maximum of $f_t$. Thus $f_t$ is a constant $c_t$ by the maximum principle and $\Phi_t = (F_t)_* \Phi_0 + c_t$.
    Finally,
    \[
        \cL(\Phi_t) = \cL((F_t)_* \Phi_0 + c_t) = c_t + \cL(\Phi_0)
        \c
    \]
    so the assumption that $\cL(\Phi_t)$ is constant shows $c_t = 0$, as desired.
\end{proof}

Recall that when we embed a test configuration $\iota: \cX \hookrightarrow \ambC$ as in \cref{thm:testconfigurationscanbeembedded}, we also obtain a lift $\hat{\lambda}$ of the degenerating action $\lambda$ to $\cO_{\amb}(1)$. When it is clear from the context, we also use $\hat{\lambda}$ to denote the corresponding lift to $\cO(-1)$.

Further recall from \cref{def:psiamb} that
\[
    e^{-p \psiamb} = e^{-h} h_\mathrm{FS}
\]
where $h_\mathrm{FS}$ is the Fubini--Study metric on $\cO(1)$ and $h$ is a smooth function on $\bbC^{N_2}$ which agrees with the cone potential $r^2/2$ of a cone metric on $\bbC^{N_2}$ for $r > 4$.

The family of biholomorphisms $F_t$ can now be used to construct a biholomorphism from $X$ to the central fiber of the test configuration:
\begin{prop}
    The family of holomorphic maps
    \begin{equation}
        G_t \coloneq \lambda_t \circ F_t: X \to X_{e^{-t}} \subseteq \amb
    \end{equation}
    has a subsequential limit $G_\infty$ in the compact-open topology, which is a biholomorphism $X \iso X_0$.
\end{prop}
\begin{proof}
    \textit{Step 1: Existence of $G_\infty$.}
    Note that $F_t$ has a canonical lift to a map $\hat{F}_t: \pm p K_X \to \pm p K_X$. %
    From the embedding $\iota:(X, -pK_X) \to (\amb, \cO(1))$ we obtain an isomorphism $\widehat{\iota}: \mp p K_X \iso \cO(\pm 1)|_X$. We define
    \[
        \hat{G}_t \coloneq \hat{\lambda}_t \circ \hat{\iota} \circ \hat{F}_t: p K_X \to \cO_\amb(-1) \cong \cO_{\bbP^{N_1}}(-1) \times \bbC^{N_2}
        \p
    \]

    Firstly, observe that the projection $\pr_2 \circ \widehat{G}: pK_X \to \bbC^{N_2}$ factors through the projection $\pr: p K_X \to X$, as can be readily seen by tracing out the definitions. Hence we can decompose $\hat{G}_t$ as $(g_t, s_t \circ \pr)$, where
    \begin{alignat}{2}
        g_t&= \pr_1 \circ \hat{G}_t : p K_X \to \cO_{\bbP^{N_1}}(-1) \\[0.1em]
        s_t &= \pr_2 \circ \lambda_t \circ \iota \circ F_t: X \to \bbC^{N_2}
        \p
    \end{alignat}
    Let $\sigma: \cO_{\bbP^{N_1}}(-1) \to \bbC^{N_1+1}$ be the blowdown map.
    \begin{claim} \label{claim:bounded}
        The sequences $s_t$ and $\sigma \circ g_t$ of vector-valued holomorphic functions are uniformly bounded.
    \end{claim}
    \begin{pf}
        Showing that $s_t$ are locally uniformly bounded amounts to showing that $s_t^* r^2$, which is the asymptotic cone metric of  $\hat{G}_t^* \psi$, is locally $t$-uniformly bounded. Since $\Phi \geq \psiamb$, we have $\Phi_t \geq \psi_t = \hat{\lambda}_t^* \psiamb$, and pulling back by $\hat{F}_t$ gives
        \begin{equation} \label{eq:psiPhieasy}
            \hat{G}_t^* \psiamb = \hat{F}_t^* \psi_t \leq \hat{F}_t^* \Phi_t = \Phi_0
            \p
        \end{equation}
        Since $\Phi_0 = \phi_\mathrm{KRS}$ is a smooth AC shrinker by assumption, it is asymptotic to a cone metric on $Y$ with potential $r_\mathrm{KRS}^2$. From \cref{thm:comparePhiR} we conclude $s_t^* r^2 \leq r_\mathrm{KRS}^2$, which is the desired $t$-uniform local bound.

        Next we show that $\sigma \circ g_t: pK_X \to \bbC^{N_1+1}$ are locally uniformly bounded.
        The metric $e^{-p \psiamb} = e^{-ph} h_\mathrm{FS}$ on $\cO_\amb(1)$ induces a metric $e^{+p \psiamb} = e^{p h} h_\mathrm{FS}^\vee$ on $\cO_\amb(-1)$.
        We can think of the metric $p\psi$ as a function $e^{+p\psi}|\cdot|^2: \cO_{\amb}(-1) \to \bbR_+$, which in this case is given by
        \[
            e^{+p \psiamb} |\cdot|^2 = e^{p h} \sigma^* |z|^2
            \c
        \]
        where $z$ is the standard coordinate function on $\bbC^{N_1+1}$.
        Using the equation above we can deduce
        \[
            |\sigma \circ g_t|^2 = \hat{G}_t^* |z|^2 = \hat{G}_t^* (e^{-p h} e^{+p \psiamb} | \cdot |^2) = e^{-p (h \circ s_t)} e^{+p \hat{G}_t^* \psiamb} |\cdot|^2
                \leq e^{+p \Phi_0} |\cdot|^2
            \c
        \]
        where we again used \cref{eq:psiPhieasy}, and also $h \geq 0$. Since the right hand side is independent of $t$, this completes the proof.
    \end{pf}

    By Montel's theorem there exists a sequence $t_i \to \infty$ such that $s_{t_i} \to s_\infty: X \to \bbC^{N_2}$ and $\sigma \circ g_{t_i} \to \bar{g}_\infty: p K_X \to \bbC^{N_1+1}$ in the compact-open topology.
    \begin{claim}
        The limit $\bar{g}_\infty$ maps $p K_X^\circ = p K_X \setminus \{\text{zero section}\}$ to $\bbC^{N_1 + 1} \setminus \{0\}$.
    \end{claim}

    \begin{pf}
        Choose any $v \in p K_X^\circ$. By \cref{claim:bounded},
        \[
            K = K(v) \coloneq \overline{\left\{\big( \hat{G}_t(v), t \big) ~\middle|~ t \in [0, \infty) \right\}} \subseteq \cO(-1)|_{\cX}
        \]
        is compact.
        Since $p \psiamb$ and $p \Phi$ both are locally bounded metrics on $\cO(-1)|_\cX$, there exists a constant $C_1 = C_1(v)$ such that $p \Phi|_K - C_1 \leq p \psiamb|_K$, and therefore
        \[
            e^{+p \psiamb}\big\lvert \hat{G}_t(v) \big\rvert^2 \geq e^{-C} e^{+p \Phi} \big\lvert \hat{G}_t(v) \big\rvert^2 = e^{-C_1} e^{+p\Phi_0}|v|^2
            \c
        \]
        again using that $\hat{G}_t^* \Phi = \hat{F}_t^* \Phi_t = \Phi_0$. Furthermore, by the definition of $\psiamb$,
        \[
            e^{p \psiamb}|\hat{G}_t(v)|^2 = e^{p h(s_t(v))} |\sigma \circ g_t(v)|^2 \leq e^{p C_2} |\sigma \circ g_t(v)|^2
        \]
        where $C_2$ is a $t$-uniform bound on $h(s_t(v))$, which we know exists by \cref{claim:bounded}. Combining these two inequalities we obtain
        \[
            |\sigma \circ g_t(v)|^2 \geq e^{-C_1 - p C_2} e^{-\Phi_0} |v|^2
        \]
        showing that the limit of the LHS is necessarily bounded away from zero.
    \end{pf}
    Since $g_{t_i}|_{p K_X^\circ}$ are equivariant for the standard $\bbC^*$-action on $p K_X^\circ$ and the standard $\bbC^*$-action on $\bbC^{N_1 + 1} \setminus \{0\}$, the same is true of their limit $\bar{g}_\infty|_{p K_X^\circ}$. Taking $\bbC^*$-quotients, we see that the induced maps $\check{g}_{t_i}: X \to \bbP^{N_1}$ converge in the compact-open topology to a map $\check{g}_\infty: X \to \bbP^{N_1}$ induced by $\bar{g}_\infty$.
    Decomposing $g_t$ as
    \[
        g_t = (\check{g}_t, \sigma \circ g_t): p K_X \to \cO_{\bbP^{N_1}}(-1) \subseteq \bbP^{N_1} \times \bbC^{N_1 + 1}
        \c
    \]
    we see that $g_{t_i} \to (\check{g}_\infty, \bar{g}_\infty)$.

    Thus since $g_{t_i}$ and $s_{t_i}$ have compact-open limits, $\hat{G}_{t_i}$ and $G_{t_i}$ have compact-open limits $\hat{G}_\infty$ and $G_\infty$.
    Since $G_t(X) \times \{e^{-t}\} = X_{e^{-t}} \times \{e^{-t}\} \subseteq \cX$ and $\cX$ is closed, we must have $G_\infty(X) \subseteq X_0$.
    Furthermore, $G_t$ and $\hat{G}_t$ are $T$-equivariant for all $t$, so $G_\infty$ and $\hat{G}_t$ are also $T$-equivariant.

    \textit{Step 2: $G_\infty: X \to X_0$ is a branched covering.}
    First consider the map
    \[
        \pr_2 \circ G_\infty : X \to \bbC^{N_2}
        \p
    \]
    Since its target is affine, $\pr_2 \circ G_\infty$ factors through $\pi: X \to Y$, which in our setting simply contracts all compact analytic subsets of $X$. Thus we have a map
    \[
        G_c: Y \to \pr_2(X_0) = Y_0 \subseteq \bbC^{N_2}
        \p
    \]
    We now show that $G_c$ is a biholomorphism.
    \begin{claim} \label{thm:strLowerbound}
        $r_\mathrm{KRS}^2 \leq C s_t^* r^2$ for some $C > 0$.
    \end{claim}
    \begin{pf}
        Let $K = \cX_\bbD \cap \{ 5 \leq r \leq 6\} \subseteq \bbP^{N_1} \times (\bbC^{N_2} \setminus \{0\}) \times \bbD$.
        First we cover $\bbP^{N_1} \times (\bbC^{N_2} \setminus \{0\}) \times \bbD$ with finitely many $J\xi$-equivariant charts $V_\alpha$ for $\amb$ such that $K_{\cX / \bbC}$
        is trivial on $\cX =\cap V_\alpha$.
        Since $\xi$ is nonvanishing on $\bbP^{N_1} \times (\bbC^{N_2} \setminus \{0\})$, on every $V_\alpha$ we can find a nonvanishing \emph{$\xi$-invariant} holomorphic section $u_\alpha$ of $-K_{\cX / \bbC}$ on $\cX \cap V_\alpha$.
        Now to any metric $\phi \in \PSH(\cX; K_{\cX / \bbC})$ we can associate a PSh function $\phi^\alpha$ by 
        \[
            e^{-\phi^\alpha} = |u_\alpha|^2_{\phi}
            \p
        \]
        Since $u_\alpha$ is $J\xi$-invariant, this has the property that $(\gamma_s^* \phi)^\alpha = \phi^\alpha \circ \gamma_s$.

        Now we apply \cref{thm:notHarnack} with $Z = V_\alpha \cap X_0$ to $(\gamma_s^* \Phi)^\alpha$ in $V_\alpha$. We obtain finitely many sets $U_{\alpha, i} \subseteq \cX_\bbD \cap V_{\alpha}$ such that $\{U_{\alpha, i}\}_{\alpha, i}$ cover $K \cap X_0$, and compact sets $K_{\alpha, i} \subseteq (V_{\alpha} \cap \cX_\bbD) \setminus X_0$, such that
        \[
            (\gamma_s^* \Phi^\alpha)|_{U_{\alpha, i}} \leq \sup_{K_{\alpha, i}} (\gamma_s^* \Phi)^\alpha = \sup_{\gamma_s(K_{\alpha, i})} \Phi^\alpha
            \c
        \]
        Without loss of generality we can take all $U_{\alpha, i}$ to be contained in $\{r \geq 4\}$, and to be precompact in $V_{\alpha}$.

        Combining the above inequality with \cref{thm:PhiBound} and the fact that $t = -\log |\tau|$ is bounded on the compact sets $K_{\alpha, i} \subseteq \cX^*$,
        we further obtain
        \begin{equation} \label{eq:scalePhi}
            (\gamma_s^* \Phi)^\alpha|_{U_{\alpha, i}}  \leq \sup_{\gamma_s(K_{\alpha, i})} \left(\psiamb^\alpha + |\tau|^\kappa r^2 + O(\log |\tau|)\right)
                \leq \sup_{\gamma_s(K_{\alpha, i})} \psiamb^\alpha + O(e^{2s})
            \p
        \end{equation}
        Recall that $e^{-\psiamb} = e^{-h} h_\mathrm{FS}$, with $h = r^2/2$ for $\{r \geq 4\}$. Thus for $\{r \geq 4\}$
        \[
            \psiamb^\alpha = r^2/2 + R_\alpha
            \c
        \]
        for $R_\alpha = \phi_\mathrm{FS}^\alpha$. By \cref{eq:FSpullbacklinear}, $\max\{\sup_{\gamma_s(K_{\alpha, i})} R_\alpha,\, \inf_{U_{\alpha, i}} R_\alpha\} \leq C_{\alpha, i} (1 + s)$ for $s \geq 0$. Thus we can continue \cref{eq:scalePhi}, for $s \geq 0$ as

        \begin{align}
            (\gamma_s^* \Phi)^\alpha|_{U_{\alpha, i}}
                &\leq \sup_{\gamma_s(K_{\alpha, i})} \psi^\alpha + C_1 e^{2s} - \inf_{U_{\alpha, i}} (\psi^\alpha \circ \gamma_s) + (\psi^\alpha \circ \gamma_s)|_{U_{\alpha, i}} \\
                &\leq (\gamma_s^* \psi)^\alpha|_{U_{\alpha, i}} + C_2 e^{2s}\\
                &\leq (\gamma_s^* \psi)^\alpha|_{U_{\alpha, i}} + \frac{C_2}{16} (\gamma_s^* r^2)|_{U_{\alpha, i}} 
                = \gamma_s^* \left.\left(\psi + \frac{C_2}{16} r^2 \right)^\alpha \right|_{U_{\alpha, i}} \c
        \end{align}
        where we have used that, since $r \geq 4$ on $U_{\alpha, i}$, we have $\gamma_s^*r^2 = e^{2s} r^2 \geq 16 e^{2s}$ on $U_{\alpha, i}$.

        Thus we have an inequality of metrics on $\cX \cap K$
        \[
            \gamma_s^*\Phi \leq \gamma_s^*(\psi + C_3 r^2)
            \p
        \]
        Since $\gamma_{[0, \infty)}(K) = \cX_\bbD \cap \{r \geq 4\}$ we obtain
        $ \Phi \leq \psiamb + C_3r^2 $
        on $\cX_\bbD \cap \{r \geq 4\}$. Now pulling back by $\hat{G}_t$ we get
        \[
            \phi_\mathrm{KRS} = \Phi_0 \leq \hat{F}_t^* \psi_t + C_3 s_t^* r^2
        \]
        and applying \cref{thm:comparePhiR} gives us $r_\mathrm{KRS}^2 \leq (1 + C_3) s_t^* r^2$, as desired.
    \end{pf}
    Now note that $\pr_2 \circ G_{t_i} \to G_c$ in the compact-open topology, hence
    \begin{equation}\label{eq:r2lowerbound}
        G_c^* r^2 = \lim_{i \to \infty} (\pr_2 \circ G_{t_i})^* r^2 = s_{t_i}^* r^2 \geq C^{-1} r_\mathrm{KRS}^2
    \end{equation}
    Thus we see that $G_c^{-1}(o_{Y_0}) = o_Y$.
    We can use this to conclude that $G_c$ contracts no positive-dimensional analytic subsets. Indeed, let $E(G_c)$ be the exceptional set of $G_c$. It is closed \cite[Thm.~L.7]{GunII}, and $(-J\xi)$-invariant because $G_c$ is $(-J\xi)$-equivariant. Since $-J \xi$ is a Reeb vector field on $Y$, the closure of any $(-J\xi)$-orbit contains $o_Y$, but we have just seen that $o_Y \not\in E(G_c)$, hence $E(G_c) = \emptyset$.

    Since $G_\infty$ agrees with $G_c$ on $X \setminus E \cong Y \setminus o_Y$, we know that any analytic set $Z \subseteq X$ contracted by $G_\infty$ must be contained in $E$. However $\hat{G}_\infty^* \cO(1) = -p K_X$ is a positive line bundle in the sense that $\big(\hat{G}_\infty^* \cO(1)\big)|_E$ is ample, and so no subvariety $Z \subseteq E$ can be contracted by $G_\infty$.

    \textit{Step 3: $G_\infty$ is a biholomorphism.}
    We know that $G_\infty$ is a finite branched holomorphic cover of degree $m$, and we wish to show $m = 1$.
    Since for holomorphic maps, the compact-open topology agrees with $C^1_\mathrm{loc}$-topology, we have $G_{t_i}^* \alpha \to G_\infty^* \alpha$ in $C^0_\mathrm{loc}$ for any $(n, n)$ test form $\alpha$ on $\amb$. Furthermore, if $\supp \alpha \subseteq \{r \leq R\}$ then
    \[
        \supp G_{t_i}^* \alpha \subseteq \{G_{t_i}^* r \leq R\} \subseteq \{r_\mathrm{KRS} \leq C^{-1/2} R\}
        \c
    \]
    where the second inclusion follows from \cref{thm:strLowerbound}.
    Thus we in fact have $G_{t_i}^*\alpha \to G_\infty^*\alpha$ in $C^0$, and therefore $\int_X G_{t_i}^* \alpha \to \int_X G_\infty^* \alpha$. Since $\alpha$ was arbitrary, we conclude that $(G_{t_i})_* [X] \weakto (G_\infty)_* [X]$ in the weak topology of currents. On the other hand, we know from \cref{thm:wcont} that $(G_{t_i})_*[X] = [X_{e^{-t_i}}] \weakto [X_0]$, hence $[X_0] = (G_\infty)_* [X]$. 
    On the other hand, the pushforward $(G_\infty)_* [X]$ can be computed on the unbranched locus of the covering $G_\infty$, so we also have $(G_\infty)_* [X] = m [X_0]$. Hence $m = 1$, as desired.
\end{proof}

\section{Examples}
\label{sec:examples}

In this section we discuss applications to polarized Fano fibrations of a special kind: Let $M$ be a compact klt Fano variety, and suppose there exists a nontrivial root of the canonical bundle of $M$. This means that there is a line bundle $q: L \to M$ and an integer $p > 1$ such that $pL = K_M$. In this case we define $X$ to be the total space of $L$ and $Y = L^\times$ to be the geometric contraction of the zero section of $L$, which exists because $-L$ is ample \cite[Prop.~8.8.2]{EGA2}. Then $\pi:X \to Y$ is the obvious contraction map. We claim that $\pi:X \to Y$ is a Fano fibration. Firstly, $Y$ is normal by \cite[Prop.~8.8.6]{EGA2}. Moreover, the short exact sequence
\begin{equation}\label{eq:tangentlinebundle}
    0 \to q^* L \to \cT_X \to q^* \cT_M \to 0
\end{equation}
shows that $-K_X \cong q^*(-K_M + L) = q^*((1-p)L)$. Since $(1-p)L > 0$ this line bundle is ample on the zero section of $X$, which is the only fiber of $\pi$ which is not a point. Thus $-K_X$ is $\pi$-ample.
Finally, $X$ is klt, since $q^{-1}(U) \cong U \times \bbC$ for any $U \subseteq M$ trivializing $L$, and $M$ is klt by assumption.

We can easily turn this into a \emph{polarized} Fano fibration in the following way.
Let $\zeta$ be the vector field on $X$ generating rotation on the fibers of $L$ with period $2\pi$. It is easy to see that if we set $\xi = \lambda \zeta$ for any $\lambda > 0$, we can consider $(\pi:X \to Y, \, \xi)$ as a polarized Fano fibration.
Moreover we can explicitly work out the full Reeb cone of $X$ in terms of $\zeta$ and the vector fields on $M$. Every holomorphic vector field $\eta$ on $M$ has a canonical lift to a vector field $\hat{\eta}$ on $-K_M$, which has the property that $q_* \hat{\eta} = \eta$. Since $-pL = -K_M$, we also obtain a canonical lift $\hat{\eta}$ to $X$. It is also clear that $[\hat{\eta}, \hat{\psi}] = \hat{[\eta, \psi]}$

While this may seem like a somewhat restrictive setting, in fact the majority of the known AC Kähler--Ricci shrinkers discussed in the introduction are of this form. See \cite[Table 1]{ZhangLi} for a complete breakdown.

\begin{lem}
    Every holomorphic vector field $\xi$ on $X$ commuting with $\zeta$ is of the form $\xi = \hat{\eta} + s \zeta$. Consequently, if $T_M$ is a maximal torus in $\Aut(M)$, then $\ft \coloneq \hat{\ft}_M \oplus \bbR \eta$ is the Lie algebra of a maximal torus in $\Aut(X)$.
\end{lem}
\begin{proof}
    Let $\xi$ be such a vector field. The main point is that the condition $[\xi, \zeta] = 0$ allows us to define a pushforward $q_*\xi$ on $M$. Indeed, we firstly see that the flow of $\xi$ must preserve the zero set of $\zeta$, namely the zero section $M$ of $L$. Hence there is a vector field $\xi_M$ on $M$ such that $(\pi_*)_p\xi = (\xi_M)_p$ for each $p \in M$.
    
    We claim that $\xi - \widehat{\xi}_M$ is a multiple of $\zeta$. It is an elementary exercise to show that any vector field $\eta$ on $X$ commuting with $\zeta$ and restricting to $0$ on the zero section of $L$ must be of the form $\eta = f\zeta$, where $f$ is the pullback of a function on $M$ to $X$. From this point it's straightforward that holomorphicity forces $f = c$ to be constant.
\end{proof}

The exact sequence \eqref{eq:tangentlinebundle} is equivariant for the natural $T_M$-actions on $\cT_X, L$ and $\cT_M$. It is also equivariant for the action of $\zeta$ by weight +1 rotation on $L$, the flow of $\zeta$ on $\cT_X$ and the trivial action on $\cT_M$. Consequently, the action of $H^0(M, \cT_M)$ on $-K_X \cong q^*((1-p)L)$ is induced by the action of $H^0(M, \cT_M)$ on $L$, and the action of $\zeta$ on $q^*((1-p)L)$ has weight $-1$ (according to our convention where weight $-1$ means that $\cL_\zeta V = -iV$).

\subsection{Computing the Futaki invariant}

\begin{lem} \label{thm:sectionsUpstairs}
    Using the isomorphism $q_* \cO^{\mathrm{alg}}_X = \bigoplus_{d \geq 0} L^{-d}$ we have
    \[
        \bigoplus_{d \geq 0} H^0_\mathrm{alg}(M, -(k + d)L) \cong H^0_\mathrm{alg}(X, -k\, q^*L),
    \]
     where the action of $\zeta$ on the right has weight $d$ and the action of $\ft_M$ is induced by the action on $L$. In particular,
    \begin{align} \label{eq:iso}
        H^0_\mathrm{alg}(X, -mK_X) &\cong \bigoplus_{d \geq 0} H^0_\mathrm{alg}(M, -(d + (p-1)m)L)
    \end{align}
    $T_M$-equivariantly and with $\zeta$ acting with weight $d-m$.
\end{lem}
\begin{proof}
    This follows from the projection formula since $q_* q^* L^{-k} = L^{-k} \otimes q_*\cO_X$, which locally is just a consequence of $\cO^{\mathrm{alg}}_{U \times \bbC} \cong \cO^{\mathrm{alg}}_U \otimes \bbC[x]$ for affine varieties $U$.
\end{proof}

\begin{prop} \label{thm:VolumeUpstairs}
    If $\dim M = n$ then
    \begin{equation}\label{eq:linebundleweightedvolume}
        \bbW_X(s \zeta + \hat{\eta}) = e^s \int_0^\infty e^{-st} (p-1 + t)^n \, \bbW_M \big((p-1 + t) \eta\big) \d t
    \end{equation}
    for any Reeb vector $s \zeta + \hat{\eta}$. In particular, this integral converges if $s \zeta + \hat\eta$ is a Reeb vector.

    Moreover, the derivatives of $\bbW_X$ can be computed by differentiating the integrand in \cref{eq:linebundleweightedvolume}. In particular, 
    \begin{align} \label{eq:linebundlefutaki}
        \mathrm{Fut}_{s_0 \zeta}(X, \, \hat{\eta}) &= e^{s_0} \int_0^\infty e^{-s_0t} (p-1+t)^n \, \mathrm{Fut}_{0}\big(M, \, ((p-1) + t)\eta\big) \d t  \\
                                      &= \mathrm{Fut}_{0}(M, \, \eta) \left(e^{s_0} \int_0^\infty e^{-s_0t} (p-1+t)^{n+1} \d t  \right) \p
    \end{align}
\end{prop}

\begin{proof}
    The claim is \emph{pointwise in $s$ and $\eta$}, so all of our estimates are allowed to depend on $s$ and $\eta$.

    By \cref{eq:iso} with $m = 0$, the vector $s \zeta + \hat{\eta} \in \bbR \zeta \oplus \ft$ being a Reeb field means
    \begin{equation} \label{eq:ReebCrit}
        k s + \langle \alpha, \eta \rangle > \epsilon ( k + |\alpha| ) \quad \forall k > 0,  \alpha \in \wt(H^0(M, -kL)),
    \end{equation}
    in particular $s$ needs to be positive.
    Plugging \cref{thm:sectionsUpstairs} into the definition of the weighted volume in \cref{eq:volume}, we see that the weighted volume functional on $X$ is
    \begin{align}
        &\bbW_X(s \zeta + \hat{\eta}) \\
        &\quad=
            \lim_{m \to \infty}
            \sum_{d \geq 0}
            m^{-1}
            \left(p-1 + \frac{d}{m}\right)^n
            \sum_{\alpha \in \ft^*}
                \exp\left(- \frac{\langle \alpha, \eta \rangle + s (d-m)}{m}\right)
                \frac{\dim H^0\left(M, -(p-1)m L - d L\right)_\alpha}{((p-1)m + d)^n}
            \\
        &\quad=
            e^s \lim_{m \to \infty}
            \sum_{d \geq 0}
            m^{-1}
            \left(p-1 + \frac{d}{m}\right)^n e^{-s \frac{d}{m}}
            \int_{\ft^*}
                e^{- \left((p-1) + \frac{d}{m}\right) \eta }
                \DH_{(p-1)m + d}
        \label{eq:RiemannSum}
        \p
    \end{align}
    Here for $\psi \in \ft$ we view $e^\psi$ as a function on $\ft^*$ by abuse of notation, and the \textit{spectral measures} $\DH_k$ of $-L$ are given by
    \[
        \DH_k = k^{-n} \sum_{\alpha \in \ft^*} \dim \left( H^0(M, -kL)_\alpha \right) \delta_{\alpha / k}\p
    \]
    By \cite[Prop.~4.1]{BWN} the spectral measures $\DH_k$, which are all supported on a compact polytope $\Delta \subseteq \ft$, converge weakly to the Duistermaat-Heckman measure $\DH$ of $-L$.
    In fact, if $\omega \in 2\pi c_1(-L)$ is a smooth metric and $\mu_\omega: M \to \ft_M^*$ is a suitably normalized moment map for $T_M \acts M$ then $\Delta = \mu_\omega(M)$ and $\DH = C \cdot (\mu_\omega)_* \omega^n$ for some normalization $C > 0$.

    Now we need an $m$-uniform tail estimate for the summands in \cref{eq:RiemannSum}. More specifically, we would like to show that the sum 
    \begin{equation} \label{eq:tail}
        \sum_{d \geq C m} 
            m^{-1}
            \left(p-1 + \frac{d}{m}\right)^n e^{-s \frac{d}{m}}
            \int_{\ft^*}
                e^{- \left((p-1) + \frac{d}{m}\right) \eta }
                \DH_{(p-1)m + d}
    \end{equation}
    is arbitrarily small for $C$ sufficiently large, independent of $m$.
    Here we use again that since $s \zeta + \hat{\eta}$ is Reeb, for $\alpha \in \Delta$
    \begin{equation}\label{eq:reebmeasureintegralcomparison}
        s \frac{d}{m} + \left\langle \alpha, \left( (p-1) + \frac{d}{m} \right) \eta \right\rangle > \epsilon\left((p-1) + \frac{d}{m} + |\alpha|\right) - s(p-1)
    \end{equation}
    so the exponential factors in \cref{eq:tail} are bounded by $\exp(C' - \epsilon d / m)$ and it is straightforward to see that the exponential decay bounds the tails.

    To handle the rest of the sum, we claim that for any given compact $I \subseteq \bbR$ and $\epsilon > 0$, we have
    \begin{equation} \label{eq:DHlimit}
        \left\lvert
            \int_{\ft^*} e^{-x \eta} \DH_k
            -
            \int_{\ft^*} e^{- x \eta} \DH
        \right\rvert < \epsilon
    \end{equation}
    for all $x \in I$ and for $k$ sufficiently large. To see this, notice that since $\left\lVert \DH_k \right\rVert$ is uniformly bounded, it is easy to see that the sequence $g_k(x) = \int_{\ft^*} e^{-x \eta} \DH_k$ is equicontinuous.
    Then since the domain is compact, it follows that the convergence of $g_k$ to the pointwise limit is uniform.

    Thus for that for any fixed $\epsilon > 0$ we have that for all $C$ sufficiently large,
    \begin{align}
        \bbW_X(s \zeta + \hat{\eta})
        &=
            e^s \lim_{m \to \infty}
            \left[
                \epsilon + \sum_{0 \leq d \leq C m} m^{-1} \left(p-1 + \frac{d}{m}\right)^n e^{-s \frac{d}{m}}
                \left(
                    \epsilon + \int_{\ft^*} e^{-\left( (p-1) + \frac{d}{m} \right)\eta} \DH
                \right)
            \right]
        \\
        &=
            e^s \lim_{m \to \infty}
            \left[
                \epsilon + \sum_{0 \leq d \leq C m} m^{-1} \left(p-1 + \frac{d}{m}\right)^n e^{-s \frac{d}{m}}
                \left(\epsilon + \bbW_M\left(\left( (p-1) + \tfrac{d}{m}\right) \eta\right)\right)
            \right]
        \\
        &=
            e^s \int_0^C \left( p-1 + t \right)^n e^{-st} \bbW_M\big(( (p-1) + t ) \eta \big) \, dt + O(\epsilon) \p
    \end{align}
    In the first line we have used the tail estimate and Claim 1. In the third line we have identified the sum as a (proper) Riemann integral. The result then follows by passing to the limit $\epsilon \to 0$, which also forces $C \to \infty$.

    The computation for the derivative follows similarly. In fact, it suffices to show that \cref{eq:tail} is small in the $C^k$ norm, since we have shown that the $d \leq Cm$ part of the sum converges to a Riemann integral on a compact domain, where we can exchange integral and differentiation.
    Since the derivatives of the summands in \cref{eq:tail} are polynomials times a decaying exponential, it is easy to bound their difference quotients. Since the derivative is again a polynomial times the same exponential, the same argument works for higher derivatives.
\end{proof}

\begin{prop}\label{thm:example:slightlymoregeneral}
    Let $M$ be a Fano manifold such that 
    \begin{enumerate}
        \item $M$ is K-unstable with vanishing Futaki invariant $\mathrm{Fut}_0 \equiv 0$,
        \item The Fano index of $M$ is strictly greater than 1, so there exists a non-trivial root $L \to M$ of $-K_M$
        \item there exists a special test configuration $(\cM, \, -K_{\cM/\bbC}, \, \eta_0)$ for $(M, \, -K_M)$ with negative Futaki invariant, such that $L$ extends to a global line bundle $\cL$ on $\cM$.
    \end{enumerate}
    If $X$ denotes the total space of $L$ and $Y$ the contraction of the zero section, then $\pi:X \to Y$ is K-unstable as a polarized Fano fibration with respect to any choice of Reeb vector field $\xi$, and hence does not admit an asymptotically conical Kähler--Ricci shrinker metric. 
\end{prop}

\begin{proof}
    Firstly, it follows from \cref{thm:VolumeUpstairs} that $\bbW_X$ is minimized at some multiple $\xi = c \zeta$, since $\bbW_M$ is minimized at $0$. Indeed, since $\bbW_M(\eta)$ is minimized at $\eta = 0$ we have by \cref{eq:linebundleweightedvolume} that for $\eta \neq 0$
    \begin{align}
        \bbW_X(s \zeta + \hat{\eta}) &= e^s \int_0^\infty e^{-st} (p-1 + t)^n \, \bbW_M \big((p-1 + t) \eta\big) \d t \\
                              &> e^s \int_0^\infty e^{-st} (p-1 + t)^n \, \bbW_M (0) \d t = \bbW_X(s \zeta)
    \end{align}
    because $e^{-st} (p-1 +t) > 0$.
    So any minimum of $\bbW_X$ must lie on the ray $\bbR_+ \zeta$. Furthermore it is easy to see that $\bbW_X(s \zeta) \to +\infty$ as $s \to 0$ and $s \to \infty$, so $\bbW_X$ has a minimum on $\bbR_+ \zeta$. Hence, the only possibility for a stable polarization lies along this ray.

    Next we use this data to build a special test configuration for $X$ whose central fiber $X_0$ is precisely the total space of $\cL|_{M_0}$. First, define $\cX$ to be the total space of $q:\cL \to \cM$, and let $\Pi_{\cX}$ be the composition of $q$ with $\Pi_{\cM}: \cM \to \bbC$. By \cite[Theorem 2]{Fujiki_blowing} applied to $(\Pi_{\cX})|_{\cM}:\cM \to \bbC$, there exists a map $\Pi:\cX \to \cY$ contracting the zero sections of the fibers of $\Pi_{\cX}$. We will consider $X$ as polarized by an arbitrary multiple $\xi = c\zeta$, and since $\cL$ extends as a line bundle we see that any such $\xi$ is well-defined on $\cX$. It is clear that the degenerating action on $\cX$ is generated by the lift $\widehat{\eta}_0$. It is straightforward to check that $(\Pi:\cX \to \cY, \, T, \, \widehat{\eta}_0)$ satisfies all four criteria in \cref{def:testconfiguration}, as we have already seen above that $(\pi_0:X_0 \to Y_0, \, \xi)$ is a polarized Fano fibration.

    Thus there is a corresponding Futaki invariant which can be computed using \cref{eq:linebundlefutaki}. We read immediately that $\Fut_{0}(M_0, \,\eta_0) < 0$ implies $\Fut_{ c \zeta}(X_0, \, \hat{\eta}_0) < 0$ for any $c > 0$ due to the positivity of the integrand. 
\end{proof}

\subsection{Example: $M = \Bl_{Q^2} Q^6$}
In this section we describe a specific example in which \cref{thm:example:slightlymoregeneral} can be applied. Our example $M$ will be six-dimensional, so that $X$ will be a 7-fold. This is the lowest dimension that we could find a relatively simple example, but it is an interesting question whether one can find one in dimensions four or five. These would be the lowest possible dimensions, since all Fano surfaces with vanishing Futaki invariant are K-stable and all K-unstable Fano threefolds have index 1 \cite{fanography}. In \cite[Theorem 1.2]{Del_Examples}, Delcroix shows that the blowup $M_{n,k} =\Bl_{Q^k} Q^n$ of a quadric $Q^n \subset \bbP^{n+1}$ along a $k$-dimensional linear subquadric is K-unstable for $2 \leq n \leq n-3$. Moreover, he shows that the Futaki invariant (at $\xi = 0$) of these manifolds vanishes, hence any Kähler--Ricci shrinker on these manifolds must be Einstein. Thus these satisfy property $(1)$ in \cref{thm:example:slightlymoregeneral}.

The Fano index of such manifolds is also straightforward to compute. Let $\sigma: M_{n,k} \to Q^n$ denote the blowup with exceptional divisor $E = \sigma^{-1}(Q_k)$. The canonical bundle is 
\begin{equation}
\begin{split}
    K_{M_{n,k}} &= \sigma^*K_{Q_n} \otimes \cO((n-k-1)E) \\
                &=\cO_{\bbP^{n+1}}(-n) \otimes \cO((n-k-1)E) ,
\end{split}
\end{equation}
hence $M_{n,k}$ satisfies $(2)$ if and only if $\mathrm{gcd}(n, n-k-1) > 1$. The main result that we show here is that, in the particular case $n = 6$ and $k = 2$, the Fano manifold $M_{6,2}$ satisfies $(3)$. In this case the Fano index is equal to $3$, and for simplicity we will focus on the specific root $L = \cO_{\bbP^{7}}(-2) + \cO(E)$ of $K_{M_{6,2}}$. Applying \cref{thm:example:slightlymoregeneral} gives:
\begin{thm}\label{thm:example}
  Let $M = M_{6,2} = \Bl_{Q^2} Q^6$, $L = \cO_{\bbP^{7}}(-2) + \cO(E)$, so that $3L = K_M$, and $X$ be the total space of $L$. Then $X$ admits no asymptotically conical Kähler--Ricci shrinker metrics.
\end{thm}

We restrict to the particular case above for computational simplicity; every $M_{n, k}$ with $\mathrm{gcd}(n, n-k-1) > 1$ and any choice of root should work the same way. In the general case, it should be possible to verify using the theory of optimal degenerations of spherical Fano varieties, which has been developed extensively by Delcroix (see \cite{Del_Examples} and the references therein). We proceed here by a more ad hoc method, working with an explicit test configuration for $M$ where we can compute by hand.

\newcommand{\blow}{\Bl_{P^3} \bbP^7}

\begin{proof}

Let $P^3 = \bbP(\bbC^4 \times \{0\}) \subseteq \bbP^7$ and note that $M = \Bl_{Q^6 \cap P^3} Q^6$.
Now consider the family $\cQ$ for $Q = Q^6$ inside $\bbP^7 \times \bbC$ given by degenerating by the action 
\[
    \lambda(t) (u, v) = (t u, v) \quad\text{for } u, v \in \bbC^4 \p
\]
Writing $q(x, y, z, w) = x^2 + y^2 + z^2 + w^2$ we can also describe this by
\[
    \cQ = \{([u, v], t) \mid q(u) + t^2 q(v) = 0 \} \p
\]
We see that the central fiber is given by $Q_0 = \{[u, v] \mid q(u) = 0\}$.
The 1-parameter subgroup $\lambda: \bbC^* \to \GL(4; \bbC) \times \GL(4; \bbC)$ preserves $P^3$, and hence induces an action on $\blow$. Therefore we have a degenerating family $\cM$ for $M$, which can be neatly compactified as follows.

Consider the vector bundle $\cV = \cO(1)^{\oplus 4} \oplus \cO^{\oplus 4}$ on $\bbP^1$, and let $\alpha \in H^0(\bbP^1, \cO(1))$ be a section vanishing at $0 \in \bbP^1$. Then
\[
    \bar\cQ = \left\{ ([u, v], \tau) \in \bbP(\cV) ~\middle|~ q(u) + \alpha(\tau)^2 q(v) = 0 \right\}
\]
is an equivariant compactification of $\cQ$. If we take the subbundle $\cS = \cO(1)^{\oplus 4}$  of $\cV$, then our compactification of $\cM$ is precisely the strict transform of
$\bar\cM$ in $\Bl_{\bbP(\cS)} \bbP(\cV)$. We denote the blowdown map by
$\sigma: \bar\cM \to \bbP(\cV)$.

\bigskip
\begin{claim}
    $\overline{\cM}$ is a special test configuration for $M$.
\end{claim}
\begin{pf}
The only thing we need to show is that $M_0$ is normal with klt singularities. Since the subquadric $Q^6 \cap P^3$ does not meet the singular locus of $Q_0$, we only need to show that $Q_0$ itself has these properties. This is standard, but we include the details for convenience of the reader. 

To see that $Q_0$ is normal, first note that any hypersurface is Cohen-Macaulay \cite[Tag 00S8]{stacks}, and thus all local rings satisfy Serre's $S_k$ for any $k$. Since $Q_0$ clearly has singularities only in codimension 3, it is normal. 

To see that $Q_0$ is klt, we check locally around any given singular point $x^* \in Q_0 \cap \{u_0 = u_1 = u_2 = u_3 = 0\}$, noting that the singular set of $Q_0$ is precisely the intersection with this 3-plane. Choosing an appropriate affine coordinate chart, we can locally view $x^*$ as the origin $0 \in \bbC^7$ and $Q_0 = \widetilde{Q} \times \bbC^3$, where
    \[  \widetilde{Q} = \{u_0^2 + u_1^2 + u_2^2 + u_3^2 = 0\} \subset \bbC^4.\]
    We will show that $\widetilde{Q}$, and therefore $Q_0$, actually has terminal singularities. Indeed, consider the blowup of $\bbC^4$ at the origin, and set $Z = \Bl_{0} \bbC^4$ with exceptional divisor $E$. Then $K_Z = \pi^*K_{\bbC^4} + 3 E$. Let $Q' = \Bl_0\widetilde{Q}$, so that $Q'$ is also the strict transform of $\widetilde{Q}$ in $Z$. Since the singularity of $\widetilde{Q}$ at $0$ has multiplicity 2, we have that $\pi^*[\widetilde{Q}] = [Q']  + 2[E]$ as divisors on $Z$. By adjunction: 
    \[ \begin{split} K_{Q'} &= (K_{Z} + \mathcal{O}_Z(Q'))|_{Q'} \\
    &= (\pi^*K_{\bbC^4} + 3E)|_{Q'} + (\pi^*\mathcal{O}_{\bbC^4}(\widetilde{Q}) - 2E)|_{Q'}\\
    &= \pi^*(K_{\bbC^4} + \mathcal{O}_{\bbC^4}(\widetilde{Q})) + E\\
    &= \pi^*K_{\widetilde{Q}} + E.  \end{split}\]

\end{pf}

\begin{claim}
    The Futaki invariant of $\overline{\cM}$ is negative.
\end{claim}

\begin{pf}
We compute this using the Odaka--Wang intersection formula \cite{OdakaSlope,WangHeightGIT}, which in the anticanonically polarized case is
\[
    \Fut_{0}(M_0, \lambda) = -\frac{1}{(n+1) V(M)} (-K_{\bar\cM / \bbP^1})^{n+1}
    \p
\]
Thus $\cM$ is destabilizing iff $(-K_{\cM / \bbP^1})^{n+1} > 0$.

This situation is amenable to the Schubert2 package in Macaulay2, which we use to compute $\Fut_{0}(M_0, \lambda)$ in this case. Note that the equation defining $\bar\cQ$ is a global section $p^*\cO_{\bbP^1}(2) \otimes \cO_{\bbP(\cV)}(2)$ where $p: \cP(\cV) \to \bbP^1$ is the structure map. Moreover since $\bar\cM$ is the strict transform of $\bar \cQ$, it is a section of $\sigma^* \left(p^*\cO_{\bbP^1}(2) \otimes \cO_{\bbP(\cV)}(2)\right)$. The following is the source code which computes that the self-intersection $(-K_{\overline{\cM}/\bbP^1})^{7}$ is given by 
\[(-K_{\overline{\cM}/\bbP^1})^{7} = 2^3 \cdot3^7 = 17496. \]
In particular this is positive, showing that $\Fut_{0}(M_0, \lambda) < 0$. Specifically we have
\[
    \Fut_{0}(M_0, \lambda) = -\frac{2}{49} \p %
\]

\lstdefinelanguage{Macaulay2}{
    sensitive=true,
    morecomment=[l]{--},
    morestring=[b]",
    morekeywords={assert,print}
}
\begin{lstlisting}[
    language=Macaulay2,
    basicstyle=\ttfamily\footnotesize,
    breaklines=true,
    columns=fullflexible,
    keepspaces=true,
    literate={^^cf^^83}{{$\sigma$}}1,
    numbers=left,
    numberstyle=\tiny\color{gray},
    frame=single,
    lastline=34,
    commentstyle=\color{gray}\itshape,
    stringstyle=\color{green!50!black},
]
needsPackage "Schubert2";

P1 = abstractProjectiveSpace 1;
O1 = OO_P1(1);
O0 = OO_P1;
V = O1 ++ O1 ++ O1 ++ O1 ++ O0 ++ O0 ++ O0 ++ O0;

-- Schubert2's projectiveBundle parametrizes rank-one subbundles.
PV = projectiveBundle V;
p = PV.StructureMap
-- "bundles PV" gives a pair (tautological bundle, quotient bundle)
OPV1 = dual first bundles PV;

-- The blowup center is the projective subbundle
--     Sbar = P(O(1)^4) -> P(V).
W = O1 ++ O1 ++ O1 ++ O1;
Sbar = projectiveBundle W;
OSbar1 = dual first bundles Sbar;
-- this describes a map Sbar \to PV over P1 such that O_{PV}(1) is pulled back to OSbar1
centerInclusion = map(PV, Sbar, OSbar1);

-- define the blowup
(BlPV, excBundle, excInclusion, σ) = blowup centerInclusion;

-- Mbar is the strict transform of a smooth quadric in the class O_PV(2) tensor p^*O(2)
Mbar = sectionZeroLocus σ^*( (OPV1)^**2 * (p^* O1)^**2 );
MbarInclusion = Mbar.StructureMap;

K_Mbar = -chern(1, tangentBundle Mbar);
-- K_{Mbar / P1} = K_Mbar - K_P1 = K_Mbar + σ^* T_P1
relK_Mbar = K_Mbar + MbarInclusion^* σ^* p^* chern(1, tangentBundle P1);

MbarIntersection = integral((-relK_Mbar)^7);
print MbarIntersection;
\end{lstlisting}

\end{pf}

The intersection number $(-K_{\overline{\cM}/\bbP^1})^{7}$ can also be computed directly by hand in this situation. If we define $H = c_1(\cO_{\bbP(\cV)}(1))$, $F = p^*c_1(\cO_{\bbP^1}(1))$, and $\cE$ to be the exceptional class of the blowup $\overline{\cM} \to \cQ$, it is not hard to work out that 
\[K_{\overline{\cM}/\bbP^1} = 6H + 2F - 3\cE. \] 
Further, one can derive relations between $H, F$ and $\cE$ which allows one to do the intersection computation directly.
\end{proof}

\linespread{1.00}
\bibliographystyle{amsalpha}
\bibliography{bib.bib}

\end{document}

%% file: preamble_v3.tex
\usepackage[utf8]{inputenc}

\usepackage{amsmath}
\usepackage{amssymb}
\usepackage{amsfonts}
\usepackage{mathtools} % for example for mathrlap

\usepackage{mathrsfs} % for \mathscr
\usepackage{mleftright} % remove extra spacing before \left and after \right
\mleftright % define \left to \mleft and \right to \mright (tighter spacing)

\usepackage[x11names]{xcolor}

\usepackage{tikz-cd}

\usepackage[breaklinks,colorlinks=true,hypertexnames=false,linkcolor=blue,citecolor=green!60!black]{hyperref}
\usepackage{hypcap}

\usepackage{zref-clever}
\zcsetup{cap=true}

\let\cref\zcref
\newcommand{\Cref}[1]{\zcref[S]{#1}}

\makeatletter
\newcommand*{\namecref}[1]{%
  \zref@extract{#1}{zc@type}%
}
\makeatother

\zcRefTypeSetup{equation}{ % reference equations with just numbers
  Name-sg =, name-sg =, Name-pl =, name-pl =, Name-sg-ab =, name-sg-ab =, Name-pl-ab =, name-pl-ab =
}

\AddToHook{env/thm/begin}{%
    \zcsetup{countertype={thm=theorem}}}
\AddToHook{env/rmk/begin}{%
    \zcsetup{countertype={thm=remark}}}
\AddToHook{env/lem/begin}{%
    \zcsetup{countertype={thm=lemma}}}
\AddToHook{env/cor/begin}{%
    \zcsetup{countertype={thm=corollary}}}
\AddToHook{env/prop/begin}{%
    \zcsetup{countertype={thm=proposition}}}
\AddToHook{env/conj/begin}{%
    \zcsetup{countertype={thm=conjecture}}}
\AddToHook{env/defn/begin}{%
    \zcsetup{countertype={thm=definition}}}
\AddToHook{env/ex/begin}{%
    \zcsetup{countertype={thm=example}}}

\zcRefTypeSetup{cd}{
  Name-sg = Diagram , name-sg = diagram , Name-pl = Diagrams , name-pl = diagrams ,
}
\AddToHook{env/cd/begin}{%
    \zcsetup{countertype={equation=diagram}}}

\zcRefTypeSetup{defnprop}{
  Name-sg = Definition/Proposition , name-sg = definition/proposition , Name-pl = Definitions/Propositions , name-pl = definitions/propositions ,
}
\AddToHook{env/defn-prop/begin}{%
    \zcsetup{countertype={thm=defnprop}}}

\zcRefTypeSetup{conj}{
    Name-sg = Conjecture , name-sg = conjecture , Name-pl = Conjectures , name-pl = conjectures ,
}

\zcRefTypeSetup{claim}{
  Name-sg = Claim , name-sg = claim , Name-pl = Claims , name-pl = claims ,
}

\zcRefTypeSetup{problem}{
  Name-sg = Problem , name-sg = problem , Name-pl = Problems , name-pl = problems ,
}

\usepackage{autonum}

\usepackage[shortlabels]{enumitem}

\usepackage{csquotes}
\MakeOuterQuote{"}

\makeatletter
\def\@cite#1#2{[{#1\if@tempswa ,\penalty400~#2\fi}]}% NEW
\makeatother

\usetikzlibrary{arrows}
\usetikzlibrary{babel}
\tikzcdset{every label/.append style = {font = \small}}
\tikzset{
    arrowdecor/.style={anchor=south, rotate=90, inner sep=.5mm}
}
\usepackage{rotating}

\makeatletter
\newcommand{\thesislanguage}[2]{%
  \@ifpackageloaded{babel}{\iflanguage{ngerman}{#1}{#2}}{#2}%
}
\makeatother

\newtheorem{thm}{\thesislanguage{Satz}{Theorem}}[section]
\newtheorem*{thm*}{\thesislanguage{Satz}{Theorem}}
\newtheorem{conj}{\thesislanguage{Vermutung}{Conjecture}}[section]
\newtheorem*{conj*}{\thesislanguage{Vermutung}{Conjecture}}

\newtheorem{lem}[thm]{Lemma}
\newtheorem*{lem*}{Lemma}

\newtheorem{prop}[thm]{Proposition}
\newtheorem*{prop*}{Proposition}

\newtheorem{cor}[thm]{\thesislanguage{Korollar}{Corollary}}

\theoremstyle{definition}
\newtheorem{defn}[thm]{Definition}
\newtheorem*{defn*}{Definition}
\newtheorem*{notation}{Notation}
\newtheorem{defn-prop}[thm]{Definition/Proposition}
\newtheorem{convention}[thm]{Convention}
\newtheorem*{convention*}{Convention}

\theoremstyle{remark}
\newtheorem{rmk}[thm]{\thesislanguage{Anmerkung}{Remark}}
\newtheorem*{rmk*}{\thesislanguage{Anmerkung}{Remark}}

\newtheorem{ex}[thm]{\thesislanguage{Beispiel}{Example}}

\newtheoremstyle{claim}
  {0.5\topsep} % Space above
  {0\topsep} % Space below
  {} % Body font
  {} % Indent amount
  {\itshape} % Theorem head font
  {:} % Punctuation after theorem head
  {.5em} % Space after theorem head
  {} % Theorem head spec (can be left empty, meaning `normal')
\theoremstyle{claim}

\newtheorem{claim}{\thesislanguage{Behauptung}{Claim}}
\newtheorem*{claim*}{\thesislanguage{Behauptung}{Claim}}
\counterwithin*{claim}{thm}

\makeatletter
\newenvironment{pf}{\par
  \pushQED{\qed}%
  \list{}{\settowidth{\leftmargin}{\itshape{Pf:~}}
          \settowidth{\labelwidth}{\itshape{Pf:~}}
          \parsep=0pt \listparindent=\parindent
          \topsep=0pt
  }
\item[\hskip\labelsep\itshape{Proof:}]
}{%
    \popQED\endlist\@endpefalse\smallskip
    \aftergroup\@afterindentfalse\aftergroup\@afterheading
}
\makeatother

\newtheoremstyle{Normal}{}{}{}{}{\bfseries}{:}{.5em}{}
\theoremstyle{Normal}

\counterwithin*{problem}{section}
\makeatletter
\let\oldendproof\endproof
\def\endproof{\oldendproof\aftergroup\@afterindentfalse\aftergroup\@afterheading}
\makeatother

\newcounter{autoItemCounter}

\newlist{xxthmlist}{enumerate}{1}
\setlist[xxthmlist]{
    label=(\roman{xxthmlisti})\protect\label{item:\theautoItemCounter\thexxthmlisti},
    ref=(\roman{xxthmlisti})
}

\usepackage{aliascnt}
\newaliascnt{cd}{equation}

\makeatletter
\newenvironment{cd}{%
    \refstepcounter{cd}%
    $$
        \begin{tikzcd}[sep=large]
}{%
    \end{tikzcd}
    \eqno \hbox{\@eqnnum}
    $$\@ignoretrue
}

\newenvironment{cd*}{%
    $$
    \begin{tikzcd}[sep=large]
}{%
    \end{tikzcd}
    $$\@ignoretrue
}
\makeatother
\makeatletter
\g@addto@macro\bfseries{\boldmath}
\makeatother
\newcommand{\noop}[1]{}

\newcommand{\surj}{\twoheadrightarrow}
\newcommand{\inj}{\hookrightarrow}
\newcommand{\weakto}{\rightharpoonup}
\newcommand{\iso}{\xrightarrow{\sim}}

\AtBeginDocument{
    \renewcommand{\phi}{\varphi}
    \renewcommand{\epsilon}{\varepsilon}
}

\renewcommand{\setminus}{\smallsetminus}

\DeclareMathSymbol{\colonrel}{\mathrel}{operators}{"3A}

\makeatletter
\newcommand{\DeclareMathActive}[2]{%
  \expandafter\edef\csname keep@#1@code\endcsname{\mathchar\the\mathcode`#1 }
  \begingroup\lccode`~=`#1\relax
  \lowercase{\endgroup\def~}{#2}%
  \AtBeginDocument{\mathcode`#1="8000 }%
}

\newcommand{\std}[1]{\csname keep@#1@code\endcsname}
\patchcmd{\newmcodes@}{\mathcode`\-\relax}{\std@minuscode\relax}{}{%
  \PackageWarning{preamble_v3}{Failed to patch \string\newmcodes@ for minus mathcode handling}%
}
\AtBeginDocument{\edef\std@minuscode{\the\mathcode`-}}

\DeclareMathSymbol{\origcolon}{\mathpunct}{operators}{"3A}%
\DeclareRobustCommand{\coloneqq}{\mathrel{\mathop{\origcolon}\mkern-6mu=}}
\DeclareRobustCommand{\coloneq}{\mathrel{\mathop{\origcolon}\mkern-6mu=}}

\DeclareRobustCommand{\eqcolon}{=\mkern-1mu\mathrel{\mathop{\origcolon}}}

\def\@mathcoloneq={\coloneq}

\DeclareMathActive{:}{\@ifnextchar={\@mathcoloneq}{\nobreak\mskip1mu\mathpunct{}\nonscript \mkern-\thinmuskip{\std{:}\mskip6muplus1mu\relax}}}

\makeatother

\expandafter\let\expandafter\originall\csname\encodingdefault\string\l\endcsname
\DeclareRobustCommand*\l
    {\ifmmode\ell\else\expandafter\originall\fi}

\newcommand{\iprod}{\mathbin{\lrcorner}}
\newcommand{\acts}{\raisebox{-0.155em}{\scalebox{0.88}{ \rotatebox[origin=lb]{90}{$\curvearrowleft$}}}\,}

\newcommand{\ideal}{\unlhd}
\usepackage{accents}
\DeclareFontFamily{U}{matha}{\hyphenchar\font45}
\DeclareFontShape{U}{matha}{m}{n}{
    <5> <6> <7> <8> <9> <10> gen * matha
    <10.95> matha10 <12> <14.4> <17.28> <20.74> <24.88> matha12
}{}
\DeclareSymbolFont{matha}{U}{matha}{m}{n}
\DeclareFontSubstitution{U}{matha}{m}{n}
\DeclareMathSymbol{\tinsubset}{3}{matha}{"80}
\let\oldbigwedge\bigwedge
\renewcommand{\bigwedge}{\textstyle\oldbigwedge}

\newcommand{\qtext}[1]{\quad\text{#1}\quad}

\newcommand{\Ric}{\mathrm{Ric}}
\newcommand{\Rm}{\mathrm{Rm}}
\newcommand{\Bl}{\mathrm{Bl}} % Blow-up

\newcommand{\GL}{\mathrm{GL}}

\newcommand{\PGL}{\mathrm{PGL}}

\DeclareMathOperator{\Aut}{Aut}

\DeclareMathOperator{\Frac}{Frac}
\DeclareMathOperator{\gr}{gr}

\DeclareMathOperator{\Lie}{Lie}

\DeclareMathOperator{\Proj}{Proj}

\DeclareMathOperator{\Spec}{Spec}

\DeclareMathOperator{\supp}{supp}
\DeclareMathOperator{\Stab}{Stab}

\renewcommand{\Re}{\operatorname{Re}}

\expandafter\let\expandafter\originald\csname\encodingdefault\string\d\endcsname
\DeclareRobustCommand*\d
    {\ifmmode\mathop{}\!d\else\expandafter\originald\fi}

\newcommand{\p}{\mathrlap{~~\text{.}}}
\expandafter\let\expandafter\originalc\csname\encodingdefault\string\c\endcsname
\DeclareRobustCommand*\c
    {\ifmmode\mathrlap{~~\text{,}}\else\expandafter\originalc\fi}

\newcommand{\demph}{\textbf}

\makeatletter
\newcommand{\ovl}[1]{%
  \overline{\raisebox{0pt}[\dimexpr\height+0.6pt\relax]{\m@th$#1$}}%
}
\makeatother

\makeatletter
\let\save@mathaccent\mathaccent
\newcommand*\if@single[3]{%
    \setbox0\hbox{${\mathaccent"0362{#1}}^H$}%
    \setbox2\hbox{${\mathaccent"0362{\kern0pt#1}}^H$}%
    \ifdim\ht0=\ht2 #3\else #2\fi
}
\newcommand*\rel@kern[1]{\kern#1\dimexpr\macc@kerna}
\newcommand*\widebar[1]{\@ifnextchar^{{\wide@bar{#1}{0}}}{\wide@bar{#1}{1}}}
\newcommand*\wide@bar[2]{\if@single{#1}{\wide@bar@{#1}{#2}{1}}{\wide@bar@{#1}{#2}{2}}}
\newcommand*\wide@bar@[3]{%
    \begingroup
    \def\mathaccent##1##2{%
        \let\mathaccent\save@mathaccent
        \if#32 \let\macc@nucleus\first@char \fi
        \setbox\z@\hbox{$\macc@style{\macc@nucleus}_{}$}%
        \setbox\tw@\hbox{$\macc@style{\macc@nucleus}{}_{}$}%
        \dimen@\wd\tw@
        \advance\dimen@-\wd\z@
        \divide\dimen@ 3
        \@tempdima\wd\tw@
        \advance\@tempdima-\scriptspace
        \divide\@tempdima 10
        \advance\dimen@-\@tempdima
        \ifdim\dimen@>\z@ \dimen@0pt\fi
        \rel@kern{0.6}\kern-\dimen@
        \if#31
        \overline{\rel@kern{-0.6}\kern\dimen@\macc@nucleus\rel@kern{0.4}\kern\dimen@}%
        \advance\dimen@0.4\dimexpr\macc@kerna
        \let\final@kern#2%
        \ifdim\dimen@<\z@ \let\final@kern1\fi
        \if\final@kern1 \kern-\dimen@\fi
        \else
        \overline{\rel@kern{-0.6}\kern\dimen@#1}%
        \fi
    }%
    \macc@depth\@ne
    \let\math@bgroup\@empty \let\math@egroup\macc@set@skewchar
    \mathsurround\z@ \frozen@everymath{\mathgroup\macc@group\relax}%
    \macc@set@skewchar\relax
    \let\mathaccentV\macc@nested@a
    \if#31
    \macc@nested@a\relax111{#1}%
    \else
    \def\gobble@till@marker##1\endmarker{}%
    \futurelet\first@char\gobble@till@marker#1\endmarker
    \ifcat\noexpand\first@char A\else
    \def\first@char{}%
    \fi
    \macc@nested@a\relax111{\first@char}%
    \fi
    \endgroup
}
\makeatother
\def\semicolon{;}
\def\applytolist#1{
    \expandafter\def\csname multi#1\endcsname##1{
        \def\multiack{##1}\ifx\multiack\semicolon
            \def\next{\relax}
        \else
            \csname #1\endcsname{##1}
            \def\next{\csname multi#1\endcsname}
        \fi
        \next}
    \csname multi#1\endcsname}

\def\calchar#1{\expandafter\newcommand\csname c#1\endcsname{{\mathcal #1}}}
\applytolist{calchar}QWERTYUIOPLKJHGFDSAZXCVBNM;
\def\bbchar#1{\expandafter\newcommand\csname bb#1\endcsname{{\mathbb #1}}}
\applytolist{bbchar}QWERTYUIOPLKJHGFDSAZXCVBNM;
\def\bfchar#1{\expandafter\newcommand\csname bf#1\endcsname{{\mathbf #1}}}
\applytolist{bfchar}qwertyuiopasdfghjklzxcvbnmQWERTYUIOPLKJHGFDSAZXCVBNM;
\def\sfc#1{\expandafter\newcommand\csname s#1\endcsname{{\sf #1}}}
\applytolist{sfc}QWERTYUIOPLKJHGFDSAZXCVBNM;
\def\fchar#1{\expandafter\newcommand\csname f#1\endcsname{{\mathfrak #1}}}
\applytolist{fchar}qwertyuopasdfghjkzxcvbnmQWERTYUIOPLKJHGFDSAZXCVBNM;
\def\scchar#1{\expandafter\newcommand\csname sc#1\endcsname{{\mathscr #1}}}
\applytolist{scchar}QWERTYUIOPLKJHGFDSAZXCVBNM;